\documentclass[11pt,reqno,a4]{amsart}

\usepackage{a4}
\usepackage{graphicx}
\usepackage{rotating}
\usepackage{amsmath,amssymb,color,bbm,ifthen,enumerate}
\usepackage{natbib}
\usepackage{hyperref}
\usepackage{fancyvrb}
\usepackage[latin1]{inputenc}

\def\BIBDIR{.}
\def\FIGDIR{.}

  \advance\oddsidemargin -1.0cm
  \advance\evensidemargin -1.0cm

\newenvironment{enum_J}
  {%
  \setlength{\leftmargini}{4em}\begin{enumerate}}
  {\end{enumerate}}

\newenvironment{enum_W}
  {%
  \setlength{\leftmargini}{4em}\begin{enumerate}}
  {\end{enumerate}}

  {\end{enumerate}}

  {\end{enumerate}}

\def\maxscale{J}
\def\nmaxscale{\maxscale_{n}}
\def\lowscale{L}
\def\nlowscale{\lowscale_{n}}
\def\upscale{U}
\def\nupscale{\upscale_{n}}

\def\rme{\mathrm{e}}

\newcommand\jmean[1]{\langle #1 \rangle}
\newcommand\cardinal[1]{|#1|}
\def\dj{u}

\def\rmd{\mathrm{d}}
\def\rmi{\mathrm{i}}
\def\ltwo{\mathrm{L}^2}

\def\cl{\stackrel{\mathcal{L}}{\longrightarrow}}

\def\hd{\hat{d}} 
\def\hs{\hat{\sigma}} 
\def\td{d_{0}} 

\def\dwt{W}
\def\dft{D}
\def\bdwt{\mathbf{\dwt}}
\def\indexset{\mathcal{I}}
\def\LWF{LWF}

\def\L{\mathrm{T}}

\def\be{\mathbf{e}}

\def\shift{\mathbf{B}}
\def\shift{\mathbf{B}}

\def\diffop{\mathbf{\Delta}}

\def\1{\mathbbm{1}}

\def\bx{\mathbf{x}}

\def\Zset{\mathbb{Z}}

\def\Rset{\mathbb{R}} 
\def\PP{\mathbb{P}} 
\def\PE{\mathbb{E}} 
\newcommand{\CPE}[3][]
{\ifthenelse{\equal{#1}{}}{\operatorname{E}\left[\left. #2 \, \right| #3 \right]}{\operatorname{E}_{#1}\left[\left. #2 \, \right | #3 \right]}}
\def\PVar{\mathrm{Var}}
\def\PCov{\mathrm{Cov}}

\def\calN{\mathcal{N}}

\def\calH{\mathcal{H}}

\def\ie{\textit{i.e.} }
\def\eg{\textit{e.g.} }
\def\argmin{\mathop{\mathrm{Argmin}}}

\newcommand{\eqdef}{\ensuremath{\stackrel{\mathrm{def}}{=}}}
\newcommand{\eqsp}{\;}

\def\contrastsymb{\mathrm{L}}
\def\contrastL{\hat{\contrastsymb}}
\def\contrast{\tilde{\contrastsymb}}

\newcommand{\AVvar}[3][]
{
\ifthenelse{\equal{#1}{}}{\mathbf{V}_{#3}(#2)}{\mathbf{V}_{#3}(#2,#1)}}
\newcommand{\AVvarJoint}[3][]
{
\ifthenelse{\equal{#1}{}}{\mathbf{W}_{#3}(#2)}{\mathbf{W}_{#3}(#2,#1)}}

\newcommand{\AsympVarWWE}[2][]
{\rho^2(#2)}

\newcommand{\AsympVarRWEopt}[2][]
{\rho^{2}_{\mathrm{\tiny opt}}(#2)}

\newcommand{\regressopt}[3][]
{\regressweights{#3}^{\mathrm{\tiny opt}}(#2)}
\newcommand{\regressAV}[3][]
{\regressweights{#3}^{\mathrm{\tiny AV}}(#2)}

\newcommand{\AVvarInv}[3][]
{\ifthenelse{\equal{#1}{}}{\mathbf{V}^{-1}_{#3}(#2)}{\mathbf{V}^{-1}_{#3}(#2,#1)}}
\newcommand{\sigmaasymp}[2][]
{
\ifthenelse{\equal{#1}{}}{\sigma(#2)}{\sigma(#2,#1)}}

\def\vjsymb{\sigma}
\newcommand{\vj}[4][]{%
\ifthenelse{\equal{#1}{}}{\vjsymb^{#4}_{#2}(#3)}{\vjsymb^{#4}_{#2}(#3,#1)}}
\newcommand{\stdj}[3][]{%
\ifthenelse{\equal{#1}{}}{\vjsymb_{#2}(#3)}{\vjsymb_{#2}(#3,#1)}}

\newcommand{\hvj}[3][]{%
\ifthenelse{\equal{#1}{}}{\hat{\vjsymb}^2_{#2}}{\hat{\vjsymb}^2_{#2}(#1)}}
\newcommand{\Kvar}[2][]{
\mathrm{K}(#2)}

\def\densletter{\mathbf{D}}
\newcommand{\bdens}[4][]{%
\densletter_{#2}({#3};#4)}

\def\densasympletter{D}

\newcommand{\bdensasymp}[4][]{%
\mathbf{\densasympletter}_{\infty,#2}({#3};#4)}

\def\intbdensletter{I}
\newcommand{\intbdens}[3][]{
\mathrm{\intbdensletter}_{#2}(#3)}

\newcommand\regressweight[1]{w_{#1}}
\newcommand\regressweights[1]{\mathbf{w}_{#1}}

\def\LRW{\ensuremath{\mathrm{LRW}}}
\def\LWW{\ensuremath{\mathrm{LWW}}}
\newcommand\nj[1]{{n}_{#1}}
\def\K{\sigma^2}

\newcommand\bandwidth[1]{m_{#1}}
\def\pooling{\ensuremath{p}}
\def\taper{\ensuremath{\tau}}
\def\difforder{\ensuremath{\delta}}
\newtheorem{theorem}{Theorem}

\newtheorem{definition}{Definition}
\newtheorem{proposition}[theorem]{Proposition}
\newtheorem{corollary}[theorem]{Corollary}
\theoremstyle{remark}
\newtheorem{remark}{Remark}

\newcommand{\iid}{i.i.d.}
\def\GPH{\ensuremath{\mathrm{GPH}}}
\def\LWF{\ensuremath{\mathrm{LWF}}}

\def\d{\mathrm{d}}

\newcommand{\note}[2]{}

\renewcommand{\hat}{\widehat}
\renewcommand{\cl}{\stackrel{\mathcal{D}}{\longrightarrow}}

\title[Wavelet estimators of Long-Memory]{Estimators of Long-Memory: Fourier versus Wavelets}

\author{Gilles Faÿ}
\address{Laboratoire Paul-Painlevé, Université Lille-1, 59655 Villeneuve-d'Ascq
  Cedex, France.}
\curraddr{Laboratoire APC, Université Paris-7, Bâtiment Condorcet, 10, rue Alice Domon et
  Léonie Duquet, 75205 Paris Cedex 13, France.}
\email{gilles.fay@univ-lille1.fr}
\author{Eric Moulines}
\author{Fran\c{c}ois Roueff}
\address{LTCI (CNRS, TELECOM ParisTech) , 46, rue Barrault, 75634 Paris Cédex 13, France.}
\email{moulines@tsi.enst.fr}
\email{roueff@tsi.enst.fr}
\author{Murad S. Taqqu}
\address{Department of Mathematics and Statistics, Boston University Boston, MA 02215, USA.}
\email{murad@math.bu.edu}



\thanks{Murad~S.~Taqqu would like to thank l'\'Ecole Normale Sup\'erieure des T\'elecom\-munications in Paris
for their hospitality.  This research was partially supported by the NSF Grants DMS--0505747 and DMS--0706786 at Boston University.}

\date{January 28, 2008}

\begin{document}

\maketitle


\begin{abstract}
There have been a number of papers written on semi-parametric estimation methods of the long-memory exponent of a time series, some applied, others theoretical.
Some using Fourier methods, others using a wavelet-based technique.
In this paper, we compare the Fourier and wavelet approaches to the local regression method and to the local Whittle method.
We  provide an overview of these methods, describe what has been done, indicate the available results and the conditions under which they hold.
We discuss their relative strengths and weaknesses both from a practical and a theoretical perspective. We also include
a simulation-based comparison. The software written to support this work is available on demand and we illustrate its use at
the end of the paper.
\end{abstract}

\newpage

 \tableofcontents

 \newpage

\section{Introduction}
\label{sec:introduction}
We study here finite variance stochastic processes $\{X_k \}_{k \geq 1}$, whose spectral density $f(\lambda)$, $\lambda \in (-\pi,\pi)$
behaves like a power function at low frequencies, that is as $|\lambda|^{-2d}$ as the frequency $\lambda \to 0+$. The case
$d > 0$ corresponds to \emph{long-memory}, $d=0$ to \emph{short-memory} and $d < 0$ is often referred to as \emph{negative dependence}.
For $X_k, k \in \Zset$ to be stationary it is necessary that $\int_{-\pi}^{\pi} f(\lambda) d \lambda < \infty$ and hence
that $d < 1/2$. We relax these restrictions in a number of ways. We shall allow the process to be non-stationary, requiring only that
it becomes stationary after it is differenced a number of times. We also suppose that the spectral density (of the differenced process)
behaves not merely like $|\lambda|^{-2d}$ but as $|\lambda|^{-2d} f^\star(\lambda)$, where $f^\star$ is regarded as a short-range
density function.

Our goal is to estimate $d$ in the presence of $f^\star$. We shall not assume that the nuisance function $f^\star$ is known, nor that it
is characterized by a finite number of unknown parameters, but merely that $f^\star(\lambda)$ is "smooth" in the neighborhood of $\lambda=0$, so that
if one focuses only on frequencies $\lambda$ that are sufficiently low, then the spectral density $f(\lambda)$ behaves essentially
like $|\lambda|^{-2d}$. What frequency cut-off should one choose will clearly become an important issue.

The estimation framework is \emph{semi-parametric}: we must estimate the unknown parameter $d$ while viewing the presence of $f^\star$ as a nuisance,
albeit one which complicates matters. The estimation method will also be \emph{local}, in that, it is necessary to focus only on frequencies $\lambda$
that are close enough to the origin, where the influence of $f^\star(\lambda)$ can be neglected.

In this paper we provide an overview and comparison of four semi-parametric estimation methods of the parameter $d$ which have all proven
to be very effective. Two of them are Fourier-based, the other two are based on wavelets. The methods are:
\begin{itemize}
    \item Geweke-Porter Hudak (GPH): Regression / Fourier,
    \item Local Whittle Fourier (\LWF): Whittle / Fourier,
    \item Local Regression Wavelets (\LRW): Regression / Wavelets,
    \item Local Whittle Wavelets (\LWW):  Whittle / Wavelets.
\end{itemize}
The Fourier methods are older and better known. They have essentially been developed by Peter Robinson in a
number of fundamental papers \cite{robinson:1995:GPH}, \cite{robinson:1995:GSE}.
If we ignore for the moment the presence of the nuisance function $f^\star$, then one has $f(\lambda)= |\lambda|^{-2d}$,
that is $\log f(\lambda) \approx -2 d \log |\lambda|$. Therefore, $d$ can be estimated by linear regression on the
periodogram. This is the Fourier-based regression method considered in \cite{geweke:porter-hudak:1983} in a parametric setting.
The semi-parametric setting was suggested by \cite{kunsch:1987} and developed by \cite{robinson:1995:GPH}.
The Fourier-based Whittle method is a pseudo-maximum likelihood method developed by \cite{fox:taqqu:1986} in a parametric setting
and extended in a semi-parametric setting by \cite{robinson:1995:GSE}.

The papers of \cite{moulines:roueff:taqqu:2007:jtsa}, \cite{moulines:roueff:taqqu:2007:fractals}, \cite{moulines:roueff:taqqu:2007:aos} and
\cite{roueff:taqqu:2007} recast the preceding Fourier-based methods in a wavelet setting. Wavelets have a number of advantages. They allow differencing
implicitly and therefore they can be used without problems when $d > 1/2$. They also automatically discount polynomial trends.
The local wavelet-based regression method was first developed by \cite{abry:veitch:1998} under the simplifying assumption that
the wavelet coefficients are uncorrelated; see also \cite{veitch:abry:1999} and the review articles \cite{abry:flandrin:taqqu:veitch:2000} and
\cite{abry:flandrin:taqqu:veitch:2003-livre}. In addition, see \cite{veitch:abry:taqqu:2003} for the automatic selection of
the cut-off frequency point and \cite{veitch:taqqu:abry:2000m} for the choice of the "scale function".
\cite{bardet:2000T} and \cite{bardet:2002} provides asymptotic result for the \LRW\ estimator in a parametric context.
\cite{bardet:lang:moulines:soulier:2000} is a first attempt to analyze the behavior of \LRW\ in a semi-parametric context by
assuming continuous-time observations. The Local Whittle wavelet method is developed in \cite{moulines:roueff:taqqu:2007:aos}.

The paper is structured as follow. In Section \ref{sec:definition-an-md}, we formalise our assumptions on $\{X_k\}$ by
defining an M($d$) process, that is, a process with memory parameter $d$. The standard ARIMA and fractional Gaussian noise
examples are introduced in Section \ref{sec:examples}. The wavelet-based semi-parametric estimators are defined
in Section \ref{sec:wavelet-estimators} and the Fourier estimators in Section \ref{sec:fourier-estimators}. The semi-parametric
setting is discussed in Section \ref{sec:semi-parametric-setting}. The asymptotic properties of the wavelet and Fourier
semi-parametric estimators are described in Sections \ref{sec:w-results} and \ref{sec:f-results}, respectively. Their properties
are discussed further in Section \ref{sec:discussion}. Section \ref{sec:monte-carlo} contains the Monte-Carlo study which compares
the effectiveness of the four methods. In Section~\ref{sec:software}, we illustrate the use of the software written in support of this work.
This software may be obtained from the authors. Section~\ref{sec:conclusion} contains concluding remarks.


\section{Definition of an M($d$) process}
\label{sec:definition-an-md}
Let  $X \eqdef \{X_k \}_{k\in\Zset}$ be a real-valued process, not necessarily stationary. Its first order difference is
$$
[\diffop X]_n \eqdef X_n - X_{n-1},\quad n\in\Zset
$$
Its $K$-th order difference  $\diffop^K X$ is defined recursively.
We suppose that the process $X$ has \emph{memory parameter} $d$, $d \in \Rset$, in short, is an M($d$)
process. We shall first define this notion for a stationary process $X$, where $d<1/2$, and then provide a
general definition for $d\in\Rset$.

Let $f^\ast$ be a non-negative even function continuous and positive at the origin.
A stationary process $X$ is said to have memory parameter $d$,
$-\infty<d<1/2$,  and short-range density function $f^\ast$, if its spectral density is given by
\begin{equation}\label{eq:spdeltaStat}
f_{X}(\lambda) \eqdef |1-\rme^{-\rmi\lambda}|^{-2d}\,f^\ast(\lambda),\quad  \lambda\in(-\pi,\pi),
\end{equation}
To allow  $d > 1/2$, we consider non--stationary processes $X$ and
extend the preceding definition, valid for stationary processes, in the following way.
\begin{definition}\label{def:-an-md1}
We say that $X$ has \emph{memory parameter} $d$, $d\in\Rset$ (in short, an M($d$) process), and \emph{short-range
density function} $f^\ast$, if $f^\ast$ is continuous and positive at the origin and, for any integer $K>d-1/2$, its $K$-th
order difference $\diffop^KX$ is stationary with spectral density function
\begin{equation}\label{eq:spdelta}
f_{\diffop^K X}(\lambda) = |1-\rme^{-\rmi\lambda}|^{2(K-d)}\,f^\ast(\lambda),\quad  \lambda\in(-\pi,\pi) \; .
\end{equation}
\end{definition}
Observe that $f_{\diffop^K X}(\lambda)$ in~(\ref{eq:spdelta}) is integrable since $-(K-d)<1/2$. Observe also that if the
process $X$ is as in Definition~\ref{def:-an-md1}, then while
$\diffop^KX$ is stationary, the process $X$ itself is stationary only when $d<1/2$. Nevertheless,
one can associate to $X$ the \emph{generalized spectral density function}
\begin{equation}
\label{eq:SpectralDensity:FractionalProcess}
f_X(\lambda) = |1 - \rme^{-\rmi \lambda}|^{-2d}  f^\ast(\lambda)
\end{equation}
\begin{remark}
This definition of M($d$) processes was proposed by \cite{hurvich:ray:1995}.
It has the advantage that $\diffop^K X$ is stationary, but
it introduces a discontinuity  at the fractional points $d=1/2, 3/2, \dots$ since $f_{\diffop^K X}$ is quite
different at these values of $d$.
In empirical work, there are typically no inherent restrictions on the value of the memory parameter $d$,
and this may cause a problems if the degree of integer differencing required to achieve stationarity must be guessed in advance.
An alternative definition of M($d$) process has been introduced by \cite{robinson:1994:NS} and later used by some authors (see \cite{tanaka:1999}, \cite{shimotsu:phillips:2005}, \cite{shimotsu:phillips:2006}).

\end{remark}
The memory parameter $d$ plays a central role in the definition of M($d$) processes because it characterizes
the behavior of the generalized spectral density $f_X(\lambda)$ at low frequencies.
Indeed, assuming that $f^\ast$ is continuous at zero, then~(\ref{eq:SpectralDensity:FractionalProcess}) implies
$ f_X(\lambda)\sim|\lambda|^{-2d}f^\ast(0)$ as $\lambda\to 0$.  Allowing $d$ to take non integer values produces
a fundamental change in the correlation structure of a fractional process, as
compared to the correlation structure of a standard time-series model, such as an ARMA($p,q$) process.

The study of M($d$) processes has recently attracted  attention amongst theorists and empirical researchers. In
applied econometric work, M($d$) processes with $d > 0$ provide sensible  models for certain macroeconomic time series (inflation, interest rates, ...)
as well as certain financial time series (volatility of financial asset returns, forward exchange
market premia,...). M($d$) models encompass both stationary and nonstationary
processes depending on the value of the memory parameter and include both
short-memory series M($0$) and unit-root M($1$) processes as special cases when the memory parameter takes on the
values zero and unity.

\section{Examples}
\label{sec:examples}

Stationarity of the increments is commonly assumed in time-series analysis. In ARIMA models, for example,~(\ref{eq:spdelta})
holds with $d=K$ integer and with $f^\ast$ equal to the spectral density
of an autoregressive moving average short-memory process.
If $d\in\Rset$ and $f^\ast \equiv \sigma^2$  in~\eqref{eq:SpectralDensity:FractionalProcess},
one gets the so-called fractionally integrated white noise process, ARFIMA(0,$d$,0). The choice $d\in\Rset$ and
\begin{equation}
\label{eq:smooth:part:ARMA}
f_{\mathrm{ARMA}}^\ast(\lambda)= \sigma^2 \frac{\left| 1 - \sum_{k=1}^q \theta_k \rme^{-\rmi \lambda k}\right|^2}
{\left|1 - \sum_{k=1}^p \phi_k \rme^{- \rmi \lambda k} \right|^2},\quad \lambda\in(-\pi,\pi)\;,
\end{equation}
with $1 - \sum_{k=1}^p \phi_k z^k \ne 0$ for $|z|=1$ and  $1 - \sum_{k=1}^p \theta_k\ne0$ (so that $f_{\mathrm{ARMA}}^\ast(0)\ne0$)
leads to the class of ARFIMA($p,d,q$) processes.

Another example is $\{ B_H(k) \}_{k \in \Zset}$,  a discrete-time version of fractional Brownian
motion (FBM) $\{ B_H(t), t \in \Rset \}$ with Hurst index $H \in (0,1)$. The latter is a centered Gaussian process with
covariance
\[
R_H(t,s) \eqdef  \PE[ B_H(t) B_H(s)] = \frac{1}{2} \left\{ |t|^{2H} + |s|^{2H} - |t-s|^{2H} \right\} \eqsp .
\]
The process $\{ B_H(k) \}_{k \in \Zset}$ is increment stationary ($K=1$) and its generalized spectral density is given up to
a multiplicative constant (see~\cite{samorodnitsky:taqqu:1994}) by
$$
f_{\mathrm{FBM}}(\lambda) \eqdef \sum_{k=-\infty}^\infty\left|\lambda+2k\pi\right|^{-2H-1},\quad \lambda\in(-\pi,\pi)\;.
$$
We can express it in the form~(\ref{eq:SpectralDensity:FractionalProcess}),
\begin{equation}
\label{eq:GenDens:FBM}
f_{\mathrm{FBM}}(\lambda)=|1-\rme^{-\rmi\lambda}|^{-2d}\,f_{\mathrm{FBM}}^\ast(\lambda) \; ,
\end{equation}
by setting $d= H+1/2\in(1/2,3/2)$ and
\begin{equation}
\label{eq:smooth:part:FBM}
f_{\mathrm{FBM}}^\ast(\lambda)=\left|\frac{2 \sin(\lambda/2)}{\lambda}\right|^{2H+1}+
\left|2\sin(\lambda/2)\right|^{2H+1}\sum_{k\neq0}\left|\lambda+2k\pi\right|^{-2H-1} \; .
\end{equation}
Observe that $f_{\mathrm{FBM}}^\ast(0)=1$ and that it is bounded on $(-\pi,\pi)$.

The process $G_H=\diffop B_H$ is fractional Gaussian noise (FGN). It is a stationary Gaussian process with spectral density
proportional
to~(\ref{eq:GenDens:FBM}), but with $d=H-1/2\in(-1/2,1/2)$.

\section{Wavelet semi-parametric estimators of the memory parameter}\label{sec:wavelet-estimators}
In this section, we introduce the wavelet setting and, based on heuristical arguments, proposed possible
semi-parametric wavelet estimators. We start with a brief summary of the basic ideas.
A wavelet $\psi(t),\;t\in\Rset$ is a function with at least one vanishing moment, that is
,$\int_\Rset\psi(t)\rmd t=0$,  and which is \emph{low-pass}, in the sense that its Fourier transform $\hat{\psi} (\xi)$
decreases as  $\xi \to \infty$. We then define
the scaled and translated versions of $\psi$, namely $\psi_{j,k}(t)=2^{-j/2}\psi(2^{-j}t-k),\;j,k\in\Zset$.
The \emph{scale index} $j$ dilates $\psi$ so that large values of $j$ correspond to coarse scales (low frequencies),
while the \emph{position index} $k$ translates the function $\psi(2^{-j}t)$ to $\psi(2^{-j}t-k)$.
The corresponding \emph{wavelet coefficients} are then defined as $\dwt_{j,k}=\int_\Rset X(t)\psi_{j,k}(t)\rmd t$ and
are used to estimate $d$.
Because $\hat{\psi}$ is low-pass, $\hat{\psi}_{j,k}$ concentrates in the low frequency region
as $j\to\infty$ and $f_X(\lambda)$ ``scales'' at low frequencies since $|1-\rme^{\rmi\lambda}|^{-2d}\sim|\lambda|^{-2d}$ as $|\lambda|\to0$.
In the above definition of wavelet coefficients, we supposed, for simplicity, that the process $\{ X(t) \}_{t \in \Rset}$
is defined in continuous time and that the integral above is well-defined.
This definition can be adapted to discrete time series $\{X_k,\;k\in\Zset\}$ by using a \emph{scale function} and also to
finite samples $X_1,\dots,X_n$ by merely restricting the set of scale and translation indices of \emph{available wavelet coefficients}.
This is in sharp contrast to Fourier analysis, where the definition of discrete Fourier coefficients
at a given frequency changes as the sample length increases.
We now turn to a more formal presentation.

\subsection{The wavelet setting}
The wavelet setting involves a \emph{scale function} $\phi \in \ltwo(\Rset)$ and a \emph{wavelet} $\psi \in \ltwo(\Rset)$,
with associated Fourier transforms
\[
\hat{\phi}(\xi) \eqdef \int_{-\infty}^\infty \phi(t) \rme^{- \rmi \xi t}\,\rmd t \quad
\text{and}
\quad
\hat{\psi}(\xi) \eqdef \int_{-\infty}^\infty \psi(t) \rme^{- \rmi \xi t}\,\rmd t \eqsp.
\]
We assume the following:
\begin{enum_W}
\item\label{item:Wreg} $\phi$ and $\psi$ are compactly-supported, integrable, and $\hat{\phi}(0) = \int_{-\infty}^\infty \phi(t)\,\rmd t = 1$ and $\int_{-\infty}^\infty \psi^2(t)\,\rmd t = 1$.
\item\label{item:psiHat}
There exists $\alpha>1$ such that
$\sup_{\xi\in\Rset}|\hat{\psi}(\xi)|\,(1+|\xi|)^{\alpha} <\infty$.
\item\label{item:MVM} The function $\psi$ has  $M$ vanishing moments, \ie\ $ \int_{-\infty}^\infty t^m \psi(t) \,\rmd t=0$ for all $m=0,\dots,M-1$
\item\label{item:MIM} The function $ \sum_{k\in\Zset} k^m\phi(\cdot-k)$
is a polynomial of degree $m$ for all $m=0,\dots,M-1$.
\end{enum_W}\def\allWA{\ref{item:Wreg}-\ref{item:MIM}}
Condition~\ref{item:psiHat} ensures that the Fourier transform $\hat{\psi}$ decreases quickly
to zero. Daubechies wavelets have $\alpha>1$ (see Table~\ref{tab:alpha} below) except for Haar wavelet which is discontinuous
and for which $\alpha=1$.
Condition~\ref{item:MVM} it ensures that $\psi$ oscillates and that its
scalar product with continuous-time polynomials up to degree $M-1$ vanishes. It is equivalent to asserting that the first
$M-1$ derivative of $\hat{\psi}$ vanish at the origin and hence
\begin{equation}
\label{eq:MVM}
|\hat{\psi}(\lambda)|=O(|\lambda|^{M})\quad\text{as}\quad\lambda\to0.
\end{equation}
And, by \cite[Theorem~2.8.1, Page~90]{cohen:2003},~\ref{item:MIM} is equivalent to
\begin{equation}
\label{eq:MIM}
\sup_{k\neq0} |\hat{\phi}(\lambda+2k\pi)|=O(|\lambda|^{M})\quad\text{as}\quad\lambda\to0.
\end{equation}
As shown below, conditions~\ref{item:MIM}-\ref{item:MVM} imply that the wavelet transform of discrete-time polynomials of degree $M-1$ vanishes.


We now describe the computation of the wavelet coefficients.
Define the family $\{\psi_{j,k}, j \in \Zset, k \in \Zset\}$ of translated and dilated functions
\begin{equation}\label{eq:psiJK}
\psi_{j,k}(t)=2^{-j/2}\,\psi(2^{-j}t-k) ,\quad j\in\Zset,\, k\in\Zset \eqsp .
\end{equation}
Consider a real-valued sequence $\bx=\{x_k,\,k\in\Zset\}$.
We need to construct a continuous--time process from a discrete--time one. Using the scaling function $\phi$, we first associate to this sequence
the continuous-time functions
\begin{equation}\label{eq:bX}
\bx_n(t) \eqdef \sum_{k=1}^n x_k \,\phi(t-k) \quad\text{and}\quad \bx(t) \eqdef \sum_{k\in\Zset} x_k\, \phi(t-k), \quad
t\in\Rset \; .
\end{equation}
The function $\bx_n$ only requires the values of $x_1,\dots,x_n$ while the function $\bx$ requires the whole sequence $\{x_k,\;k\in\Zset\}$.
Without loss of generality we may suppose that the supports of the scaling function $\phi$ and
of the wavelet function $\psi$ are included in $[-\L,0]$ and  $[0,\L]$, respectively, for some
integer $\L\geq1$. This implies that $ \bx_n(t)=\bx(t)$ for all $t\in[0, n-\L+1]$ and that
the support of $\psi_{j,k}$ is included in the interval $[2^j k, 2^j(k+\L)]$.
The wavelet coefficient $\dwt^\bx_{j,k}$ at scale $j\geq0$ and location $k\in\Zset$ is defined as
\begin{equation}\label{eq:coeffN}
\dwt^\bx_{j,k} \eqdef \int_{-\infty}^\infty \bx(t) \psi_{j,k}(t)\,\rmd t
=\int_{-\infty}^\infty \bx_n(t) \psi_{j,k}(t)\,\rmd t,
 \quad  j \geq 0, k \in \Zset \;.
\end{equation}
The second equality holds when $[2^j k, 2^j(k+\L)]\subseteq [0, n-\L+1]$, that is, for all $(j,k)\in\indexset_n$, where
\begin{equation}
\label{eq:deltan}
\indexset_n \eqdef \{(j,k):\,j\geq0, 0\leq k  < \nj{j} \}\quad\text{with}\quad \nj{j}= \lfloor 2^{-j}(n-\L+1)-\L+1 \rfloor \eqsp.
\end{equation}
In other words $\indexset_n$ denotes the set of indices $(j,k)$ for which the wavelet coefficients $\dwt_{j,k}$
depend only on $x_1,\dots,x_n$. If the sample size $n$ increases, these wavelet coefficients remain unchanged and new ones
can be computed. Thus the definition of wavelet coefficients does not depend on the sample length,
in contrast to Fourier coefficients.
The wavelet coefficient $\dwt^\bx_{j,k}$ can be computed explicitly by using discrete convolution and downsampling, namely,
\begin{equation}\label{eq:down}
\dwt^{\bx}_{j,k}=\sum_{l\in\Zset} x_l\,h_{j,2^j k-l}=(h_{j,\cdot}\star \bx)_{2^jk}= (\downarrow^j [h_{j,\cdot}\star \bx])_{k}
,\quad j\geq0,k\in\Zset \;,
\end{equation}
where, for all $j\geq0$, the impulse response $h_{j,\cdot}$ is defined by
\begin{equation}\label{eq:FilterJ}
 h_{j,l} \eqdef 2^{-j/2} \int_{-\infty}^\infty \phi(t+l)\psi(2^{-j}t)\,\rmd t,\quad l\in\Zset \; ,
\end{equation}
where '$\star$' denotes the convolution  of discrete sequences and $\downarrow^j$ is the $j$--power downsampling operator defined, for any sequence $\{ c_k \}_{k \in \Zset}$, by $(\downarrow^j c)_k = c_{k 2^j}$.
Since $\phi$ and $\psi$ have compact support, the associated transfer function $H_j$ is a trigonometric polynomial,
\begin{equation}\label{eq:FilterJ-FIR}
H_j(\lambda)=\sum_{l\in\Zset} h_{j,l}\;\rme^{-\rmi\lambda l}
=\sum_{l=-\L(2^j+1)+1}^{-1} h_{j,l}\;\rme^{-\rmi\lambda l} \; .
\end{equation}

Under assumption \ref{item:MIM}, $t \mapsto \sum_{l \in \Zset} \phi(t+l) l^m$ is a polynomial of degree $m$ and \ref{item:MVM} therefore
implies that, for all $j\geq0$ and all $m=0,\dots,M-1$,
\begin{equation}
\label{eq:polyHj}
\sum_{l\in\Zset} h_{j,l}\,l^m= 2^{-j/2} \int_{-\infty}^\infty \psi(2^{-j}t) \sum_{l\in \Zset} \phi(t+l) l^m dt =0 \eqsp.
\end{equation}
Now consider $P_j(x)=\sum_{l \in \Zset} h_{j,l} \, x^l$ and observe that~(\ref{eq:polyHj}) implies $P_j(1)=0$, $P_j'(1)=0$,
..., $P_j^{(M-1)}(1)=0$, and hence $H_j(\lambda)=P_j(\rme^{-\rmi\lambda})$ factors as
\begin{equation}
\label{eq:Htilde}
 H_{j}(\lambda)= (1-\rme^{-\rmi \lambda})^{M} \, \tilde{H}_j(\lambda) \; ,
\end{equation}
where $\tilde{H}_j(\lambda)$ is also a trigonometric polynomial. The wavelet coefficient \eqref{eq:down} may therefore be computed as
\begin{equation}
\label{eq:op}
\dwt^{\bx}_{j,k}= (\downarrow^j [\tilde{h}_{j,\cdot}\star \diffop^M \bx])_{k}
\end{equation}
where $ \{ \tilde{h}_{j,l} \}_{l \in \Zset}$ are the coefficients of the trigonometric polynomial $\tilde{H}_j$ and
$\diffop^M \bx$ is the $M$-th order difference of the sequence $\bx$.
In other words, the use of a wavelet and a scaling function
satisfying \ref{item:MIM} and \ref{item:MVM} implicitly perform a
$M$-th order differentiation of the time-series. Therefore, we may
work with an M($d$) processes $X$ beyond the stationary regime
($d>1/2$) possibly contaminated by a polynomial trend of degree $K$
without specific preprocessing, provided that $d-M < 1/2$ and $M\geq
K+1$. It is perhaps less known that wavelets can be used with
non-invertible processes ($d \leq -1/2$) thanks to the decay property
of $\hat{\psi}$ at infinity assumed in \ref{item:psiHat}.
\note{f-EM}{we prefer to omit this : ``The choice
of the wavelet requires the prior knowledge of an upper bound and a
lower bound for $d$."}

\subsection{Choice of the wavelets}

In this paper, we do not assume that $\psi_{j,k}$ are orthonormal in $\ltwo(\Rset)$ nor
that they are associated to a multiresolution analysis (MRA). We may therefore use other convenient choices for $\phi$ and $\psi$
as long as~\allWA\ are satisfied. A simple choice is for instance, for some integer $M\geq2$,
\begin{equation}
 \label{eq:splinewavelet}
\phi(t) \eqdef \1_{[0,1]}^{\star M}(t)\quad\text{and}\quad\psi(t) \eqdef c_M\;\frac{\rmd^M}{\rmd t^M} \1_{[0,1]}^{\star 2M}(2t),
\end{equation}
where $\1_A$ is the indicator function of the set $A$ and for an integrable function $f$, $f^{\star M}$ denotes the $M$-th  self-convolution of
$f$,
$$
f^{\star M}=\underbrace{f\star\dots\star f}_{M\text{ times }},\text{ with } (f\star g) (t) = \int_{-\infty}^\infty f(t-u)\;g(u)\;\rmd u \;,
$$
and $c_M$ is a normalizing positive constant such that $\int_{-\infty}^\infty\psi^2(t)\;\rmd t=1$.

Scaling and wavelet functions associated to an MRA present two important features: 1) they give raise orthonormal
$\ltwo(\Rset)$ bases $\{\psi_{j,k}\}$; 2) a recursive algorithm, the so--called pyramidal algorithm,
is available for performing the convolution/downsampling operations at all scale $j$.
The complexity of this algorithm is $O(n)$ for a sample of length $n$, see \cite{mallat:1998}. In other words, in an MRA, the
computation~(\ref{eq:down}) can be made recursively as $j$ grows and it is not
necessary to explicitly compute the filters $h_{j,\cdot}$ defined in~(\ref{eq:FilterJ}).

Assumptions \allWA\ are  standard in the context of an MRA, see for instance~\cite{cohen:2003}.
Common wavelets are Daubechies wavelet and Coiflets (for which the scale function
also has vanishing moments). What matters in the asymptotic
theory of wavelet estimators presented below is the number of vanishing moments $M$ and the decay
exponent $\alpha$, which both determine the  frequency resolution of $\psi$. For standard wavelets, $M$ is always known and
\cite[Remark~2.7.1, Page~86]{cohen:2003} provides a sequence of lower bounds $(\alpha_k)$ tending to
$\alpha$ as $k\to\infty$. Daubechies wavelet are defined by their number $M$ of vanishing moments, for any $M\geq1$
(the case $M=1$ corresponds to the so called Haar wavelet). An analytic formula for their decay exponent $\alpha$
is available, see \cite[Eq~(7.1.23), Page~225 and the table on Page~226]{daubechies:1992} and
note that our $\alpha$ equals the $\alpha$ of \cite{daubechies:1992} plus 1.
A simpler lower bound  $\alpha\geq (1-\log_2(3)/2) M$ holds for Daubechies wavelets, see~\cite{daubechies:1992}. Although it
is not sharp, it shows that the number of vanishing moments $M$ and decay exponent
$\alpha$ of Daubechies's wavelets can be made arbitrarily large at the same time.
Table~\ref{tab:alpha} provides some values of $\alpha$ for Daubechies wavelets
and the lower bound $\alpha_{k}$ with $k=10$ for Coiflets  with given number of vanishing moments $M$ ranging from 2 to 10.

\begin{table}[ht]
\begin{tabular}{l|c|c|c|c|c|c|c|c|c} 
$M$&2 &3 &4 &5 &6 &7 &8 &9 &10 \\ 
\hline 
$\alpha$ (DB)&1.3390&1.6360&1.9125&2.1766&2.4322&2.6817&2.9265&3.1676&3.4057\\ 
\hline 
$\alpha_{10}$ (Coif.) & 1.6196&N.A. &1.9814&N.A. &2.5374&N.A. &3.0648&N.A. &3.5744\\ 
\end{tabular}

\caption{The decay exponent $\alpha$ or its lower bound $\alpha_{10}$ of $|\hat{\phi}(\xi)|$ (and hence of  $|\hat{\psi}(\xi)|$) with  $M$ vanishing moments. First line: $M$; second line: $\alpha$ for Daubechies wavelet; third line: the lower bound $\alpha_{10}$ for the Coiflet. N.A. stands for \emph{not available} (Coiflets are defined for $M$ even).}

\label{tab:alpha}
\end{table}

In view of Table~\ref{tab:alpha}, one can observe that the decays of Coiflets are slightly faster than the ones of Daubechies
for given $M$'s. On the other hand the Daubechies wavelets have shorter support, since it is of length $\L=2M$ while
it is of length  $\L=3M$ for Coiflets. The support length  impacts on the number of available wavelet coefficients:
given a sample size $n$, the greater the support length $\L$ the smaller the cardinality of the set $\indexset_n$ defined in~(\ref{eq:deltan}).



We should also mention the so-called Shannon wavelet $\psi_S$ whose Fourier transform $\hat{\psi}_S$ satisfies
\begin{equation}\label{eq:psihatShannon}
|\hat{\psi_S}(\xi)|^2= \begin{cases}
1 & \quad \text{for}\quad |\xi|\in[\pi,2\pi]\\
0 & \quad \text{otherwise}.
\end{cases}
\end{equation}
This wavelet satisfies~\ref{item:psiHat}--\ref{item:MIM} for arbitrary large $M$ and $\alpha$ but
does not have compact support, hence it does not satisfy~\ref{item:Wreg}. We may therefore not choose this wavelet in our analysis.
It is of interest, however, because it gives a rough idea of what happens when $\alpha$ and $M$ are large since one can always
construct a wavelet $\psi$ satisfying~\ref{item:Wreg}--\ref{item:MIM} which is arbitrarily close to the Shannon wavelet.

\subsection{The local regression wavelet (\LRW) estimator of $d$}\label{sec:reg}
For any integers $n$, $j_0$ and $j_1$, $j_0 \leq j_1$, the set of all available wavelet coefficients from $n$
observations $X_1,\dots,X_n$ having scale indices between $j_0$ and $j_1$ is
\begin{equation}\label{eq:InDef}
\indexset_n(j_0,j_1) \eqdef \left\{ (j,k):\,j_0 \leq j \leq j_1, \; 0\leq k < \nj{j} \right\} \;,
\end{equation}
where $\nj{j}$ is given in \eqref{eq:deltan}.
Consider two integers   $\lowscale<\upscale$ satisfying
\begin{equation}\label{eq:J0J1n}
0\leq \lowscale< \upscale \leq \maxscale_n \eqdef \max\{j:\,\nj{j}\geq1\} \eqsp.
\end{equation}
The index $\nmaxscale$ is the maximal available scale index for the sample size $n$;
$\lowscale$ and $\upscale$ will denote, respectively, the lower and upper scale indices used in the estimation.
For an M($d$) process, under regularity
conditions on the short-memory part $f^\ast$ and for appropriately chosen scale function and wavelet $\phi$ and $\psi$, it may be shown that as $j \to \infty$,
$\vj{j}{d,f^{\ast}}{2} \eqdef \PVar [ \dwt^X_{j,0}]  \asymp \sigma^2 2^{2 d j}$
 and the empirical variance
\begin{equation}\label{eq:hvj}
\hvj{j}{\nj{j}} \eqdef \nj{j}^{-1} \sum_{k=0}^{\nj{j}-1} \left( \dwt_{j,k}^X \right)^2 \eqsp,
\end{equation}
is a consistent sequence of estimator of $\vj{j}{d,f^{\ast}}{2}$ (see Proposition \ref{prop:approximation:spectral:wavelet}).
A popular semi-parametric estimator of the memory parameter $d$ is the
\emph{local regression wavelet} (\LRW) estimator of \cite{abry:veitch:1998}, defined as the least squares estimator in the
"linear regression
model"
\[
\log \left[ \hvj{j}{\nj{j}} \right] = \log \sigma^2 + d j \{ 2 \log(2) \} + u_j \eqsp,
\]
where $u_j = \log[ \hvj{j}{\nj{j}} / \sigma^2 2^{2 d j} ]$. This regression problem can be solved in closed form:
\begin{equation}
\label{eq:definition:estimator:regression}
\hd_n^{\LRW}(\lowscale,\upscale,\regressweights{}) \eqdef \sum_{j=\lowscale}^{\upscale} \regressweight{j-\lowscale} \log \left( \hvj{j}{n_j} \right) \eqsp ,
\end{equation}
(in short $\hd^{\LRW}$) where the vector $\regressweights{} \eqdef [\regressweight{0},\dots,\regressweight{\upscale-\lowscale}]^T$  of weights satisfies
\begin{equation}
\label{eq:propertyw}
\sum_{j=0}^{\upscale-\lowscale} \regressweight{j}  = 0 \quad \text{and} \quad 2 \log(2) \sum_{j=0}^{\upscale-\lowscale} j \regressweight{j}  = 1 \eqsp.
\end{equation}
For $\upscale-\lowscale=\ell\geq1$, one may choose, for example, $\regressweights{}$ corresponding to the weighted least-squares
regression vector, defined by
\begin{equation}
  \label{eq:generalFormweights}
\regressweights{} \eqdef D B (B^TDB)^{-1} \mathbf{b}
\end{equation}
where
\begin{equation}
  \label{eq:bAndB}
\mathbf{b}\eqdef \left[\begin{matrix}
0 \\
(2\log(2))^{-1}
\end{matrix}
\right],
\quad
B \eqdef \left[\begin{matrix}
1 & 1 & \dots & 1 \\
0 & 1 & \dots & \ell
\end{matrix}\right]^T
\end{equation}
is the design matrix and $D$ is an arbitrary positive definite matrix.
We will discuss the choice of the regression weights $\regressweights{}$ after stating Theorem~\ref{thm:LRW:CLT},
which provides the
asymptotic variance of  $\hd^{\LRW}(\lowscale,\upscale,\regressweights{})$.

\subsection{The local Whittle wavelet (\LWW) estimator of $d$}
\label{subsec:LWWestimator} Let $\{c_{j,k},\,(j,k)\in\indexset\}$ be
an array of centered independent Gaussian random variables with
variance $\PVar(c_{j,k})=\vjsymb_{j,k}^2$, where $\indexset $ is a
finite set.  The negative of its log-likelihood  is $ (1/2)
\sum_{(j,k) \in \indexset} \left\{ c_{j,k}^2/\vjsymb_{j,k}^2 +
\log(\vjsymb_{j,k}^2)\right\}$ up to a constant additive term. The
\emph{local Whittle  wavelet} (\LWW) estimator uses such a contrast process
to estimate the memory parameter $\td$ by choosing
$c_{j,k}=\dwt^X_{j,k}$ and
$\indexset=\indexset_n(\lowscale,\upscale)$ as defined in
\eqref{eq:InDef} for \emph{appropriately chosen} lower and upper
scale indices $\lowscale$ and $\upscale$, so that the corresponding
wavelet coefficients $\dwt^X_{j,k}$ are computed from
$X_1,\dots,X_n$. The scaling property $\vj{j}{d,f^{\ast}}{2} \asymp
\sigma^2 2^{2 d j}$ and weak dependence conditions of the wavelet
coefficients then suggest the following \emph{pseudo} negative
log-likelihood
\[
\contrastL_{\indexset} (\K,d)
=   \frac1{2 \K} \sum_{(j,k) \in \indexset}2^{-2 d j} (\dwt^X_{j,k})^2 + \frac{\cardinal{\indexset}}{2} \log(\K2^{2\jmean{\indexset}\, d}) \eqsp,
\]
where $ \cardinal{\indexset}$ denotes the cardinal of  $\indexset$ and
$\jmean{\indexset}$ is the average scale, $\jmean{\indexset}\eqdef \cardinal{\indexset}^{-1} \sum_{(j,k) \in \indexset} j$.
Define
$\hs^2_\indexset(d)\eqdef\argmin_{\K>0}\contrastL_{\indexset} (\K,d) =
\cardinal{\indexset}^{-1} \sum_{(j,k) \in \indexset} 2^{-2 d j} (\dwt^X_{j,k})^2 \eqsp \eqsp.
$
The pseudo maximum likelihood estimator of the memory parameter is then equal to the minimum of the negative profile log-likelihood,
\begin{equation}\label{eq:LWWdef}
\hd^{\LWW}(\lowscale,\upscale) \eqdef
\argmin_{d\in [\Delta_1,\Delta_2]} \contrastL_{\indexset_n(\lowscale,\upscale)} (\hs^2_\indexset(d),d)
\end{equation}
where $[\Delta_1,\Delta_2]$  is an interval of admissible values for $d$ and
\begin{equation}
\label{eq:TildeJdef}
\contrast_{\indexset} (d) \eqdef \log \left(\sum_{(j,k) \in \indexset}\,2^{2d(\jmean{\indexset}-j)}(\dwt^X_{j,k})^2 \right) \eqsp.
\end{equation}
This estimator has been proposed for analyzing  noisy data in~\cite{wornell:oppenheim:1992},
and was then considered by several authors, mostly in a parametric context, see \eg \cite{kaplan:kuo:1993} and
\cite{mccoy:walden:1996}.
If $\indexset$ contains at least two different scales then $\contrast_{\indexset} (d)\to\infty$ as
$d\to\pm\infty$, and thus $\hd$ is finite. 
This contrast is strictly convex, and the minimum is unique: it can be found using any one-dimensional convex
optimization procedure. In contrast to the the local regression wavelet estimator, the definition of LWWE does not relies on particular weights.
An important issue for both estimators is the choice of the scale indices $\lowscale$ and $\upscale$.
The asymptotic theory developed for these estimators in a semi-parametric context sheds some light
on the role played by these quantities, as will be explained in Section \ref{sec:w-results}.

\section{Fourier semi-parametric estimators of the memory parameter}\label{sec:fourier-estimators}
\subsection{The periodogram}
\label{sec:notations}
Given $n$ observations $X_1,\cdots,X_n$, the \emph{discrete Fourier transform} (DFT) and the \emph{periodogram} are respectively defined as
\begin{equation}
\label{eq:DX}
\dft^X(\lambda) \eqdef (2\pi n)^{-1/2}  \sum_{t=1}^n X_t \rme^{\rmi t \lambda} \eqsp, \quad I^X(\lambda) \eqdef |\dft^X(\lambda)|^2 \eqsp.
\end{equation}
These quantities are computed at the Fourier frequencies
\begin{equation}
\label{eq:Ff}
\lambda_j \eqdef 2 \pi j /n \eqsp, \quad \text{for $k \in \{1, \dots, \tilde{n}\}$, where $\tilde{n}= \lfloor (n-1)/2 \rfloor$} \eqsp.
\end{equation}
For stationary and invertible $M(d)$ processes (see for example \cite{lahiri:2003}), the DFT coefficients at
Fourier frequencies are known to be \emph{approximately} asymptotically independent outside a shrinking neighborhood of zero. Thus the Fourier transform performs
a \emph{whitening} of the data, and, as a consequence, Fourier methods have neat asymptotic statistical properties.
\subsubsection{Differencing and Tapering}
\label{sec:tapering}
To overcome the presence of polynomial or smooth trends (see \cite{hurvich:lang:soulier:2005}), or to estimate the
memory parameter of an $M(d)$ process beyond the stationary regime ($d > 1/2$), some adjustments are necessary.
For instance, it has been suggested to apply a data taper either to the time-series $X$  or to its ${\difforder}$-th order
difference $\diffop^{\difforder} X$.
A taper is a non-random weight function (with certain desired properties) that is multiplied to the time-series
(or its difference) prior to Fourier transformation.
Tapering was originally used in nonparametric spectral analysis of short memory $(d=0)$ time series in order to reduce bias due to
frequency domain leakage, where part of the spectrum "leaks" into adjacent frequencies. The leakage is due to the discontinuity caused by the finiteness of the sample and is
reduced by using tapers which smooth this discontinuity. But such a bias reduction inflates the variance as will be seen
later in a special case in the context of long memory semi-parametric estimation (see section \ref{sec:f-results}).

The idea of applying taper directly to the observations $X$ was proposed by
\cite{velasco:1999,velasco:1999:JE,velasco:robinson:2000}, who considered
several tapering schemes such as the cosine bell and the Zurbenko-Kolmogorov
tapers \citep{zurbenko:1979}. These tapers $h_{t}$, $t=1,\dots,n$
have the property of being orthogonal to polynomials up to a given order,
for a subset of Fourier frequencies,
\begin{equation}
\label{eq:orthogonality}
\sum_{t=1}^n  (1 + t + \dots + t^{{\difforder}-1}) h_{t} \rme^{\rmi t \lambda_j} = 0 \eqsp, \quad j  \in \mathcal{J}_{{\difforder},n} \subset \{1,\dots, \tilde{n}\} \eqsp,
\end{equation}
where $\lambda_j$ and $\tilde{n}$ are defined in \eqref{eq:Ff}.
A problem with this approach is that   the efficiency loss due to these
tapers may be quite substantial, because the set $\mathcal{J}_{{\difforder},n}$ can be fairly small when ${\difforder}$ is large.

In this contribution, we rather focus on the construction suggested in \cite{hurvich:ray:1995} and later developed in
\cite{hurvich:chen:2000}, which consists in differencing before tapering.  Differencing is a very widely used technique for
detrending and inducing stationarity. The ${\difforder}$-th order difference will convert the memory parameter of a $M(d)$ process to
$d-{\difforder}$, and will completely remove a polynomial trend of degree ${\difforder}-1$. To apply this technique, an upper bound to the memory
parameter (or to the degree of the polynomial trend) should be known in advance. But if only an upper bound is known, $\difforder$ may be chosen too large and consequently
the ${\difforder}$-th order difference may be non-invertible,
that is one may have $d-{\difforder} \leq -1/2$. This situation, referred to as \emph{over-differencing} which may cause  difficulties in spectral inference (see \cite{hurvich:ray:1995}).
As was suggested by these authors, the use of a data taper can alleviate the detrimental effect of overdifferencing.
A main drawback with this approach is that tapering inflates the variance of the estimator. To minimize this  effect, the tapers should be chosen carefully.

\cite{hurvich:chen:2000} have defined a family of data taper depending on a single parameter $\taper$, referred to as the
\emph{taper order}. Set $ h_{t} = 1 - \rme^{2 \rmi \pi t / n}$
and, for any integer $\taper \geq 0$,  define the tapered DFT of order $\taper$ of the sequence $\bx= \{x_k, k \in \Zset \}$
as follows
\begin{equation}
\label{eq:hurvichtaper}
\dft^{\bx}_{\taper}(\lambda) \eqdef (2\pi n a_\taper)^{-1/2} \sum_{t=1}^n h_{t}^\taper x_t \rme^{\rmi t \lambda},
\ \ I^{\bx}_{\taper}(\lambda) \eqdef | \dft^{\bx}_{\taper}(\lambda)|^2
\end{equation}
where the subscript $\taper$ denotes the taper order and $a_\taper\eqdef n^{-1} \sum_{t=1}^n | h_{t}|^{2\taper}$ is a normalization factor.
As shown in \cite[Lemma~0]{hurvich:chen:2000},
the decay of the discrete Fourier transform of the taper of order $\taper$ is given by
$$
\left| (2\pi n a_\taper)^{-1/2} \sum_{t=1}^n h_{t}^\taper \rme^{\rmi t \lambda}\right| \leq  C\frac{n}{(1+n|\lambda|)^\taper},\quad\lambda\in(-\pi,\pi)\;.
$$
This property means that higher-order tapers control the leakage more effectively.

Note that the Fourier transform of the taper
may be expressed as a finite sum of shifted Dirichlet kernels,
\begin{equation}
\label{eq:dft-hurvich-taper}
\sum_{t=1}^n h_t^\taper \rme^{\rmi t \lambda} =
\sum_{k=0}^\taper \left\{ \sum_{t=1}^n \rme^{\rmi t (\lambda + \lambda_k)} \right\}
\end{equation}
Since $\sum_{t=1}^n \rme^{\rmi t \lambda_k}= 0$, this relation implies that for $j \in \{1, \dots, \tilde{n}-\taper\}$, $\sum_{t=1}^n h_{t}^\taper \rme^{\rmi t  \lambda_j}= 0$
so that the tapered Fourier transform (evaluated at Fourier frequencies) is invariant to shift in the mean. This
shift-invariance is achieved without restricting attention
to a coarse grid of Fourier frequencies, as is necessary for the Zhurbenko-Kolmogorov
taper (\cite{velasco:1999}).

However, the construction of theoretically-justified memory estimators using the tapered periodogram  may require dropping some Fourier frequencies as will be
seen in the definition of the pooled periodogram in~(\ref{eq:pooled-periodogram}). This is related to the following
observation.
For $\taper=0$ (no taper), the DFT coefficients $\dft^Z_{\taper}(\lambda_j)$ of a white
noise $\{Z_t \}$ at Fourier frequencies $\lambda_k$, $\lambda_j$, $k\ne j \in \{1,\dots, \tilde{n} \}$ are uncorrelated.
This property is lost by tapering.  For $\taper\geq1$, the correlation $\PE[ \dft^Z_{\taper}(\lambda_j)
\overline{\dft^Z_{\taper}(\lambda_k)} ]$  is equal to $0$ if $|k-j| > \taper$, and
$(2 \pi a_\taper)^{-1}  (-1)^k \left(^{2\taper}_{\taper+k}\right)$ if $|k| \leq \taper$.

\subsubsection{Pooling}
\label{sec:pooling}
Let $I^{\bx}_{\taper}(\lambda)$ be the tapered periodogram introduced in~(\ref{eq:hurvichtaper}).
When considering non-linear transformations of the periodogram such as taking logarithm,
variance reduction can be obtained by pooling groups of finitely many, say $\pooling$, consecutive
$I^{\bx}_{\taper}(\lambda_j)$, which results in a \emph{pooled periodogram}
(see \cite{hannan:nicholls:1977} and \cite{robinson:1995:GPH}). To understand
why pooling may be helpful, recall that if $Z$ is a white Gaussian noise,
then the variance of the running mean of the log-periodogram $\pooling^{-1/2} \sum_{k=1}^{\pooling} \log
I^{Z}(\lambda_{j+k})$ is $\psi'(1)$ whereas the variance of the logarithm of the running mean $\pooling^{-1/2} \log\left(
  \sum_{k=1}^{\pooling}  I^{Z}(\lambda_{j+k}) \right)$ is $\pooling \psi'(\pooling)$, where $\psi(z) = \Gamma'(z)/\Gamma(z)$
is the digamma function (see for instance \cite{johnson:kotz:1970}). The quantity $\pooling \psi'(\pooling)$ decreases from
$\pi^2/6$ to $1$ as $\pooling$ goes from $1$ to $\infty$.
For $\pooling=3$,  $\pooling \psi'(\pooling)$ is 1.1848 and its value  changes slowly thereafter as $\pooling
\psi'(\pooling)= 1 + 1/(2 \pooling) + O( \pooling^{-2})$. Nonetheless, this shows that the
variance of the local average of the log-periodogram is larger than the variance of the logarithm of local average and
explains why typical values of $\pooling$ are $\pooling=3,4$.

As seen above for the tapered DFT coefficients of a white noise in a non-asymptotic context, in order to guarantee asymptotic
independence of the \emph{tapered} periodogram ordinates, if $\pooling$ successive values of the periodogram are pooled, then, at
the end of the block, $\taper$ DFT coefficients are dropped, where $\taper$ is the taper order. More
precisely, set $K(\pooling,\taper) = [(n-1)/2(\pooling+\taper)]$ and for $k \in \{
1, \cdots, K(\pooling,\taper) \}$, define the \emph{pooled periodogram} as follows
\begin{equation}
\label{eq:pooled-periodogram}
\bar{I}^{\bx}_{\pooling,\taper}(\tilde{\lambda}_k) \eqdef \sum_{j=(\pooling+\taper)(k-1) +
1}^{(\pooling+\taper)(k-1) + \pooling} I^{\bx}_{\taper}(\lambda_j).
\end{equation} 
where $\tilde{\lambda}_k \eqdef \pooling^{-1} \sum_{j=(\pooling+\taper)(k-1) +
  1}^{(\pooling+\taper)k} \lambda_j= \left( 2(\pooling + \taper) (k-1) +
  \pooling + \taper +1 \right)\pi/n$. The definition of the central frequency
$\tilde{\lambda}_k$ seems somehow arbitrary. Our choice is motivated by the fact
that the Chen and Hurvich's taper actually mixes together $\taper$
adjacent periodogram ordinates, so that $(\pooling + \taper)$ Fourier
frequencies are mixed in each pooled and tapered periodogram ordinate. Note
that the bias of the \GPH\ estimator defined below is very sensitive to the
definition of this central frequency.



\subsection{The Geweke--Porter--Hudak  (\GPH) estimator of $d$}
\label{sec:defgph}
Assume that the differencing order ${\difforder}$
induces stationarity, \ie\ $d < {\difforder} +1/2$, and the taper order $\taper$ is larger than ${\difforder}$. \note{f-EM}{slightly rephrased}
Then, for certain sequences $\{ \ell_n \}$ and $ \{ \bandwidth{n} \}$ which increase slowly with $n$, and
under  smoothness conditions on the short memory component $f^\ast$,
the ratios of the pooled periodogram of the ${\difforder}$-th order difference of $Y\eqdef \diffop^{\difforder} X$ divided by  its spectral density
$\bar{I}^{Y}_{\pooling,\taper}(\tilde{\lambda}_k) / f_{Y} ( \tilde{\lambda}_k)$, $\ell_n \leq k \leq \bandwidth{n}$,
can be regarded as approximately independent and identically distributed (i.i.d.) in a sense that can
be rigorously characterized; \cite{robinson:1995:GPH} and \cite{lahiri:2003}.

Based on this heuristics, a popular semiparametric estimate of $d$ is the log-periodogram estimate of \cite{geweke:porter-hudak:1983},
defined here [in the manner of \cite{robinson:1995:GPH}] as the least squares estimate in the \emph{linear regression model}
\begin{equation}
\log \left[ \bar{I}^{Y}_{\pooling,\taper}(\tilde{\lambda_k}) \right] = \log f^\ast(0) + (d-{\difforder}) g(\tilde{\lambda}_k) + u_k
\eqsp, \quad 1 \leq k \leq \bandwidth{} \eqsp,
\end{equation}
where $g(\lambda) \eqdef -2 \log |1 -\rme^{\rmi \lambda}|$ and $u_k$ is "approximately" equal to $\log \left[ \bar{I}^Y_{\pooling,\taper}(\tilde{\lambda}_k) / f_X(\tilde{\lambda}_k)  \right]$.
This regression equation can be solved in closed form~:
\begin{equation}
\label{eq:definition-GPH-estimator}
\hat{d}^{\GPH}(\bandwidth{})  = \sum_{k=1}^{\bandwidth{}} \frac{\left\{g(\tilde{\lambda}_k)- \bandwidth{}^{-1}\sum_{k=1}^{\bandwidth{}} g(\tilde{\lambda}_k)\right\}}
{\sum_{k=1}^{\bandwidth{}} \left\{g(\tilde{\lambda}_k)- \bandwidth{}^{-1}\sum_{k=1}^{\bandwidth{}} g(\tilde{\lambda}_k) \right\}^2}
 \log [\bar{I}^Y_{\pooling,\taper}(\tilde{\lambda}_k)] + {\difforder}\eqsp.
\end{equation}
To simplify the notations, we have made the dependence in the differencing, tapering and pooling orders  implicit.
The choice of these orders ${\difforder}$, $\taper$ and $\pooling$ will be discussed in Section~\ref{sec:f-results}.
This estimator has been introduced by \cite{geweke:porter-hudak:1983} and was later used in many empirical works.
\begin{remark}
\label{re:trimming}
In the definition of the \GPH\ estimator in \cite{robinson:1995:GPH}, the first $\ell_n$ DFT coefficients are eliminated. $\ell_n$
is referred to as the \emph{trimming number}. Trimming is sometimes required to eliminate deterministic trend \citep{hurvich:lang:soulier:2005}, or
to deal with non-Gaussian processes \citep{velasco:2000:NG}.
\end{remark}
\subsection{The local Whittle Fourier (\LWF) estimator of $d$}
\label{sec:defgse}
Since, as mentioned above, $I^Y_{\taper}(\lambda_k)/f_Y(\lambda_k)$ can be regarded as approximately \iid\ and $f_Y(\lambda)
\approx C |1- \rme^{\rmi \lambda}|^{-2(d-{\difforder})}$
in the neighborhood of zero,  using the same arguments as in Section \ref{subsec:LWWestimator},
we may approximate the negated likelihood as follows:
\begin{equation}
\label{eq:Whittle}
\hat{L}_{\taper,\bandwidth{}}(C,d) = \bandwidth{}^{-1} \sum_{k=1}^{\bandwidth{}} \left \{ \log (C |1-\rme^{\rmi \lambda_k}|^{-2(d-{\difforder})}) +
  \frac{I^Y_{\taper}(\lambda_k) }{C |1 - \rme^{\rmi \lambda_k}|^{-2(d-{\difforder})}} \right \} \eqsp.
\end{equation}
Note that pooling is here irrelevant because non non-linear transformation is involved. This estimator was originally proposed by \cite{kunsch:1987} and later studied in
\cite{robinson:1995:GSE}. After eliminating $C$ by maximizing the contrast \eqref{eq:Whittle}, we get $\tilde{L}^{\LWF}_{\taper,\bandwidth{}}$ the \emph{profile likelihood}, defined as
\begin{equation}
\label{eq:profile:likelihood}
\tilde{L}^{\LWF}_{\taper,\bandwidth{}}(d)
\eqdef \log \left( \bandwidth{}^{-1} \sum_{k=1}^{\bandwidth{}} I^Y_{\taper}(\lambda_k) |1 - \rme^{\rmi \lambda_k}|^{2(d-{\difforder})} \right)
- 2 (d-{\difforder}) \bandwidth{}^{-1} \sum_{k=1}^{\bandwidth{}} \log(|1- \rme^{\rmi \lambda_k}|) \eqsp,
\end{equation}
and we define the \emph{Local Whittle Fourier} estimator (\LWF) as the minimum
\begin{equation}
\label{eq:definition-dLWF}
\hat{d}^{\LWF}(\bandwidth{}) \eqdef \argmin_{\bar{d} \in \Rset} \tilde{L}_{\taper,\bandwidth{}}(\bar{d}) + {\difforder} \eqsp.
\end{equation}
The function $d \to \tilde{L}_{\taper,\bandwidth{}}(d)$ is convex and thus admit a single global minimum, which
can be obtained numerically by using a  standard one-dimensional convex optimization algorithm.
In \cite{robinson:1995:GSE}, the minimization in \eqref{eq:definition-dLWF} is performed over a closed interval which
was supposed to include the true value of the parameter. But, because the contrast is strictly convex, there is in fact no need to impose such a restriction.

\subsection{The exact  and the non-stationary extended local Whittle estimators}
To conclude this section, let us mention two recent works on the estimation of the memory parameter
for non-stationary $M(d)$ processes, $d \geq 1/2$, which are not covered in details in this contribution because they are derive
under slightly different conditions and henceforth do not compare well with wavelet estimators.

\cite{shimotsu:phillips:2005} introduced an \emph{exact local Whittle} estimator. It is applicable when the $M(d)$
series is generated by a linear process and when the domain of $d$ is not wider than 9/2.
Their estimator is based on fractional differencing of the data and the complexity of their algorithm is of the order
$n^2$, where $n$ is the number of observations. In contrast, the complexity of the estimators we consider  is  of the
order of $n \log_2(n)$;
see the discussion in \cite{moulines:roueff:taqqu:2007:aos}.
Note also that the model considered by~\cite{shimotsu:phillips:2005} is not an $M(d)$ process in the sense given above
and is not time-shift invariant, see their Eq.~(1).
In addition, their estimator is not invariant upon addition of a constant in the data, a drawback which is not easily dealt with, see their Remark~2.

\cite{abadir:distaso:giraitis:2007}  propose to extend the local Whittle estimator to $d \in (-3/2,\infty)$
calling it the fully extended local Whittle  estimator. This estimator is based on an extended definition of the DFT $\tilde{\dft}^X(\lambda_j,d)$, which
include correction terms, \ie\ $\tilde{\dft}^X(\lambda_j,d) \eqdef \dft^X(\lambda_j) + k^X(\lambda_j,d)$. The correction term $k^X(\lambda_j,d)$, which takes constant
values on the intervals $d \in [p-1/2,p+1/2)$, $p=0,1,\dots$  is defined as $k^X(\lambda_j,d)=0$ if $d \in (-1/2,1/2)$ and
\begin{equation}
\label{eq:correction-term}
k^X(\lambda_j,d)= \rme^{-\rmi \lambda_j} \sum_{r=1}^p (1- \rme^{-\rmi \lambda_j})^{-r} Z_{n,r} \eqsp, \quad d \in [p-1/2,p+1/2) \eqsp, \quad  p=1,2, \dots \eqsp,
\end{equation}
where $Z_{n,r} \eqdef (2 \pi n)^{-1/2} (\diffop^{r-1} X_n - \diffop^{r-1} X_0)$, $r=1,\dots,p$. The corresponding corrected periodogram $\tilde{I}^X(\lambda_j,d)$ is
given by $\tilde{I}^X(\lambda_j,d) = |\tilde{\dft}^X(\lambda_j,d)|^2$. The extended local Whittle estimator is then defined as
$\hat{d}^{\LWF} = \argmin_{\bar{d} \in [d_{\min},d_{\max}]} L^{\LWF}_{\bandwidth{}}(\bar{d})$ where $L^{\LWF}(d)$ is
\[
L^{\LWF}_{\bandwidth{}}(d) \eqdef \log \left( \bandwidth{}^{-1} \sum_{k=1}^{\bandwidth{}} \tilde{I}^X(\lambda_k,d) |1 - \rme^{\rmi \lambda_k}|^{2d} \right)
- 2 d \bandwidth{}^{-1} \sum_{k=1}^{\bandwidth{}} \log(|1- \rme^{\rmi \lambda_k}|) \eqsp.
\]
Note that to compute $k^X(\lambda_j;d)$, we have to involve additional observations $X_{-p+1}, \dots, X_n$, where $p \eqdef 0 \vee \lfloor d_{\max}-1/2 \rfloor$.
Compared to the \cite{shimotsu:phillips:2005} estimator, this estimator is easy to evaluate numerically, but the approximation of the extended local Whittle function
is not continuous at $d=1/2, 3/2, \dots$, which does not allow one to obtain limit theorems at these points (and, in the finite sample case, causes disturbances in the neighborhood of
these values). In addition,
this estimator is not robust to the presence of polynomial trends in the data (if $d$ is the memory parameter, the method tolerate a polynomial trend of degree at most $\lfloor d+1/2 \rfloor$).


\section{The semi-parametric estimation setting}
\label{sec:semi-parametric-setting}
The theory of semi-parametric Fourier estimators was developed in two fundamental papers by Robinson, \cite{robinson:1995:GPH} and
\cite{robinson:1995:GSE}, which establish, under suitable conditions, the asymptotic normality of the GPH and the
LWE estimators in the stationary case. These results were later extended to non-stationary M($d$) processes for different versions of
the memory estimator and under various sets of assumptions.
The theory of semi-parametric wavelet estimators was developed much more
recently in \cite{moulines:roueff:taqqu:2007:jtsa} and \cite{moulines:roueff:taqqu:2007:aos} (some preliminary results are in \cite{bardet:lang:moulines:soulier:2000} and \cite{bardet:2002}).
To allow for comparison  the wavelet and the Fourier approaches, the asymptotic properties of the estimators are presented
under a common set of assumptions. Because the theory of wavelet estimators is much less developed than the theory of Fourier estimators,
these assumptions  can often be relaxed in the context of Fourier estimators.

There are two types of additional assumptions that enter into play in an asymptotic theory. First, the
semi-parametric rates of convergence depends on the smoothness of the short-memory component in a neighborhood of zero frequency.
The most common assumption, introduced in \citep{robinson:1995:GPH}, is a Hölder condition on the short-memory component of the spectral density $f^\ast$ in \eqref{eq:spdeltaStat}.
\begin{definition}
\label{def:HolderClass}
For any $0 < \beta \leq 2$, $\gamma>0$ and $\varepsilon\in(0,\pi]$,  $\calH(\beta,\gamma,\varepsilon)$ is the set of all
non-negative and even function $g$ that satisfies $g(0) > 0$ and for all  $\lambda\in(-\varepsilon,\varepsilon)$
\begin{equation}\label{eq:Hbeta}
|g(\lambda)- g(0) | \leq \gamma \, g(0) \, |\lambda|^\beta \eqsp .
\end{equation}
\end{definition}
The larger the value of $\beta$, the smoother the function at the origin.
Observe that if $f^\ast$ is even -- as assumed -- and if it is in addition infinitely differentiable, then
$f^{\ast'}(0)=0$ and hence, by a Taylor expansion,~(\ref{eq:Hbeta}) holds with $\beta=2$, that is, in this case, one has
$f^\ast\in\calH(2,\gamma,\varepsilon)$.
\cite{andrews:guggenberger:2003} extend this definition to the case $\beta > 2$ by considering even functions satisfying
\begin{equation}
\label{eq:extension-phi}
\left|  g(\lambda) - g(0) -\sum_{k=1}^{\lfloor \beta/2 \rfloor} \varphi_k \lambda^{2k}/(2k!) \right| \leq \gamma g(0) |\lambda|^{\beta - 2 \lfloor \beta/2 \rfloor} \eqsp,
\end{equation}
with $|\varphi_k| \leq \gamma g(0)$, for any $k \in \{1, \dots, \lfloor \beta/2 \rfloor\}$. To take advantage of this more refined smoothness assumption when $\beta>2$, bias reduction techniques must be applied (see for example \cite{andrews:guggenberger:2003,robinson:henry:2003, andrews:sun:2004}). We will only consider $\beta\leq2$
for sake of brevity.

Second, the definition of an M($d$) process accounts only for the spectral (or equivalently  covariance) structure, which
specifies the distribution of the process if $X$ is Gaussian.
To extend the results in the non-Gaussian context, it is necessary to specify the distribution of the process beyond
its second-order properties. The most complete asymptotic theory has been developed so far for linear processes.
\begin{definition}\label{def:strongLinearMd}
We say that $X$ is a \emph{strong linear M($d$) process} if
there exist a $\ltwo(\Zset)$-sequence $\{a_s\}_{s \in \Zset}$ and an \iid\ sequence $\{Z_s \}_{s \in \Zset}$ satisfying
$\PE[Z_0]= 0$, $\PE[Z_0^4] < \infty$, and for any $t \in \Zset$,
\begin{equation}
\label{eq:strong-linear}
(\diffop^K X)_t = \sum_{s \in \Zset} a_s Z_{t-s} \eqsp, \quad K \eqdef \lfloor d + 1/2 \rfloor \eqsp.
\end{equation}
\end{definition}

\begin{remark}
\label{rem:weak-linear-process}
According to the standard terminology, \emph{strong} here refers to the fact that
$\{Z_t\}$ is \iid\ (or a strong white noise). This assumption can often be relaxed by supposing that $\{Z_t\}$ is a
martingale
difference sequence ( $\CPE{Z_{t}}{\mathcal{F}_{t-1}}= 0$
where $\mathcal{F}_t = \sigma(Z_s, s \leq t)$ is the natural filtration of the process) satisfying various additional
conditional moment assumptions. For example, \cite{robinson:1995:GSE}, and many authors after that,  assume
that $\{Z_t^2 - \PE[Z_t^2]\}$ is a square integrable martingale difference.
\end{remark}

The following theorem, established in \cite{giraitis:robinson:samarov:1997}, provides a lower bound for the estimation error.
\begin{theorem}\label{theo:lowerboundRate}
Let $d_{\min}<d_{\max}$ in $\Rset$, $\varepsilon\in(0,\pi]$, $\beta\in(0,2]$ and $\gamma >0$.
There exists a constant $c>0$ such that,
\begin{equation}
\label{eq:minimaxgrs} \liminf_{n\to\infty} \inf_{\hat d_n} \sup_{d_{\min} \leq
d \leq d_{\max}} \sup_{f^*\in\calH(\beta,\gamma,\varepsilon)}
\PP_{d,f^*}\left(n^{\beta/(2\beta+1)}|\hat d_n - d|\geq c \right)>0,
\end{equation}
where the infimum $\inf_{ \hat d_n }$ is taken over all possible
estimators based on $\{ X_1,\cdots, X_n \}$ and $\PP_{d,f^*}$ denotes
the distribution of a Gaussian M($d$) process $\{X_t\}_{t\in\Zset}$ with generalized spectral
density of the form \eqref{eq:SpectralDensity:FractionalProcess}.
\end{theorem}
We shall see in the sequel that the best possible rate $n^{\beta/(2\beta+1)}$ is achieved by both  wavelet and Fourier estimators when $X$ is a strong linear M($d$) process.

The theory of the semi-parametric estimation of $d$ for several \emph{non-linear} M($d$) processes, used in particular in
financial econometric, and in teletrafic modeling,  have also been investigated in the literature (see the recent surveys
\cite{deo:hsieh:hurvich:soulier:2006} and \cite{teyssiere:abry:2007} and the references therein).
The stochastic volatility model (a special instance of the signal plus noise model) has been considered in
\cite{hurvich:moulines:soulier:2005,deo:hurvich:lu:2006}, which establish consistency, rate of convergence
and asymptotic normality of an appropriately modified LWF estimator.  \cite{dalla:giraitis:hidalgo:2006}
provide general conditions under which the LWF estimator of the memory parameter of a stationary process is consistent and examines
its rate of convergence. This class of processes include, among others, signal plus noise processes,
nonlinear transforms of a Gaussian process, and exponential generalized autoregressive,
conditionally heteroscedastic (EGARCH) models, etc...
\cite{fay:roueff:soulier:2007} provide the consistency and rate of convergence of the LWW estimator for the
infinite source Poisson process. The results are in general not as complete as in the
linear case, and the required assumptions  are specific to each considered model (abstract assumptions like in \citep[Eq. (8)]{dalla:giraitis:hidalgo:2006}
or \citep{moulines:roueff:taqqu:2007:aos} can be used, but checking these still require model-dependent conditions).
In addition, the asymptotic behavior of the estimators are model-dependent and might depart significantly from the results
obtained in the linear case.

\section{Asymptotic properties of the wavelet estimators \LRW\ and \LWW}
\label{sec:w-results}
\subsection{The between-scale process}
%
Before stating the main known results on the asymptotic behavior of \LWW\ and \LRW\ estimators, some
additional definitions related to the spectral density of wavelet coefficients are required.

If the process $X$ is an M($d$) process, as defined in Section \ref{sec:definition-an-md}, and  $M>d-1/2$, then $\diffop^MX$ is weakly stationary.
It follows that the process $\{\dwt^X_{j,k}\}_{k\in\Zset}$ of wavelet coefficients at scale
$j\geq0$ is weakly stationary in $k$. However the two--dimensional process $\{[\dwt^X_{j,k},\,\dwt^X_{j',k}]^T \}_{k\in \Zset}$ of wavelet
coefficients at two different scales $j$ and $j'$, with $j\neq j'$, is not weakly stationary.
This is why we consider the \emph{between-scale} process
$\{[\dwt^X_{j,k},\,\bdwt^X_{j,k}(j-j')^T]^T \}_{k\in \Zset}\;,$ where
\begin{equation}\label{eq:Defbd}
\bdwt^X_{j,k}(\dj) \eqdef \left[\dwt^X_{j-\dj,2^{\dj}k},\,\dwt^X_{j-\dj,2^{\dj}k+1},\,\dots,  \dwt^X_{j-\dj,2^{\dj}k+2^{\dj}-1}\right]^T.
\end{equation}
The index $\dj$ in~(\ref{eq:Defbd}) denotes the scale difference $j-j'\geq0$ between the finest  scale $j'$ and the coarsest scale  $j$. If $\dj=0$, that is $j=j'$, then $\bdwt^X_{j,k}(0)$ is the scalar $\dwt^X_{j,k}$.
For $\dj>0$, the second component of the between--scale process is a vector whose entries are the wavelet coefficients with scale index $j-\dj=j'$
and translation indices $2^{j-j'}k,2^{j-j'}k+1,\dots,2^{j-j'}(k+1)-1$. Using~\eqref{eq:down}, it follows that
\begin{align*}
\bdwt^X_{j,k}(\dj)&=\left[(h_{j-\dj,\cdot}\star \bx)_{2^{j}k+l\;2^{j-\dj}},\;l=0,\dots,2^{\dj}-1\right]^T\\
&= \left[(\downarrow^{j}\circ\;\shift^{-l\;2^{j-\dj}} [h_{j-\dj,\cdot}\star \bx])_{k},\;l=0,\dots,2^{\dj}-1\right]^T\;,
\end{align*}
where $\shift$ is the shift operator defined, for any sequence $\{ c_k \}_{k \in \Zset}$, by $(\shift c)_k = c_{k-1}$.
The between--scale process is stationary in $k$ because all its entries can be expressed with the same downsampling operator $\downarrow^{j}$ applied to
jointly stationary time series, namely, $h_{j,\cdot}\star \bx$ and $\shift^{-l\;2^{j-\dj}} [h_{j-\dj,\cdot}\star \bx]=[\shift^{-l\;2^{j-\dj}}h_{j-\dj,\cdot}]\star \bx,\;l=0,\dots,2^{\dj}-1$.
One can therefore write, for all $0 \leq \dj \leq j$,
$$
\PCov_f(\dwt^X_{j,k},\bdwt^X_{j,k'}(\dj)) = \int_{-\pi}^\pi \rme^{\rmi\lambda(k-k')} \,
\bdens[\phi,\psi]{j,\dj}{\lambda}{f} \, \rmd\lambda \; ,
$$
where $\bdens[\phi,\psi]{j,\dj}{\lambda}{f}$ is the cross-spectral density function of the between-scale process. The case $\dj=0$
corresponds to the spectral density  of the \emph{within-scale} process $ \{ \dwt^X_{j,k} \}_{k \in \Zset}$.
When $X$ is an M($d$) process with spectral density function~(\ref{eq:SpectralDensity:FractionalProcess}), we often denote
$\bdens[\phi,\psi]{j,\dj}{\lambda}{f}$ by $\bdens[\phi,\psi]{j,\dj}{\lambda}{d,f^\ast}$.

\subsection{Generalized fractional Brownian motion}
We shall approximate the within- and between-scale spectral densities $\bdens[\phi,\psi]{j,\dj}{\lambda}{d,f^\ast}$ of the process $X$ with
memory parameter $d\in\Rset$ by the corresponding spectral densities of the \emph{generalized fractional Brownian motion} $B_{(d)}$.
This process is parametrized by a family $\Theta_{(d)}$ of smooth test functions $\theta(t)$,
$t\in\Rset$ and is defined as follows: $\{B_{(d)}(\theta),\,\theta\in\Theta_{(d)}\}$ is a mean zero Gaussian process
with covariance
\begin{equation}\label{eq:CovFBMGen}
\PCov\left(B_{(d)}(\theta_1),\,B_{(d)}(\theta_2)\right) = \int_{\Rset}
|\lambda|^{-2d}\,\hat{\theta_1}(\lambda)\,\overline{\hat{\theta_2}(\lambda)}\,\rmd\lambda\;.
\end{equation}
The finiteness of the integral $\int_\Rset|\lambda|^{-2d}\,|\hat{\theta}(\lambda)|^2\,\rmd\lambda$ provides a constraint on the
family $\Theta_{(d)}$. For instance, when $d>1/2$, this condition requires that $\hat{\theta}(\lambda)$ decays sufficiently
quickly at the origin and, when $d<0$, it requires that $\hat{\theta}(\lambda)$ decreases sufficiently rapidly at
infinity. Hence, under \allWA, $\theta$ can be a wavelet $\psi$ if $d\in(1/2-\alpha,M+1/2)$. The discrete wavelet transform of $B_{(d)}$ is defined as
\begin{equation}\label{eq:DWTFBMGen}
\dwt_{j,k}^{(d)} \eqdef B_{(d)}(\psi_{j,k}),\quad (j,k)\in\Zset \times \Zset\;.
\end{equation}
As shown in~\cite[Relation~(35)]{moulines:roueff:taqqu:2007:jtsa}, for all $j$,
$k$ and $k'$ in $\Zset$, and $\dj\geq0$, one has
\begin{align}
\label{eq:DasympAndCov}
\PCov\left(\dwt_{j,k}^{(d)},\,\bdwt^{(d)}_{j,k'}(\dj)\right)
 = 2^{2dj}\,\int_{-\pi}^\pi \bdensasymp[\psi]{\dj}{\lambda}{d} \, \rme^{\rmi \lambda (k-k')}\,\rmd\lambda \; ,
\end{align}
where $\bdensasymp[\psi]{\dj}{\lambda}{d}$  does not involve $j$ and is given by
\begin{equation}\label{eq:bDpsi}
\bdensasymp[\psi]{\dj}{\lambda}{d} =
\sum_{l\in\Zset} |\lambda+2l\pi|^{-2d}\,\be_{\dj}(\lambda+2l\pi) \,
\overline{\hat{\psi}(\lambda+2l\pi)}\hat{\psi}(2^{-\dj}(\lambda+2l\pi))\;,
\end{equation}
where, for all $\dj\geq0$,
$
\be_\dj(\xi) \eqdef 2^{-\dj/2}\, [1, \rme^{-\rmi2^{-\dj}\xi}, \dots,
\rme^{-\rmi(2^{\dj}-1)2^{-\dj}\xi}]^T,\quad\xi\in\Rset
$.
As mentioned above, $\dwt_{j,k}^{(d)}$ is well defined under \allWA\ when  $d\in(1/2-\alpha,M+1/2)$.
Under the same condition, using the decay condition of $\hat{\psi}$ in~\ref{item:psiHat} and~(\ref{eq:MVM}), one easily gets
that for any $\dj\geq0$, $\bdensasymp[\psi]{\dj}{\lambda}{d}$ is continuous $(-\pi,\pi)\setminus\{0\}$ and
$|\bdensasymp[\psi]{\dj}{\lambda}{d}|=O(|\lambda|^{2(M-d)})$ as $\lambda\to0$. If moreover $d\leq M$,
one has (see Relation~(72) in~\cite{moulines:roueff:taqqu:2007:aos})
\begin{equation}
  \label{eq:bdensasympBound}
\sup_{\dj\geq0}2^{\dj(2d-1/2)}\;\int_{-\pi}^\pi\left|\bdensasymp[\psi]{\dj}{\lambda}{d}\right|^2\;\rmd\lambda  <\infty \; .
\end{equation}

The variance $ \vj{j}{d,f^{\ast}}2 \eqdef \PVar [ \dwt^X_{j,0}] $
can be interpreted as a \emph{scale spectrum}, in words, the \emph{power} at scale $j$. It is approximated by
$f^\ast(0)\Kvar[\psi]{d} \, 2^{2jd}$ where the constant $\Kvar[\psi]{d}$ is given by
\begin{equation}\label{eq:Kpsi}
\Kvar[\psi]{d}\eqdef \int_{-\infty}^{\infty} |\xi|^{-2d}\,|\hat\psi(\xi)|^2\,d\xi \eqsp.
\end{equation}

\subsection{Uniform bounds}
Using~\cite[Theorem~1]{moulines:roueff:taqqu:2007:jtsa} and~\cite[Theorem~1]{moulines:roueff:taqqu:2007:aos}, the following result holds.

\begin{proposition}\label{prop:approximation:spectral:wavelet}
Let $X$ be an M($d$) process with $d\in\Rset$ and
$f^\ast\in\calH(\beta,\gamma,\varepsilon)$ for some $\beta,\gamma>0$
and $\varepsilon\in(0,\pi]$. Assume \allWA\ with
$(1+\beta)/2-\alpha<d<M+1/2$,
\begin{equation}\label{eq:VjApprox}
\left| \vj{j}{d,f^\ast}2 - f^\ast(0) \, \Kvar[\psi]{d} \, 2^{2jd} \right| \leq C\, f^\ast(0) \, L \, 2^{(2d-\beta)j}
\end{equation}
If moreover $\varepsilon=\pi$ and $d\leq M$, then, for all $\lambda\in(-\pi,\pi)$, $j\geq \dj \geq 0$,
\begin{equation}\label{eq:DjApprox}
\left|\bdens[\phi,\psi]{j,\dj}{\lambda}{d,f^\ast} - f^\ast(0)\, \bdensasymp[\psi]{\dj}{\lambda}{d} \, 2^{2j d}\right|
\leq C\, f^\ast(0) \, L \, 2^{(2d-\beta)j}
\end{equation}
where  $|y|$ denotes the Euclidean norm of any vector $y$.
\end{proposition}

In \eqref{eq:VjApprox} and \eqref{eq:DjApprox}, the constant $C$ only depends on $d,\beta$ and
on the wavelets $\psi$ and $\phi$. It can be made independent of $d$ on any compact set included in $((1+\beta)/2-\alpha,M+1/2)$   for  \eqref{eq:VjApprox} and
in $((1+\beta)/2-\alpha,M]$ for \eqref{eq:DjApprox}.
This proposition shows that the covariance properties of the wavelet coefficients of an M($d$) process resemble those of the
generalized FBM $B_{d}$ at large scales. The latter are not, in general, decorrelated as sometimes heuristically assumed (see \cite{abry:veitch:1998} and
\cite{veitch:abry:1999} for example).  Exact decorrelation occurs  but in very specific cases: if $\{\psi_{j,k}\}$ is an orthonormal basis of  $\ltwo(\Rset)$  and $d=0$,
see~\cite[Remark~7]{moulines:roueff:taqqu:2007:jtsa}.
Due to its very specific spectral property, the \textit{ideal} Shannon wavelet coefficients (defined in \eqref{eq:psihatShannon}) of  $B_{d}$ satisfy partial
independence for $d\neq0$. Indeed, applying \eqref{eq:psihatShannon} in~(\ref{eq:bDpsi}), we get, for all $\lambda\in(-\pi,\pi)$,
\begin{equation}
  \label{eq:ShannonBdensAsym}
  \bdensasymp[\psi_S]{\dj}{\lambda}{d}=
\begin{cases}
0 & \quad \text{for}\quad \dj\geq1\\
(2\pi-|\lambda|)^{-2d} & \quad \text{otherwise},
  \end{cases}
\end{equation}
implying that $\dwt_{j,k}^{(d)}$ and $\dwt_{j',k'}^{(d)}$ are
uncorrelated for $j\neq j'$ and that
$\dwt_{j,k}^{(d)}$ and $\dwt_{j,k'}^{(d)}$ are uncorrelated only if $d=0$.

The asymptotic behavior of wavelet estimators
$\hd^{\LRW}(\nlowscale,\nupscale,\regressweights{n})$ and
$\hd^{\LWW}(\nlowscale,\nupscale)$ defined in
\eqref{eq:definition:estimator:regression} and \eqref{eq:LWWdef}
will be derived for specific lower and upper scale sequences
$(\nlowscale)$ and $(\nupscale)$. In the semiparametric framework,
the lower scale sequence $(\nlowscale)$ governs the rate of
convergence of the memory estimator. There are two possible settings
as far as the upper scale sequence $(\nupscale)$ is concerned:
\begin{enum_J}
\item\label{it:caseJ1-J0fixed} $\nupscale-\nlowscale$ is fixed, equal to $\ell>0$
\item\label{it:caseJ1isJ} $\nupscale\leq \nmaxscale$ for all $n$ and $\nupscale-\nlowscale\to\infty$ as
$n\to\infty$, where $\nmaxscale$ is the largest available scale defined in \eqref{eq:J0J1n}.
\end{enum_J}

\subsection{Asymptotic properties of the \LRW\ estimator}
We will use the following definition, for all $i,j \geq 0$,
\begin{equation}
\label{eq:definitionAVvar}
\AVvar{d}{i,j} \eqdef \frac{4 \pi 2^{2d |i-j| } 2^{i \wedge j}}{\Kvar{d}^2} \int_{-\pi}^{\pi} \left|
  \bdensasymp{|i-j|}{\lambda}{d} \right|^2 \, d \lambda \eqsp ,
\end{equation}
with $\bdensasymp{\dj}{\lambda}{d}\in\Rset^{2^\dj}$ defined by~(\ref{eq:bDpsi}).
The following result, adapted from \cite{moulines:roueff:taqqu:2007:fractals}, applies to Gaussian M($d$) processes;
it has recently been extended to strong linear M($d$) processes in \cite{roueff:taqqu:2007}.
\begin{theorem}\label{thm:LRW:CLT}
Let $X$ be a Gaussian M($d$) process with generalized spectral density given by~\eqref{eq:SpectralDensity:FractionalProcess}
with $d\in\Rset$ and $f^\ast\in\calH(\beta,\gamma,\varepsilon)$ for some $\gamma>0$, $\beta\in(0,2]$
and $\varepsilon\in(0,\pi]$. Assume \allWA\ with
\begin{equation}
\label{eq:DrangeLRW}
(1+\beta)/2-\alpha<d\leq M \; .
\end{equation}
Let $(\nlowscale)$ be a  sequence satisfying
\begin{equation}
\label{eq:mAndnconditionForCLT}
\lim_{n \to \infty} \left\{ n2^{-(1+2\beta)\nlowscale}+ (n2^{-\nlowscale})^{-1} \right\} = 0 \eqsp.
\end{equation}
and $\regressweights{}$ be a weight of length $\ell+1$ satisfying \eqref{eq:propertyw}.
Then,  as $n\to\infty$,
\begin{equation}
\label{eq:AbryVeitchEstimatorCLTGaussianFixedWeights}
\sqrt{n2^{-\nlowscale}} \left(\hd^{\LRW}(\nlowscale,\nlowscale+\ell,\regressweights{}) - d \right) \cl \calN\left(0,
\sum_{i,j=0}^{\ell} \regressweights{i}
\AVvar{d}{i,j} \regressweights{j}\right) \eqsp .
\end{equation}
\end{theorem}

\begin{remark}
This result is stated in \cite{moulines:roueff:taqqu:2007:fractals} with $\varepsilon=\pi$. The case $\varepsilon<\pi$
can be obtained by using \cite[Corollary~2]{moulines:roueff:taqqu:2007:aos}.
\note{f}{mettre a jour le numero du corollaire en fonction des épreuves envoyées par aos}
\end{remark}

\begin{remark}
The larger the value of $\beta$, the smaller the size of the allowed range for $\td$ in~(\ref{eq:DrangeLRW}) for a given
decay exponent $\alpha$ and number $M$ of vanishing moments. Indeed the range in~(\ref{eq:DrangeLRW})
has been chosen so as to obtain a bound on the bias which corresponds to the best possible rate under the condition
$f^\ast\in\calH(\beta,\gamma,\varepsilon)$. If~(\ref{eq:DrangeLRW}) is replaced by the weakest condition $\td\in(1/2-\alpha,M]$,
which does not depend on $\beta$, the same CLT~(\ref{eq:DrangeLRW}) holds but $\beta$ in
Condition~(\ref{eq:mAndnconditionForCLT}) must be replaced by $\beta'\in(0,\beta]$. This $\beta'$ must satisfy
$1/2-\alpha<(1+\beta')/2-\alpha<\td$, that is $0<\beta'<2(\td+\alpha)-1$. When $\beta'<\beta$ one gets a slower achievable
rate in~(\ref{eq:AbryVeitchEstimatorCLTGaussianFixedWeights}).  
\end{remark}

\begin{remark}
The condition $n2^{-(1+2\beta)\nlowscale}\to0$ guarantees that the bias is negligible in the limit. The optimal rate
$n^{\beta/(1+2\beta)}$ given by Theorem \ref{theo:lowerboundRate}
is obtained with $n2^{-(1+2\beta)\nlowscale}\asymp 1$, in which case the squared bias and the variance are of the same
order of magnitude, see Theorem~3 in~\cite{moulines:roueff:taqqu:2007:jtsa} where uniform bounds of the mean square error and an exact
equivalent of the variance are given. The asymptotic equivalent of the variance is $(n2^{-\nlowscale})^{-1}\sum_{i,j=0}^{\ell} \regressweights{i}
\AVvar{d}{i,j} \regressweights{j}$ as can be expected from~(\ref{eq:AbryVeitchEstimatorCLTGaussianFixedWeights}).
\end{remark}

\begin{remark}
Theorem~\ref{thm:LRW:CLT} applies only to the setting  \ref{it:caseJ1-J0fixed} since $\nupscale=\nlowscale+\ell$.
Under the setting~\ref{it:caseJ1isJ},  $\nupscale-\nlowscale\to\infty$ as $n\to\infty$, one has to replace the
weights  $\regressweights{}$ by a sequence of weights  $(\regressweights{n})$ of lengths $\nupscale-\nlowscale+1$.
Using \cite[Proposition~4]{moulines:roueff:taqqu:2007:jtsa}, if $\varepsilon=\pi$, it is possible to extend the Theorem \ref{thm:LRW:CLT} to this setting, provided that $\regressweight{n,i}\to\regressweight{\infty,i}$
as $n\to\infty$ for all $i$ and
\begin{equation}\label{eq:Condregresseursinfinis}
\lim_{\ell\to\infty}\sum_{i>\ell}\sup_{n}|\regressweight{n,i}|2^{i/2} = 0\;,
\end{equation}
in which case one has
$$
\sqrt{n2^{-\nlowscale}} \left(\hd^{\LRW}(\nlowscale,\nupscale,\regressweights{n}) - d \right) \cl
\calN\left(0,\sum_{i,j=0}^{\infty} \regressweight{\infty,i}
\AVvar{d}{i,j} \regressweight{\infty,j}\right) \eqsp.
$$
The variance in the right-hand side of the last display is finite as a consequence of~(\ref{eq:bdensasympBound})
and~(\ref{eq:Condregresseursinfinis}).
\end{remark}
The standard theory of linear regression shows that, for any fixed $\ell\geq1$, the optimal design matrix is
$D=\AVvarInv{d,\ell}{}$, where $\AVvar{d,\ell}{}=\left[\AVvar{d}{i,j}\right]_{0\leq i,j,\leq\ell}$ is a $(\ell+1)\times(\ell+1)$
matrix. By~(\ref{eq:generalFormweights}), the
corresponding weights read
\begin{equation}
  \label{eq:wopt}
\regressopt{d,\ell}{} \eqdef  \AVvarInv{d,\ell}{} B (B^T \AVvarInv{d,\ell}{} B)^{-1}  \mathbf{b} \; ,
\end{equation}
where $B$ and $\mathbf{b}$ are defined by~\eqref{eq:bAndB} and  the associated limiting variance is
\begin{equation}
  \label{eq:VARwopt}
\AsympVarRWEopt[\psi]{d,\ell}\eqdef
{\regressopt{d,\ell}{}}^T \AVvar{d,\ell}{} \regressopt{d,\ell}{} =
\mathbf{b}^T(B^T \AVvarInv{d,\ell}{} B)^{-1} \mathbf{b} \;.
\end{equation}
Since the regression vector weights $\regressweight{}$ of length $\ell$ can be viewed as regression vector weights
of length $\ell+1$ with zero as last coordinate, we have that $\AsympVarRWEopt[\psi]{d,\ell}$ decreases as $\ell$ increases and we will denote its limit by
$\AsympVarRWEopt[\psi]{d,\infty}$. Figure~\ref{fig:AsymptoticVarianceOPT} shows that
the limit is approximately attained for $\ell\geq7$ for a standard choice of wavelet.

The optimal regression vector $\regressopt{d,\ell}{}$ cannot be used directly since it depends on
unknown the memory parameter $d$, but a plug-in method can be used as suggested by \cite{bardet:2002} in a similar context: a
preliminary consistent estimator, say
$\hd^{(1)}$, is used to estimate the optimal weights $\widehat{\regressweights{}}=\regressopt{\hd^{(1)},\ell}{}$
and then the estimator $\hd^{\LRW}(\lowscale,\upscale,\widehat{\regressweights{}})$ is applied.

A different choice of weights is suggested by \cite{abry:veitch:1998} (in a parametric context).
This choice relies on the approximation
that $\bdensasymp{\dj}{\lambda}{d}/\Kvar{d}$ is nearly zero for $\dj>0$ and nearly
constant equal to $(2\pi)^{-1}$ for $\dj=0$. This provides a diagonal approximation of $\AVvarInv{d,\ell}{}$
yielding a diagonal design matrix $D$ with diagonal entries $D_{i,i}=2^{-i}$,
$i=0,\dots,\ell$, up to a multiplicative constant. By straightforward computations, this design matrix define the following
\emph{Abry--Veitch} weights.
\begin{equation}
  \label{eq:ABweigths}
  \regressAV{\ell}{i} \eqdef \frac{(i-\eta_{\ell})2^{-i}}{2\log(2)\kappa_\ell(2-2^{-\ell})}  ,\quad i=0,\dots,\ell \;,
\end{equation}
where
\begin{equation}
\label{eq:eta-k}
\eta_\ell\eqdef   \sum_{j=0}^\ell  j \frac{2^{-j}}{2-2^{-\ell}}
\quad\text{and}\quad
\kappa_\ell\eqdef\sum_{j=0}^{\ell} (j-\eta_{\ell})^2 \frac{2^{-j}}{2-2^{-\ell}}\;.
\end{equation}

This choice, while not optimal,  is closed to it in
practice, at least for not too large values of $d$ and with a standard choice of wavelets, see
Figure~\ref{fig:AsymptoticVarianceAVweightsVSOPTweights}. One
advantage of this choice of regression vector weights stems from the fact that it does not require the use of a pilot estimator, since the weights do not depend on the unknown parameter $d$.

Let us denote, for all $\dj\geq0$ (see~\cite{moulines:roueff:taqqu:2007:aos}),
\begin{equation}\label{eq:bDint}
\intbdens[\psi]{\dj}{d}  \eqdef
\int_{-\pi}^{\pi} \left|\bdensasymp[\psi]{\dj}{\lambda}{d} \right|^2 \, \d \lambda
 =  (2\pi)^{-1}\sum_{\tau\in\Zset} \PCov^2\left(\dwt_{0,0}^{(d)},\,\dwt^{(d)}_{-\dj,\tau}\right) \;,
\end{equation}
with $\bdensasymp{\dj}{\lambda}{d}\in\Rset^{2^\dj}$ defined by~(\ref{eq:bDpsi}). In the case where the weights
$\regressweights{}$ are chosen as proposed by \cite{abry:veitch:1998} given by~\eqref{eq:ABweigths}, the asymptotic variance
in the right-hand side of~(\ref{eq:AbryVeitchEstimatorCLTGaussianFixedWeights}) reads
$$
\AsympVarWWE[\psi]{d,\ell}\eqdef\sum_{i,j=0}^{\ell} \regressAV{d,\ell}{i}
\AVvar{d,\ell}{i,j} \regressAV{d,\ell}{j} \;,
$$
where $\ell+1$ is the number of scales used in the regression. Inserting~(\ref{eq:ABweigths}) and~(\ref{eq:bDint}), we get,
for any $\ell\geq1$,
\begin{multline}
\label{eq:varLimite-k}
\AsympVarWWE[\psi]{d,\ell} = \frac{\pi}{(2-2^{-\ell})\kappa_\ell(\log(2)\Kvar[\psi]{d})^2} \times \\
\left\{\intbdens[\psi]{0}{d} +
 \frac2{\kappa_\ell}
\sum_{\dj=1}^\ell \intbdens[\psi]{\dj}{d} \, 2^{(2d-1) \dj}\,\sum_{i=0}^{\ell-\dj}\frac{2^{-i}}{2-2^{-\ell}}(i-\eta_\ell)(i+\dj-\eta_\ell) \right\}\eqsp ,
\end{multline}
where $\Kvar[\psi]{d}$ is defined in \eqref{eq:Kpsi}, and $\eta_\ell$ and $\kappa_\ell$ in \eqref{eq:eta-k}.
When $\ell$ is large, the last display can be approximated by its limit as $\ell\to\infty$, namely,
\begin{equation}
\label{eq:varLimite}
\AsympVarWWE[\psi]{d,\infty} \eqdef
\frac{\pi}{[2\log(2)\Kvar[\psi]{d}]^2} \left\{\intbdens[\psi]{0}{d} +
2\sum_{\dj=1}^\infty \intbdens[\psi]{\dj}{d} \, 2^{(2d-1) \dj} \right\} \eqsp .
\end{equation}

\subsection{Asymptotic properties of the \LWW\ estimator}

Let us now consider the \LWW\ estimator $\hd^{\LWW}(\nlowscale,\nupscale)$ defined in \eqref{eq:LWWdef}.
The following results was first proved for a Gaussian M($d$) process in \cite{moulines:roueff:taqqu:2007:aos} and then
extended to linear processes in \cite{roueff:taqqu:2007}.

\begin{theorem}\label{thm:LWW:CLT}
Let $X$ be a strong linear M($d$) process with generalized spectral density given by~\eqref{eq:SpectralDensity:FractionalProcess}.
with $d\in\Rset$ and $f^\ast\in\calH(\beta,\gamma,\varepsilon)$ for some $\gamma>0$, $\beta\in(0,2]$
and $\varepsilon\in(0,\pi]$.
Assume \allWA\ with Condition~(\ref{eq:DrangeLRW}). Let $(\nlowscale)$ be a sequence satisfying
\begin{equation}
\label{eq:mAndnconditionForCLT:LWW}
\lim_{n \to \infty} \left\{ n2^{-(1+2\beta)\nlowscale}+ \nlowscale^2(n2^{-\nlowscale})^{-1/4} \right\} = 0
\end{equation}
and $(\nupscale)$ be a sequence such that either \ref{it:caseJ1-J0fixed} or~\ref{it:caseJ1isJ} holds.
Then, as $n\to\infty$,
\begin{equation}\label{eq:clt}
(n2^{-\nlowscale})^{1/2}\,(\hd^{\LWW}(\nlowscale,\nupscale)-d)\cl\calN \left[ 0,\AsympVarWWE[\psi]{d,\ell} \right] \; ,
\end{equation}
where $\ell=\lim_{n\to\infty}(\nupscale-\nlowscale)\in\{1,2,\dots,\infty\}$ and where $\AsympVarWWE[\psi]{d,\ell}$ is given by \eqref{eq:varLimite-k} for $l < \infty$ and \eqref{eq:varLimite} for $l= \infty$.
\end{theorem}

\begin{remark}
As in Theorem~\ref{thm:LRW:CLT}, the condition $n2^{-(1+2\beta)\nlowscale}\to0$ guarantees that the bias is negligible in the
limit by imposing a sufficiently fast growth for $\nlowscale$. The condition  $\nlowscale^2(n2^{-\nlowscale})^{-1/4}\to0$
means that $n2^{-\nlowscale}$ has to grow faster than $\nlowscale^8$ which is at most of order $\log^8(n)$ and hence always
holds in the typical regime where $n2^{-\nlowscale}\asymp n^\gamma$ with $\gamma\in(0,1)$.
Relation \eqref{eq:clt} is an asymptotic normality result. If we are interested merely in the rate of convergence of the
\LWW\ estimator $\hd^{\LWW}$, then we can relax condition \eqref{eq:mAndnconditionForCLT:LWW}.
It follows from \cite[Theorems~1 and~3]{moulines:roueff:taqqu:2007:aos} that if
under the assumptions of Theorem \ref{thm:LWW:CLT}, \eqref{eq:mAndnconditionForCLT:LWW} is replaced by
$$
\lim_{n \to \infty} \left\{ \nlowscale^{-1}+ \nlowscale^2(n2^{-\nlowscale})^{-1/4} \right\} = 0\;,
$$
then
$$
\hd^{\LWW}(\nlowscale,\nupscale)=d+O_\PP\left((n2^{-\nlowscale})^{-1/2}+2^{-\beta\nlowscale}\right) \; .
$$
Hence, for $2^{\nlowscale}\asymp n^{1/(1+2\beta)}$, we get the optimal rate $n^{\beta/(1+2\beta)}$,
stated in Theorem \ref{theo:lowerboundRate}.
As shown in \cite{moulines:roueff:taqqu:2007:aos}, this result holds for a class of
"weak" linear M($d$) processes (see remark \ref{rem:weak-linear-process}).
\end{remark}
The following observations, which follow directly from Theorem \ref{thm:LWW:CLT} seems to have be unkwnown
so far:
\begin{corollary}\label{rem:compLWRasympVar}
The \LWW\ estimator has the same asymptotic
variance as the \LRW\ estimator with Abry--Veitch weights \eqref{eq:ABweigths}.
\end{corollary}
\begin{corollary}
\label{cor:optimality-LRW}
Among all  wavelet estimators of the memory parameter $d$ presented in this paper, for a given choice of
wavelet and scales involved in the estimates, the estimator with optimal asymptotic variance is the \LRW\
estimator using the optimal weights defined in~\eqref{eq:wopt}.
\end{corollary}
As explained above, the optimal \LRW\ in Corollary \ref{cor:optimality-LRW} estimator requires plugging a preliminary consistent estimator of $d$.

\subsection{Asymptotic variances}
The asymptotic variances of both the \LRW\ and the \LWW\ estimators
depend on true value of the memory parameter $d$ and on the wavelet $\psi$, as they are all expressed in terms of
$\bdensasymp{\dj}{\lambda}{d}$, defined in~\eqref{eq:bDpsi}.
In practice, one estimates the limiting variance $\AsympVarWWE[\psi]{d,\ell}$
by $\AsympVarWWE[\psi]{\hd,\ell}$ in order to construct asymptotic confidence intervals. The continuity of
$\AsympVarWWE[\psi]{\cdot,\ell}$ and the consistency of $\hd$ justify this procedure.

\begin{figure}[p]
  \centering
    \includegraphics[width=0.9\textwidth]{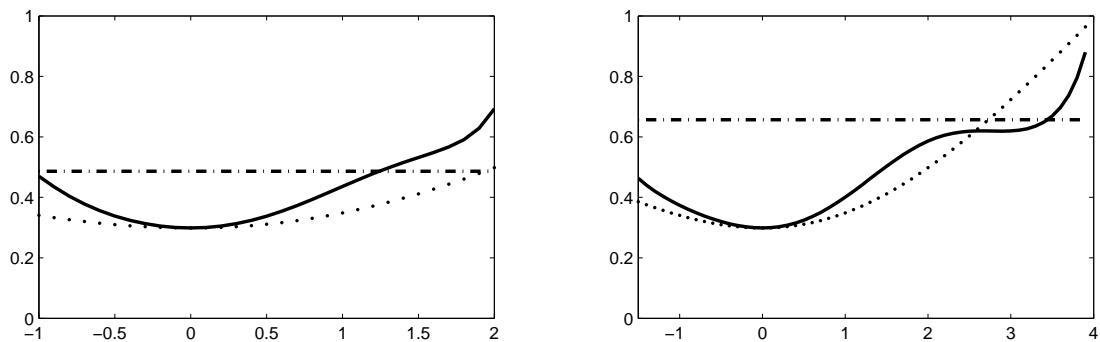}
  \caption{Comparison of the asymptotic variances of \LWF\ and  \LWW\ estimators as functions
    of $d$. The dot/dash line displays the variance \eqref{eq:definition-Phi} of the \LWF\ estimator
    with taper order \taper; the plain curve displays the variance $\AsympVarWWE[\psi]{d,\ell}$ in \eqref{eq:varLimite-k} with $\ell=7$ of the \LWW\ estimator using Daubechies wavelets of order $M$; the dotted
    curve displays the variance \eqref{eq:varLimite} of the \LWW\ estimator using the ideal Shannon wavelet. Left panel: $\taper=2$ for the \LWF, $M=2$ for the \LWW. Right panel: $\taper=4$ for the \LWF, $M=4$ for \LWW.}
  \label{fig:AsymptoticVarianceFourierVsWavelet}
\end{figure}

\begin{figure}
  \centering
    \includegraphics[width=0.9\textwidth]{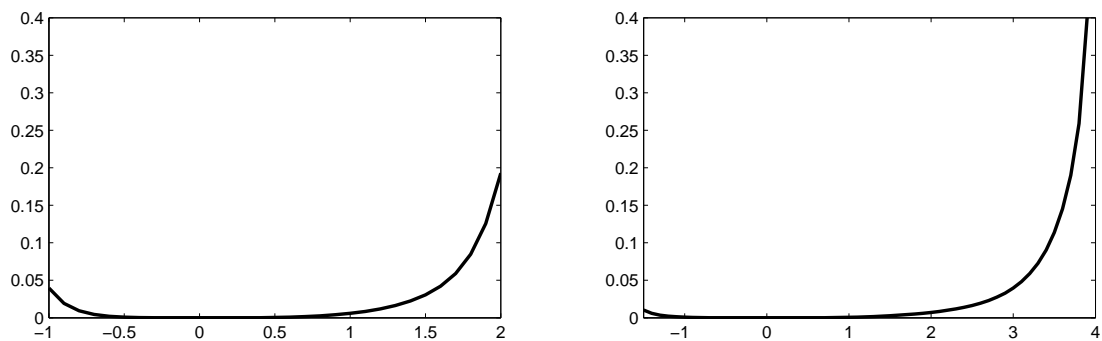}
  \caption{Comparison of the asymptotic variance of the \LRW\ estimator using Abry-Veitch weights
    given by \eqref{eq:ABweigths}
    with  the \LRW\ estimator using optimal weights given by \eqref{eq:wopt}. We plot
    $\AsympVarWWE[\psi]{d,\ell}-\AsympVarRWEopt[\psi]{d,\ell}$ as a function of $d$ with $\ell=7$.
   We used  Daubechies wavelets for two different values for $M$. Left panel: $M=2$. Right panel: $M=4$.}
  \label{fig:AsymptoticVarianceAVweightsVSOPTweights}
\end{figure}

\begin{figure}
  \centering
    \includegraphics[width=0.9\textwidth]{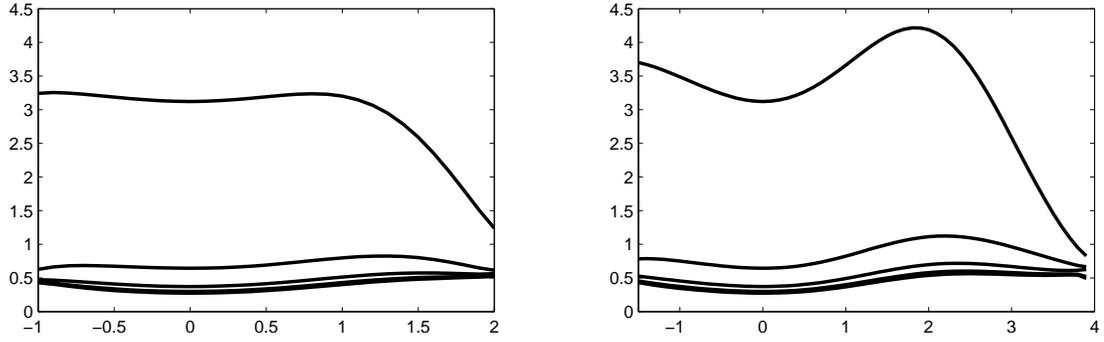}
  \caption{Comparison of the asymptotic variance \eqref{eq:varLimite-k} of the \LRW\ estimator using Abry-Veitch weights
    given by \eqref{eq:ABweigths}
    with  the one of \LRW\ estimator using optimal weights given by \eqref{eq:wopt}. We plot $\AsympVarRWEopt[\psi]{d,\ell}$ as a
    function of $d$ for successive values of  $\ell=1,3,5,7,9$ (from top to bottom).
    We used  Daubechies wavelets for two different values for $M$. Left panel: $M=2$. Right panel: $M=4$.}
  \label{fig:AsymptoticVarianceOPT}
\end{figure}

A comparison  of the asymptotic variances $\AsympVarWWE[\psi]{d,\ell}$ for several wavelets  can be found in
\cite{moulines:roueff:taqqu:2007:aos}, \note{f}{mettre a jour le numero du corollaire en fonction des épreuves envoyées par aos}
(see also Figure~\ref{fig:AsymptoticVarianceFourierVsWavelet})
In particular,   as Figure~1 in~\cite{moulines:roueff:taqqu:2007:aos} indicates, the choice of wavelets does not
matter much (provided that $(1+\beta)/2-\alpha<d\leq M$ holds)
and a sensible approximation can be obtained by using the Shannon wavelet, for which a simple expression of the
asymptotic variance can be obtained thanks to~\eqref{eq:ShannonBdensAsym}.
Using the Shannon wavelet in~(\ref{eq:bDpsi}), we get, for all $\lambda\in(-\pi,\pi)$,
$\bdensasymp[\psi_S]{\dj}{\lambda}{d}=0$ for $\dj\geq1$ and
$\bdensasymp[\psi_S]{0}{\lambda}{d}=(2\pi-|\lambda|)^{-2d}$ so that, for all $d\in\Rset$,~(\ref{eq:varLimite-k}) becomes
\begin{equation}
  \label{eq:shannonApprox}
 \AsympVarWWE[\psi_S]{d,\ell} = \frac{\pi \; g(-4d)}
{2(2-2^{-\ell})\kappa_\ell\log^2(2) \, g^2(-2d)}
\quad\text{where}\quad g(x)=\int_\pi^{2\pi} \lambda^x \,\d\lambda \;.
\end{equation}

\section{Asymptotic properties of the Fourier estimators \GPH\ and \LWF}
\label{sec:f-results}
\subsection{Asymptotic properties of the \GPH\ estimator}

The consistency and asymptotic normality of the \GPH\ have been established by
\cite{robinson:1995:GPH}
for stationary invertible Gaussian $M(d)$ process $-1/2 < d < 1/2$ with no data taper ($\taper=0$) and no differencing (${\difforder}=0$).
As shown by \cite{velasco:1999:JE}, the \GPH\ estimator with ($\taper=0$ and ${\difforder}=0$) exhibits non-standard behavior when $d > 1/2$.
Although it is  consistent for $d \in (1/2,1]$ and asymptotically normally distributed for $d \in (1/2,3/4)$,
the \GPH\ estimator has a non-normal limit distribution for $d \in [3/4,1]$,
and for $d > 1$, it converges to 1 in probability and is inconsistent. Hence, the interest in applying the \GPH\
estimator under differencing ${\difforder}>0$ and tapering $\taper> 0$.

The following result is adapted from \cite{moulines:soulier:2003}.
To state the results, some additional notations are required. If $\{Z_t\}$ is a Gaussian white noise,  $\bar{I}^{Z}_{\pooling,\taper}(\tilde{\lambda}_k)$ is distributed as
$\|G_{\pooling,\taper}\|^2/2$ where $G_{\pooling,\taper}= [G_{\pooling,\taper}^{(1)}, \dots, G_{\pooling,\taper}^{(2\pooling)}]$
is a $2\pooling$-dimensional zero-mean Gaussian vector with covariance matrix $\Sigma_{\pooling,\taper}$, whose expression is given in \cite{hurvich:moulines:soulier:2002}. Define
\begin{equation}
\label{eq:defgammampsigmamp}
\gamma_{\pooling,\taper} = \PE \left[\log(\|W_{\pooling,\taper}\|^2/2) \right] \eqsp, \quad
\sigma_{\pooling,\taper}^2 = \PVar \left[\log(\|W_{\pooling,\taper}\|^2/2) \right] \eqsp.
\end{equation}
Numerical expressions for these quantities are given in \cite{hurvich:moulines:soulier:2002}.
\begin{theorem}
\label{theo:gphgaus}
Assume that $X$ is a Gaussian $M(d)$ process and $f^\ast\in \calH^\ast(\beta,\gamma,\varepsilon)$
for some $\beta \in (0,2]$, $\gamma>0$, $\varepsilon \in (0,\pi]$.
Let ${\difforder}\geq 0$ be the differencing order, $\taper \geq 0$ be the tapering order, and $\pooling \geq 1$
be the pooling order. Let $\bandwidth{n}$ be a non-decreasing sequence of integers such that
\begin{equation}
\label{eq:ratemgph}
\lim_{n\rightarrow\infty} ( \bandwidth{n}^{-1} + \bandwidth{n}^{2 \beta+1} n^{-2 \beta} ) = 0.
\end{equation}
Then, for any $d$ satisfying
\begin{equation}
\label{eq:range-condition-memory-parameter}
{\difforder}-\taper-1/2<d<{\difforder}+1/2  \eqsp,
\end{equation}
the \GPH\ estimator defined in \eqref{eq:definition-GPH-estimator} satisfies,
\begin{equation}
\label{eq:result-GPH-estimator}
\sqrt{\bandwidth{n}} ( \hat{d}^{\GPH}(\bandwidth{n}) - d ) \cl \calN(0,\sigma_{\pooling,\taper}^2/4) \eqsp.
\end{equation}
where $\sigma^2_{\pooling,\taper}$ is defined in \eqref{eq:defgammampsigmamp}.
\end{theorem}
Compared to the wavelet estimators, the size of the confidence
intervals does not depend on $d$, which may be seen as a significant advantage. On the other hand, there is an inflation of the variance over all the interval $({\difforder}-\taper-1/2,{\difforder}+1/2)$ and it can be greater than the limiting variance obtained using by the wavelet
estimator, at least  for certain values of the memory parameter.

  The definition of the pooled periodogram \eqref{eq:pooled-periodogram} implies that the number of Fourier
  frequencies used to evaluate $\hat{d}^{\GPH}(\bandwidth{n})$  is equal to $(\pooling+\taper)\bandwidth{n}$.
  Hence the efficiency ratio between two \GPH\ estimators using \emph{same} number of Fourier frequencies but
  two \emph{different} pooling numbers, say $\pooling$ and $\pooling'$, may be expressed as
  $(\pooling+\taper)\sigma_{\pooling,\taper}^2/(\pooling'+\taper) \sigma^2_{\pooling',\taper}$. As shown in \cite[Theorem 1]{hurvich:moulines:soulier:2002},
  the function $\pooling \mapsto (\pooling+\taper)\sigma_{\pooling,\taper}^2$ is decreasing showing that pooling increases asymptotic efficiency.
  In addition,  for any fixed $\taper$, $\lim_{\pooling\to\infty} (\pooling+\taper)\sigma_{\pooling,\taper}^2 = \Phi(\taper)$, where
  \begin{equation}
  \label{eq:definition-Phi}
  \Phi(\taper) = \frac{\Gamma(4\taper+1) \Gamma^4(\taper+1)}{\Gamma^4(2\taper+1)} \eqsp.
  \end{equation}
  As seen below, this efficiency bound is achieved by the local Whittle estimator.
Therefore, the order of pooling can be made arbitrarily large, and at least asymptotically, an increase in the pooling order will result in a decrease of the asymptotic variance; see \cite{hannan:nicholls:1977} and Theorem \ref{theo:gphgaus}.
In practice, of course, this is not possible and since the improvements in asymptotic efficiency happen quickly, there is no need to use a very large \pooling; $\pooling=3,4,5$ are typical values.

\begin{remark}
As observed in \eqref{eq:range-condition-memory-parameter}, the differencing order \difforder\ and the taper order \taper\ control the
range of values of the memory parameter $d$ which can be inferred.
The number of differentiation ${\difforder}$ controls the upper bound for $d$,
while the taper order \taper\ controls the range. These two parameters are independent, by choosing the
number of differentiations ${\difforder}$ and the taper order
\taper\ we can therefore cover any intervals of admissible values
for $d$ (the same comment apply to the \LWF\ estimator). If
$\taper=0$, the interval over which the \GPH\ estimator is
consistent and asymptotically normal is
$({\difforder}-1/2,{\difforder}+1/2)$. If $\taper > 0$, the range
is $[{\difforder}-\taper-1/2,{\difforder}+1/2]$, as indicated in Theorem \ref{theo:gphgaus}. Note that these ranges
may not be optimal: for $\taper= 0$, using the sharpened results from
\cite{velasco:1999:JE}, the range over which the memory parameter is asymptotically normal can be shown to be
$({\difforder}-3/4,{\difforder}+3/4)$.
\end{remark}

\begin{remark}
It is possible to replace the Gaussian assumption by the weaker assumption that the process $X$ is a strong linear $M(d)$
process by adding moment and regularity conditions on the distribution of the driving noise in the definition
\eqref{eq:strong-linear}. In this case however, the estimator $\hat{d}^{\GPH}(\bandwidth{n})$ should be slightly
modified to avoid a number of Fourier frequencies near 0. In addition,  tapering and  pooling are then required, even if
the process $X$ is stationary and invertible; see \cite[Theorem 3]{velasco:2000:NG} and \cite{fay:moulines:soulier:2004}.
\end{remark}

\subsection{Asymptotic properties of the \LWF\ estimator}

The consistency and asymptotic normality of the \LWF\ estimator have been established by
\cite{robinson:1995:GSE}
for stationary invertible linear $M(d)$ process $-1/2 < d < 1/2$ with no data taper ($\taper=0$) and no differencing (${\difforder}=0$) (under the weaker assumption
that $\{Z_t\}$ in \eqref{eq:strong-linear} are martingale differences, whose squares, centered at their expectation, are also
weakly stationary martingale differences). \cite{velasco:1999} has shown that the \LWF\ estimator with $\taper=0$ and ${\difforder}=0$ was consistent for $d \in (-1/2,1]$
and asymptotically $\calN(0,1/4)$ for $d \in (-1/2,3/4)$ under the same assumptions than \cite{robinson:1995:GSE}.
\cite[Theorem 2]{hurvich:chen:2000} established Theorem \ref{theo:LWEestimator} for $\taper=1$ and ${\difforder}=1$.
This result was later extended in \cite{moulines:soulier:2003} to general $\taper$ and ${\difforder}$. The consistency of the \LWF\
was established (with ${\difforder}=0$ and $\taper=0$) for $-1/2 < d < 1/2$ for a general class of non-linear processes in \cite{dalla:giraitis:hidalgo:2006}.

\begin{theorem}
\label{theo:LWEestimator}
Assume that $X$ is a strong linear $M(d)$ process for some $d \in \Rset$ and $f^\ast\in \calH^\ast(\beta,\gamma,\varepsilon)$
for some $\beta \in (0,2]$, $\gamma>0$, $\varepsilon \in (0,\pi]$.  Let ${\difforder}$ be the differencing order and
$\taper$ be the taper order. Let $\bandwidth{n}$ be a non decreasing sequence of integers such that
\begin{equation}
\label{eq:ratemgse}
\lim_{n\rightarrow\infty} \left( \bandwidth{n}^{-1} + \bandwidth{n}^{2\beta+1} n^{-2 \beta}  \right) = 0 \eqsp.
\end{equation}
Then, the \LWF\ estimator defined in \eqref{eq:definition-dLWF} satisfies, for any
$d$ satisfying 
\begin{equation}
  \label{eq:rangeLWF}
{\difforder}-\taper-1/2<d<{\difforder}+1/2\;,  
\end{equation}
\begin{equation}
\label{eq:result-LWF-estimator}
\sqrt{\bandwidth{n}} (\hat d^{\LWF}(\bandwidth{n}) - d) \cl \calN(0,\Phi(\taper) /4) \eqsp,
\end{equation}
where $\Phi(\taper)$ is defined in \eqref{eq:definition-Phi}
\end{theorem}
The quantity $\Phi(\taper)$ quantifies the loss of efficiency due to tapering.
As the taper order increases, the limiting variance inflates:
$\Phi(0) = 1$ (no tapering), $\Phi(1) = 1.5$, $\Phi(2) = 35/18$, etc.
Since the \LWF\ estimator is based on linear functionals of the periodogram, pooling is irrelevant. Recall that,
in the definition \eqref{eq:definition-dLWF}, the standard periodogram is used and not the pooled one. This is why the
pooling order does not appear in the conditions of Theorem~\ref{theo:LWEestimator}.

\begin{remark}
The condition on the bandwidth \eqref{eq:ratemgse} is slightly weaker than Assumption A$4'$ in
\cite{robinson:1995:GSE} and \cite{hurvich:chen:2000}. It seems that the $\log^2(\bandwidth{n})$ term in
these assumptions is superfluous (see \cite[Comments of Assumption 4]{andrews:sun:2004} for the required adaptation of the
proof).
\end{remark}


\section{Discussion}\label{sec:discussion}
\label{sec:disc-concl-remarks}

The Fourier and wavelet estimators of the memory parameter present similar characteristics and distinctive advantages.
In Table~\ref{tab:FvsW}, we summarize the main features of the estimators considered in this paper.

\begin{table}[h]

\begin{tabular}{c|c|c}
                 & Fourier & Wavelets \\
\hline
\begin{tabular}{c}
Non-stationarity \\
($d$ large)
\end{tabular}
&
\begin{tabular}{c}
Pre-apply $(I-B)^{\difforder}$ \\
with ${\difforder}> d-1/2$
\end{tabular}
&
\begin{tabular}{c}
take $M\geq d$\\
(a sufficient number \\
of vanishing moments)
\end{tabular}
\\
\hline
Polynomial trends  & same as above & same as above \\
of degree $K$  &   with ${\difforder}\geq K+1$ &  with $M\geq K+1$ \\
\hline
\begin{tabular}{c}
Leakage \\
($d$ small)
\end{tabular}
&
\begin{tabular}{c}
Use taper $(1-\rme^{2\rmi\pi k/n})^p$ \\
with taper order $p>1/2-d$  \\
\end{tabular}
&
\begin{tabular}{c}
take $\alpha> (1+\beta)/2-d$ \\
(sufficiently\\
 smooth wavelet).
\end{tabular}\\
\hline
Rate of convergence & $n^{\beta/(1+2\beta)}$ ($\beta\leq2$) & $n^{\beta/(1+2\beta)}$  ($\beta\leq2$) \\
of $\hat{d}$ to $d$ & with $\bandwidth{n}\asymp n^{\beta/(1+2\beta)}$& with $2^{\nlowscale}\asymp n^{1/(1+2\beta)}$ \\
\hline
\begin{tabular}{c}
Asymptotic \\
variance
\end{tabular}
&
\begin{tabular}{c}
depends on \\
taper order $p$ only;\\
\GPH's $\downarrow$ \LWF's,\\
as pooling order $\to\infty$.
\end{tabular}
&
\begin{tabular}{c}
depends on $d$ and $\psi$;\\
\LRW's with $\regressweights{}^{\mathrm{AV}}$ $=$ \LWW's\\
$\geq$ \LRW's with $\regressopt{\hat{d}^{(1)},\ell}{}$.
\end{tabular}
\end{tabular}

\bigskip

\caption{Fourier VS wavelets: trends, non-stationarity, non-invertibility. In the \emph{Wavelets} column, $M$ and $\alpha$
  are defined in \allWA.}
\label{tab:FvsW}

\end{table}

To allow comparison between wavelet and Fourier estimators, we must first link the normalization factors, $(n2^{-\lowscale_n})^{1/2}$
for wavelet estimators and  $\bandwidth{n}^{1/2}$ for Fourier estimators.
A Fourier estimator with bandwidth $\bandwidth{n}$ projects the observations $[X_1\,\,\dots\,\,X_n]^T$ on the space generated by the vectors
$\{\cos(2\pi k \cdot /n),\sin(2\pi k \cdot /n) \}$, $k=1,\dots,\bandwidth{n}$, whose
dimension is $2\bandwidth{n}$; on the other hand, the wavelet coefficients $\{\dwt_{j,k}, j\geq \lowscale, k=0,\dots,\nj{j}-1\}$ used in the
wavelet estimator correspond to a projection on a space whose dimension is at most
$\sum_{j=\nlowscale}^{\nmaxscale}\nj{j} \sim \sum_{j=\nlowscale}^\infty n 2^{-j}  \sim 2n2^{-\lowscale_n}$.
Hence, for $\bandwidth{n}$ or
$n2^{-\lowscale_n}$ large, it makes sense to consider $n2^{-\lowscale_n}$ as an analog of the bandwidth parameter $\bandwidth{n}$.

We shall now compare the asymptotic variances in the CLT's in Theorems \ref{thm:LRW:CLT}, \ref{thm:LWW:CLT}, \ref{theo:gphgaus} and \ref{theo:LWEestimator}.
While the asymptotic variance of the Fourier estimators is a constant, the variance of the wavelet estimators is a function of the memory parameter,
which can be numerically computed. For the Fourier estimators, the  allowed range of the memory parameter $d$ is given by (see
Theorems~\ref{theo:gphgaus} and~\ref{theo:LWEestimator})
\begin{equation}
  \label{eq:allowedRangeFourier}
{\difforder}-\taper-1/2<d<{\difforder}+1/2 \; .
\end{equation}
The length of this range equals $\taper+1$, while the differentiation
order ${\difforder}$ allows to shift it towards large values of $d$. For instance, if one wishes to shift the upper boundary
of the range towards large values of $d$ while keeping the lower boundary unchanged, one has to increase both $\taper$ and
$\difforder$ by the same factor. As shown in Theorem \ref{theo:LWEestimator}, increasing $\taper$ inflates the asymptotic
variance of the estimator. For wavelet estimators, the  allowed range of the memory parameter $d$ is given by
\begin{equation}
  \label{eq:allowedRangeWavelet}
1/2-\alpha<d\leq M \; ,
\end{equation}
[see Theorems~\ref{thm:LRW:CLT} and~\ref{thm:LWW:CLT}; here we took $\beta$ arbitrarily small, since we focus on the asymptotic
variance in this discussion].
Of course one may still shift this range by a factor ${\difforder}$ to the right by differentiating the series $X$ at the
order ${\difforder}$ before processing the wavelet transform. This will also eliminate polynomial trends up to degree $M+\difforder-1$.

Observe that the higher  the $\alpha$ in~(\ref{eq:allowedRangeWavelet}), the more negative $d$ is allowed to be. This is because
the higher  the $\alpha$, the smoother the wavelet $\psi_{j,k}$ and hence the better the spectral resolution of the
wavelet. This matters particularly when $d<0$ because $f(\lambda)$ is then very small around the origin making it harder to
estimate $d$. In Fourier (see~(\ref{eq:allowedRangeFourier}), it is $\taper$ that plays a role similar to $\alpha$.

It is important to note that, for a given wavelet family such as Daubechies and Coiflets, increasing $M$ yields a larger
$\alpha$, so that the allowed range is effectively increased by increasing $M$. In contrast to Fourier methods, by increasing
the number of vanishing moments $M$, say of a Daubechies wavelet, the asymptotic variance converges to the asymptotic
variance obtained with the Shannon wavelet, presented in
\eqref{eq:shannonApprox}.
Thus, for a given $d$, there is asymptotically no price to pay for increasing the number of vanishing moments $M$ and the
number of available scales.
This is an argument  in favor of wavelet estimators as compared to tapered Fourier estimator.
This should, however, \note{f}{I removed :nevertheless} be interpreted with  care. For a given sample size $n$,
an increase of $M$ causes an increase of the support of the wavelet and a decrease in the number of available
scales. While this does not influence the asymptotic variances, it  affects the performance on finite
samples\footnote{The same
  remark apply to the pooling order for the \GPH\ estimator: the asymptotic variance decreases as $\pooling\to\infty$ but in
  practice, one takes small values, \eg $\pooling=3,4$.}.

The plots in Figure~\ref{fig:AsymptoticVarianceFourierVsWavelet}
indicate that the asymptotic variance of the
\LWW\ estimator is lower than the one obtained using the tapered version of \LWF\ estimator, for most values of the memory
parameter and this difference increases as $\taper$ increases in order to adapt to larger ranges for $d$.
The asymptotic variance $\AsympVarWWE[\psi]{d,\ell}$ in \eqref{eq:varLimite-k} of the \LWW\ estimator is
displayed for  $\ell=7$ using the Daubechies wavelets with $M=2$ (Left) to $M=4$ (Right).
For these choices of wavelets, the corresponding $\alpha$'s are $1.34$ and $1.91$, and
the allowed intervals for $d$ are $[-0.84,2]$ and $[-1.41,4]$, respectively. The asymptotic variances
$\AsympVarWWE[\psi]{d,\infty}$ in \eqref{eq:varLimite}  for $M=2$ and $4$ can also be compared to the one
of the Shannon wavelet on this plot. The asymptotic variance $\Phi(\taper)$ in \eqref{eq:definition-Phi} of the \LWF\ estimator is constant in $d$ but increases
when the taper order $\taper$ increases from $\taper=2$ to $\taper=4$, these values corresponding to intervals
lengths for $d$ close to those of the $M=2,4$ wavelet estimators (a bit larger for the former: 3 versus 2.84,  and smaller for the latter: 5 versus 5.41).

Using wavelet to estimate the memory parameter has several additional benefits compared to using Fourier estimators.
The wavelets present a rich time/frequency representation of the process, which can be more informative than that of the classical Fourier analysis,
as discussed in \cite{serroukh:walden:percival:2000}, \cite{stoev:taqqu:park:2006} and \cite{percival:walden:2006}.
Wavelets can be used to detect the presence of outliers or jumps in the mean.
The short-range dependence of the wavelet coefficients suggests construction of bootstrap confidence intervals for functionals of the wavelet coefficients, a procedure referred to as
\emph{wavestrapping}. This technique, which still is not rigorously justified, may be used to construct bootstrapped confidence interval
for the memory parameter; see for example \cite{percival:sardy:davison:2000}.

\section{A Monte-Carlo study}
\label{sec:monte-carlo}
In this section, we present some Monte-Carlo simulation results that compare
the root-mean square error performance of our four  estimators for finite samples. The
four estimators are denoted \GPH\ (Geweke-Porter-Hudak), \LWF\ (local Whittle
Fourier), \LWW\ (local Whittle wavelet) and \LRW\ (local regression wavelet).
We consider three models and several parameter combinations for each model:
\begin{enumerate}
\item The ARFIMA models, introduced by \cite{granger:joyeux:1980}, and
  generalized here to any value of the memory parameter $d$. We considered
  the ARFIMA(0,$d$,0) and ARFIMA(1,$d$,0) with $d$ in $\{-1.2, 0, 0.3, 1.5, 2.5,
  3.5\}$ and sizable lag 1 AR coefficient equal to 0.8. The innovation is assumed to be Gaussian.
  The short-memory component $f^*$ of the spectral density satisfies $f^* \in \mathcal H(2,\gamma,\pi)$, where $\mathcal{H}$
  is defined in Definition~\ref{def:HolderClass}.
\item The DARFIMA models, defined in \cite{andrews:sun:2004}, is an ARFIMA-like
  process that has a discontinuity in its spectral density at frequency $\lambda=\lambda_0$.
  The DARFIMA(1,$d$,0) has the spectral density of an ARFIMA(1,$d$,0) on the interval $[-\lambda_0,\lambda_0]$
  and is zero for $|\lambda| \in [\lambda_0,\pi]$.  It is obtained
  by low-pass filtering of an ARFIMA(1,$d$,0) trajectory by a truncated \textit{sinc} function in the time domain.
  We chose $\lambda_0 = \pi/2$ and Gaussian innovations.
\item The third model is a non-linear function of a  Gaussian sequence: $X_t = G(Y_t)$ where $\{Y_t\}$
is a stationary Gaussian sequence with zero-mean and variance 1 and $G: \Rset \to \Rset$ is a measurable function such that
$\PE [G^2(Y_0)] < \infty$. Then, $X_t$ may be expressed as the sum $X_t= c_0 + \sum_{k=k_0}^\infty (c_k / k!) H_k(Y_t)$,
where $H_k(\cdot)$ is the $k$-th Hermite polynomial and $c_k= \PE[G(Y_t) H_k(Y_t)]$. The minimal integer $k_0 \geq 1$ such that
$c_{k_0} \ne 0$ is called the Hermite rank of $G$. If  $Y$ is a $M(d)$ process with memory parameter
$d_Y \leq 1/2$, then  $X$ is also an $M(d)$ process with memory parameter
$d_X = \tfrac 1 2 \left(1 - k_0(1-2d_Y) \right)$ (see \citep[p. 229, Eq. (55)]{dalla:giraitis:hidalgo:2006} for details).
In simulations, we use $G(x) = \exp(x)$ (for which  $k_0=1$) and
$G(x) = H_2(x)= x^2-1$ (for which $k_0=2$) and denote those models SUBORD1 and SUBORD2, respectively.
\end{enumerate}
For the estimators \LWW\ and \LWF, we have used
a convex minimization procedure of the contrast functions~(\ref{eq:TildeJdef}) and~(\ref{eq:profile:likelihood}). In all
cases, 1000 simulation runs for each value of $d$ are used. This produces simulation standard errors that are roughly 3\%.

The tuning parameters of each estimation procedure have been chosen to allow a
fair comparison of those methods in a realistic setting, where  the order of magnitude of
the memory parameter $d$ is only loosely known and where one may be in the presence of high-order polynomial trends.
In order  to cover all the values of $d$ above ($-1.2\leq d \leq 3.5$), we have used a Daubechies wavelet with $M=4$ vanishing
moments  for the wavelet estimators (hence $\alpha\approx1.91$, see Table~\ref{tab:alpha}) and we have differenced the series 
${\difforder}=4$ times and have used a taper order $\taper=5$ for the Fourier estimators. The corresponding admissible ranges
are (see~(\ref{eq:allowedRangeWavelet}) and~(\ref{eq:allowedRangeFourier})) $(-1.41,4]$ and $(-1.5, 4.5)$, respectively.
For the \GPH\ estimator, we took $\pooling=4$ in Relation~\eqref{eq:pooled-periodogram} defining the 
pooled periodogram. This reduces the number of frequencies by a factor $\taper+\pooling=9$ (see
Relation~(\ref{eq:pooled-periodogram})).\note{f}{i think the factor is $\taper+\pooling=9$ and not $2(\taper+\pooling)$} 
In the case of the \LRW\ estimator, the computation of the
optimal weights in the least-square criterion is numerically quite  involved  so we ran
the simulations using the weights suggested by \cite{abry:veitch:1998}. The
difference in the results becomes significant only when $d$ gets close to the boundaries of the admissible range, see
Figure~\ref{fig:AsymptoticVarianceAVweightsVSOPTweights} for the asymptotic variance.   
The remaining free parameters are the number of frequencies (\LWF) or blocks of
frequencies (\GPH\ with pooling) denoted $m_n$, and the minimal (\ie finest)
dyadic wavelet scale $L_n$ (for the \LRW\ and the \LWW\ estimators); the highest (\ie coarsest) scale is
chosen to be the highest available ($U_n = J_n$).

\begin{sidewaysfigure}[p]
  \centering
  \includegraphics[width=1.0\textwidth]{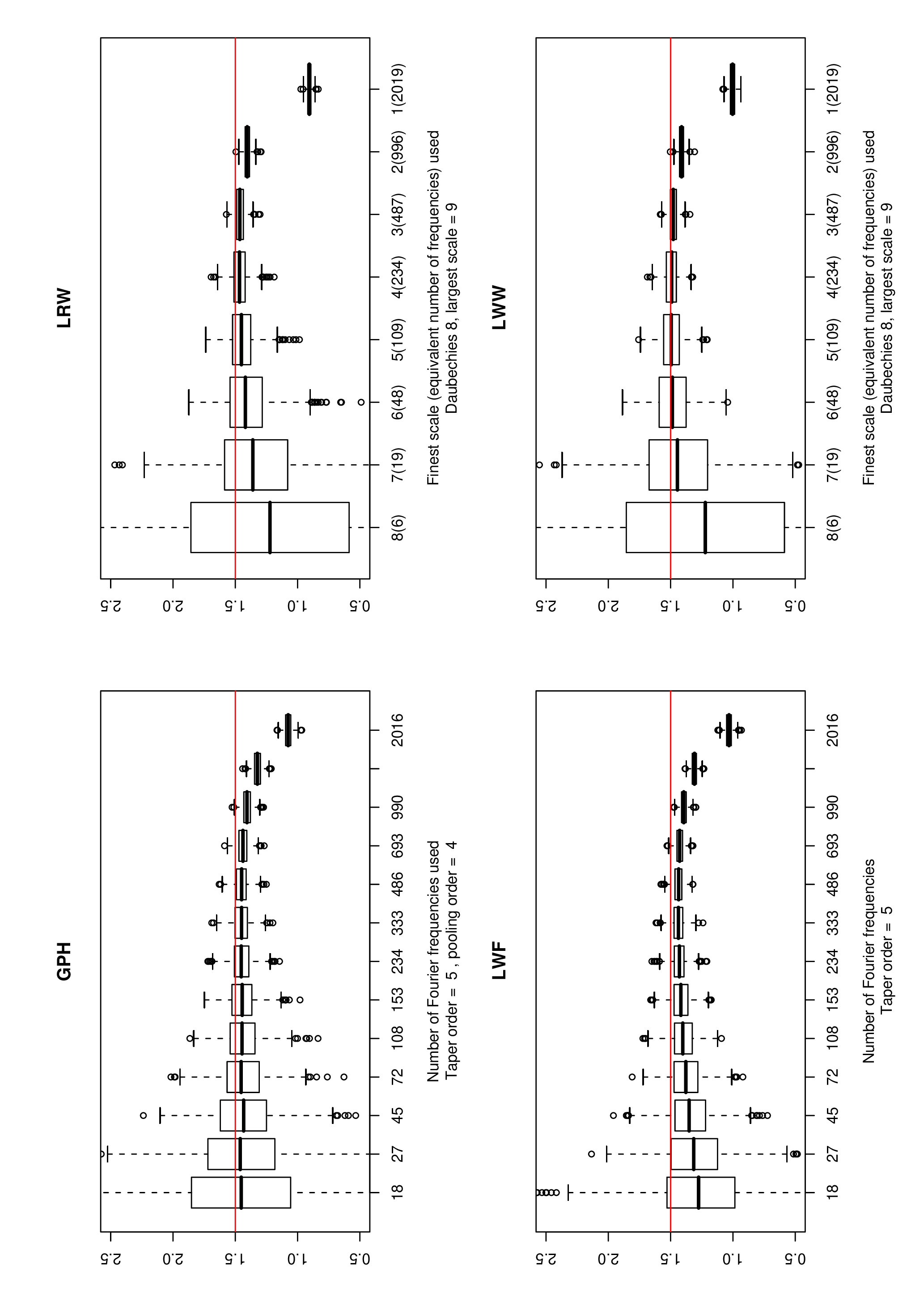}
  \caption{Distribution of the estimators for an ARFIMA(1,1.5,0)}
  \label{fig:boxplot1}
\end{sidewaysfigure}

\begin{sidewaysfigure}[p]
  \centering
  \includegraphics[width=1.0\textwidth]{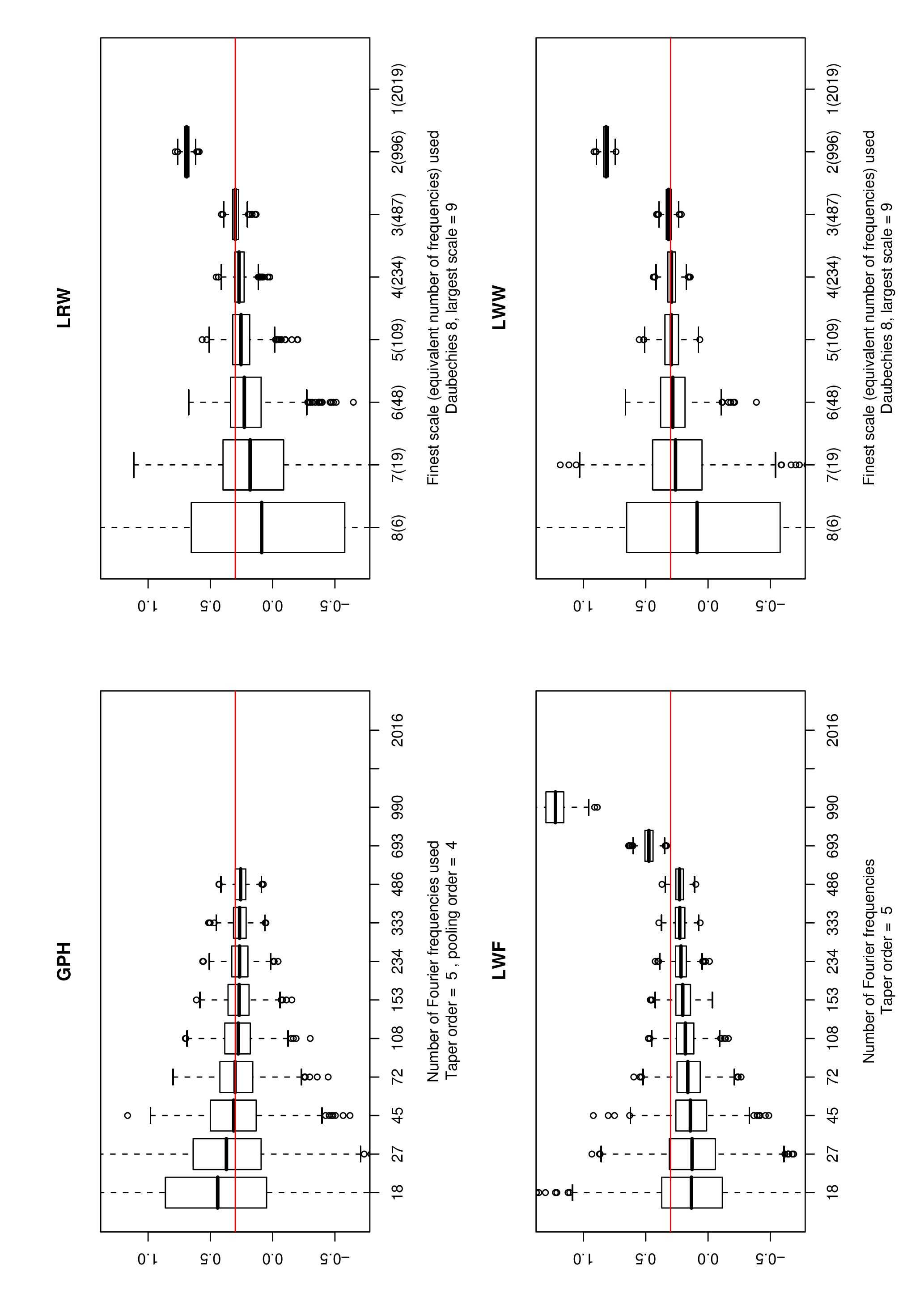}
  \caption{Distribution of the estimators for a DARFIMA(1,0.3,0)}
  \label{fig:boxplot2}
\end{sidewaysfigure}

The box and whisker plots of the estimators are displayed in
Figures~\ref{fig:boxplot1} and~\ref{fig:boxplot2} for different values of the
bandwidth (Fourier methods) and the finest scale (wavelet methods). In
Figure~\ref{fig:boxplot1}, the model is an ARFIMA(1,$d$,0) with $d=1.5$.  In
Figure~\ref{fig:boxplot2}, the model is an DARFIMA(1,$d$,0) with $d=0.3$. The
AR coefficient is 0.8 in both cases.  These figures illustrate the
bias-variance trade-off inherent to semi-parametric methods (the variance
decreases as the bandwidth or the number of scales increases, but then the bias
increases). In general, the standard deviation of the 1000 runs of the wavelet
methods is comparable to  that of the Fourier methods.

Tables~\ref{tab:results512} and \ref{tab:results4096} give the bias, variance,
and RMSE (root mean square error) for models 1,2,3 for sample sizes 512 and
4096, respectively. Those quantities are computed for the optimal bandwidth
$m_n$ (Fourier methods) or the optimal finest scale $L_n$ (wavelet methods) in
the RMSE-sense, whose values are displayed in the fourth column. For each
model, the lowest RMSE among the four methods appears in boldface. Note that
all the possible values of finest scale $L_n$ have been considered, but only a
subset of the many possible values of the bandwidth $m_n$.  The standard
deviations of the \LWW\ estimator are lower than those of the \LRW\ estimator,
which is consistent with our theoretical findings. Also, the standard
deviations of the Fourier methods remain approximately constant for the
different values of the memory parameter, whereas the variance of the wavelet
methods increase with $|d|$. Also, as predicted by the expressions of the
limiting variance, the variance of the wavelet methods are lower than those of
the Fourier methods, especially when $d$ is small.  The reported values of the
standard deviations agree with our theoretical findings for the sample size
$n=4096$.

For the non-linear processes, the results suggest that the wavelet estimator
remains consistent. However, the presence of non-linearity worsens the behavior
of the estimator at a given finite sample and a larger sample size is required
to achieve a prescribed accuracy.

\begin{figure}[p]
  \centering
  \includegraphics[width=0.8\linewidth]{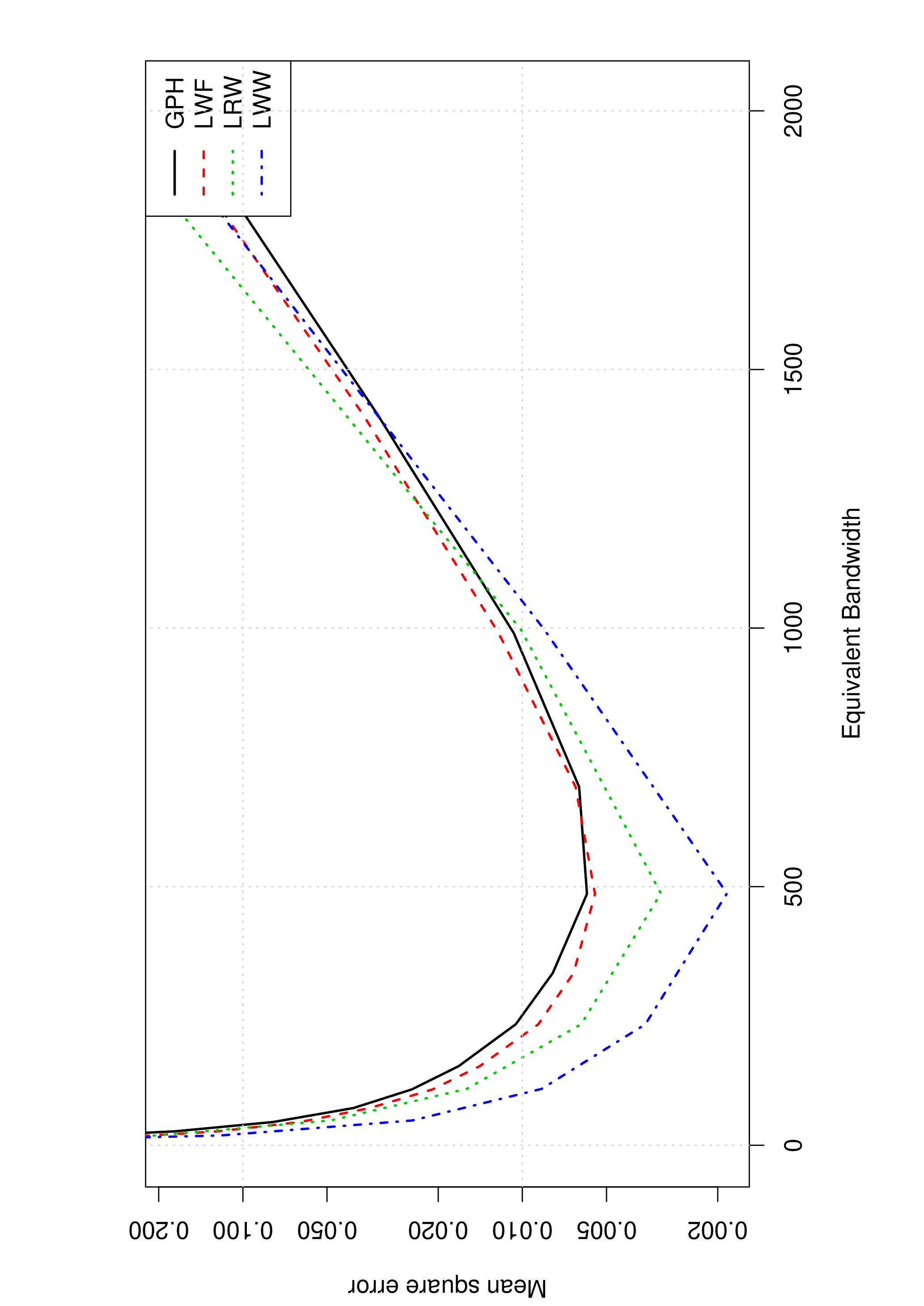}  \includegraphics[width=0.8\linewidth]{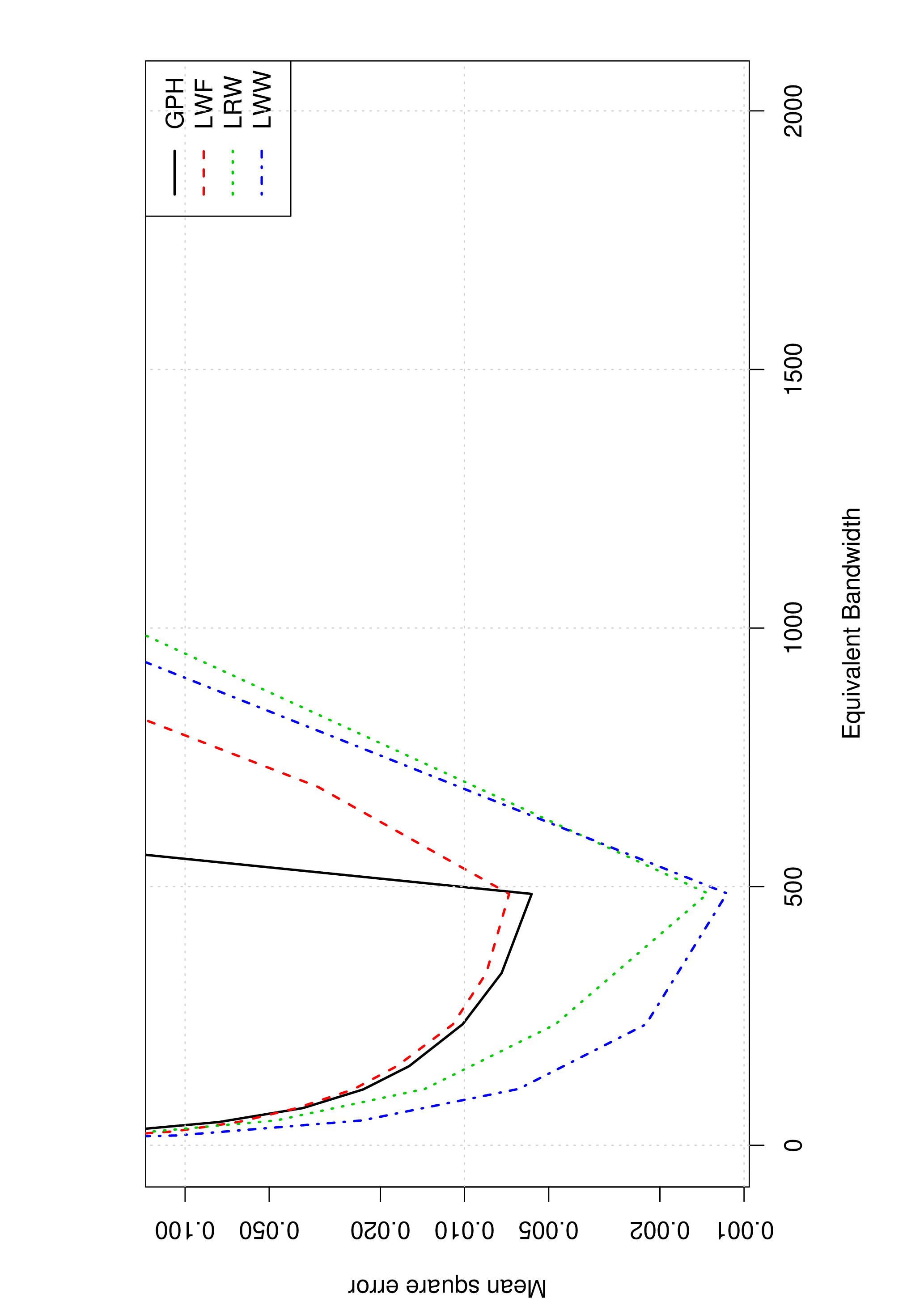}
  \caption{MSE comparison; Top panel: ARFIMA(1,$d$,0) model, $d=1.5$. Bottom
    panel: DARFIMA(1,$d$,0), $d=0.3$, $\lambda_0 = \pi/2$}
  \label{fig:MSE}
\end{figure}

The root mean square error is shown in some particular cases in
Figure~\ref{fig:MSE}. The MSE of the Fourier criteria is plotted against the
value of the bandwidth, that is, $m_n$ for the \LWF\ estimator and $m_n\times(\pooling+\taper)$ for the
\GPH\ estimator. For the wavelet methods, the somehow arbitrary ``equivalent bandwidth''
abscissa is half the number of wavelet coefficients used by the estimators: $m_n^{\text{equiv}} = \tfrac 1 2 \sum_{j=L_n}^{U_n}n_j$ (see
the discussion on the comparison of the asymptotic variances in the previous
section).


\section{Software}
\label{sec:software}

\def\TOOLBOX{../ToolboxLRD/}

\fvset{frame=none,xleftmargin=1cm,xrightmargin=1cm,baselinestretch=1,fontsize=\footnotesize}

The software used to perform the estimation of the long-memory parameter was
written in \textsc{Matlab/Octave} and may be obtained from the authors.
It includes the four estimators (\LWF, \GPH, \LWW and \LRW), basic random processes
generators and some other utilities such as the pyramidal algorithm for
computing wavelet coefficient.

\indent\textbf{Basic installation}. After downloading the tar archive
  (\texttt{toolboxLRD.tar}) and expanding it in \eg \texttt{/home/user/octave},
  one has to add the directory to the search path:
\begin{Verbatim}
  addpath(genpath('/home/user/octave/ToolboxLRD'));
\end{Verbatim}
This line can be added to the \texttt{.octave} or \texttt{.matlab} file.  Some
demos are available in \texttt{ToolboxLRD/Examples}.

\indent\textbf{Loading the data}. Use the \texttt{load} command to load a data set
(time series) into some vector, say \texttt{x}. One can also synthesize some
trajectories using one's personal generator or the one present in the \texttt{Utils}
subdirectory. For instance, a 4096 long trajectory of the ARFIMA model $(1-B)^d
(X_t - \alpha X_{t-1}) =Z_t$ with $d=1.4$, $\alpha=0.8$ and Gaussian i.i.d
sequence $Z$ can be obtained by setting:
\begin{Verbatim}
n = 4096;
x = randARFIMA(1.4,[0.8],[],n);
\end{Verbatim}
The argument \texttt{[]} above means that the MA part of the generated ARFIMA is a weak white noise (MA($0$)).  
This example is used in the following to describe the package.

\indent\textbf{Estimating the long-memory parameter}.  We shall now obtain the
\LRW, \LWW, \GPH\ and \LWF\ estimators of the memory parameter $d$ of the
series as well as an estimate of their standard deviation using the asymptotic
variance given in Theorems \ref{thm:LRW:CLT}, \ref{thm:LWW:CLT},
\ref{theo:gphgaus}, and \ref{theo:LWEestimator}. The standard deviations can be
used to build asymptotic confidence intervals.  

\begin{enumerate}[1.]
\item 
The \textit{Geweke-Porter Hudak} (\GPH) estimator is obtained as follows:
\begin{Verbatim}
param.taper=5; param.pooling= 4;
param.bandwidth=[6 12 24 50]; param.difforder = 4;
[d, stds]=GPH(x,param)
\end{Verbatim}
where \texttt{param.bandwidth} is a vector giving the different values for the
upper Fourier frequency $\bandwidth{}$ on which the regression is to be
performed (Theorem~\ref{theo:gphgaus}). The taper order $\taper$, pooling order
$\pooling$ and differentiation order $\difforder$ are specified by
\texttt{param.taper}, \texttt{param.pooling} and \texttt{param.difforder}, respectively.  One obtains
\begin{Verbatim}
d =
    1.2744    1.2736    1.3670    1.3930
stds =
    0.1803    0.1215    0.0840    0.0576
\end{Verbatim}
Note that \texttt{d} and \texttt{stds} are vectors. 
In the above example, they have four components, corresponding respectively to
the bandwidths $\bandwidth{} = 6, 12, 24, 50$.
\item 
The \textit{Local Whittle Fourier} (\LWF) estimator is invoked in the following way:
\begin{Verbatim}
param.taper=5; param.difforder=4;
param.bandwidth=[50 100 200 500]; 
[d, stds]=LWF(x,param,[])
\end{Verbatim}
One gets :
\begin{Verbatim}
d =
    1.3114    1.2974    1.2905    1.3422
stds =
    0.1206    0.0853    0.0603    0.0381
\end{Verbatim}
Here the minimization of~(\ref{eq:profile:likelihood}) is done over the whole real line. As for the \texttt{LWW} function, one may specify the range where to
optimize the \LWF\ contrast  function~(\ref{eq:profile:likelihood}) by replacing the third argument \texttt{[]} by an interval $[\Delta_1,\Delta_2]$.
For instance, if one wants to restrict this minimization to the set~(\ref{eq:allowedRangeFourier}) of admissible
values of $d$ for the CLT Theorem~\ref{theo:LWEestimator} to hold, 
\begin{Verbatim}
range= [ param.difforder-param.taper-0.5,  param.difforder+0.5];
[d, stds]=LWF(x,param,range)
\end{Verbatim}
In this specific case, the output is unchanged since the minimizing values of $d$ are within the corresponding interval $[-1.5,4.5]$.
\item 
The \textit{Local Whittle Wavelets} (\LWW) estimator is obtained as follows:
\begin{Verbatim}
LU = [6 9; 5 9; 4 9; 3 9]; 
[phi, M, alpha] = scalingfilter('Daubechies',4);
[d, stds] = LWW(x,LU,phi,[])
\end{Verbatim}
where  \texttt{phi} indicates the scaling function, \texttt{M} the corresponding number of vanishing moments, \texttt{alpha}
the Fourier decay exponent (see \allWA),  \texttt{x} contains a finite set of observations  and 
\texttt{LU} is a two column matrix, containing scales limits $\lowscale$ and $\upscale$
in the \LWW\ objective function~(\ref{eq:TildeJdef}). If  \texttt{LU} is a one column vector then it contains different
values of the lower scale $\lowscale$ and $\upscale$ is taken equal to the maximal available scale index $\maxscale$ defined
in~(\ref{eq:J0J1n}). The argument \texttt{[]} above means 
that  the interval, denoted $[\Delta_1,\Delta_2]$ in
Definition~(\ref{eq:LWWdef}), over which the contrast function~(\ref{eq:TildeJdef}) is minimized is $(-\infty,\infty)$.  
It can be replaced by an interval $[\Delta_1,\Delta_2]$, if one wants to restrict the minimization to a particular range, for
instance to the one
given by~(\ref{eq:allowedRangeWavelet}) which corresponds to admissible values of $d$ for the CLT Theorem~\ref{thm:LWW:CLT}
to hold. One gets :
\begin{Verbatim}
d =
    1.3183    1.3366    1.4209    1.3788
stds =
    0.1138    0.0719    0.0479    0.0324
\end{Verbatim}
\item 
The \textit{Local Regression Wavelet}
(\LRW) estimator and the standard deviation is invoked in the following way:
\begin{Verbatim}
L = [6;5;4;3];
[d,stds] = LRW(x,L,phi)
\end{Verbatim}
The three first argument of \texttt{LRW} are the same as \texttt{LWW} but  \texttt{LU} has been replaced by
\texttt{L}, a column vector containing different choices for the lower scale $\lowscale$ used in the regression, see
Eq.~(\ref{eq:definition:estimator:regression}). In this case, the upper scale is the
largest scale available ($\upscale = \maxscale_n$). If one wants to take different values for $\upscale$, a two-columns
matrix must replace \texttt{L}, for example, by the \texttt{LU} in the \texttt{LWW} case.
Here the \LRW\ estimator is computed using Abry--Veitch weights (see~\ref{eq:ABweigths}) and one gets the output:
\begin{Verbatim} 
d =
    1.4161    1.3815    1.4305    1.3882
stds =
    0.1020    0.0675    0.0461    0.0317
\end{Verbatim}
But the \LRW\ estimator can also be obtained using the optimal weights. In fact, the following additional output variables
are available :  
\begin{enumerate}
\item the value of a log-regression multiplicative constant $c$ so that  $\hvj{j}{\nj{j}} \approx c \; 2^{2 d
    j}$. Equivalently, $\log \hvj{j}{\nj{j}} \approx \log c +2 d j$, where $\log c$ is the intercept of the regression line;
\item new estimates of $d$ using the two-step procedure based on the optimal weights~(\ref{eq:wopt}). These weights are
  computed using the preliminary value of $d$ estimated with the Abry-Veitch weights;
\item the standard deviations of the new estimates of $d$;
\item the corresponding values of the log-regression multiplicative constant.
\end{enumerate}
Thus if the \texttt{LRW} call of the last example is
replaced by
\begin{Verbatim}
[d,stds, c, dopt,stdopt, copt] = LRW(x,L,phi)
\end{Verbatim}
the following additional output is added to the previous one :
\begin{Verbatim}
c =
    0.0074    0.0105    0.0067    0.0094
dopt =
    1.4199    1.3852    1.4451    1.3834
stdopt =
    0.1011    0.0666    0.0453    0.0311
copt =
    0.0073    0.0103    0.0060    0.0097
\end{Verbatim}

Alternatively, one could also use a different preliminary estimate of $d$, say $d=1.3823$, a value obtained as the first output of the
\texttt{GPH} routine above :
\begin{Verbatim}
[d,stds, c, dopt,stdopt, copt] = LRW(x,L,phi,1.3823)
\end{Verbatim}
The last three output values are then replaced by
\begin{Verbatim}
dopt =
    1.4199    1.3852    1.4441    1.3835
stdopt =
    0.1011    0.0666    0.0451    0.0310
copt =
    0.0073    0.0103    0.0061    0.0097
\end{Verbatim}
\end{enumerate}

\textbf{Obtaining confidence intervals.} 
A routine for obtaining the asymptotically normal  \textit{confidence intervals} for $d$ at a given level has also been
included. It works as follows :
\begin{Verbatim}
p=0.95; d=dopt; stds=stdopt;
[I]=ConfidenceInterval(d,stds,p)
\end{Verbatim}
where \texttt{d} and \texttt{stds} are any outputs of the above procedures and \texttt{p} is the confidence level. Here we
used the last displayed estimates with $p=0.95$ and get
\begin{Verbatim}
I =
    1.2217    1.6180
    1.2546    1.5158
    1.3557    1.5325
    1.3227    1.4443
\end{Verbatim}
In this particular example we can see that the true value of $d$, namely $1.4$,
belongs to the four intervals.
 
\textbf{Obtaining the theoretical asymptotic variances.} 
It is possible also to obtain directly the asymptotic variances $\AsympVarRWEopt[\psi]{d,\ell}$ and
$\AsympVarWWE[\psi]{d,\ell}$ defined in~(\ref{eq:VARwopt}) and~(\ref{eq:varLimite-k}).
These are the asymptotic variances of the \LRW\ estimator (Theorem~\ref{thm:LRW:CLT}), when, respectively, the optimal
weigths~(\ref{eq:wopt}) are chosen or when the Abry-Veitch weights~(\ref{eq:ABweigths}) are chosen. The 
 asymptotic variance of the \LWW\ estimator is also $\AsympVarWWE[\psi]{d,\ell}$ (Theorem~\ref{thm:LWW:CLT}).
The approximation  of $\AsympVarWWE[\psi]{d,\ell}$  given in~(\ref{eq:shannonApprox}) and obtained by replacing $\psi$ by the
Shannon wavelet is also available. This is how to get these asymptotic variances :
\begin{Verbatim}
d=1.4; l=5;
[v, vs, vopt, wopt]= AsymptoticVarianceLRW(phi,d,l)
\end{Verbatim}
Here $\AsympVarWWE[\psi]{d,\ell}$ and  $\AsympVarRWEopt[\psi]{d,\ell}$ are computed for $d=1.4$ and $\ell=5$ and 
stacked in the output \texttt{v} and \texttt{vopt} respectively. The output \texttt{vs} corresponds to the Shannon
approximation~(\ref{eq:shannonApprox}) and \texttt{wopt} to the optimal weights~(\ref{eq:wopt}), of length $\ell+1=6$. For these
values one gets 
\begin{Verbatim}
v =
    0.5848
vs =
    0.4949
vopt =
    0.5698
wopt =
   -0.2693    0.0546    0.0827    0.0587    0.0410    0.0322
\end{Verbatim}
We observe that for this value of $d$ the optimal variance ($0.5698$) is sensitively lower than the one obtained with
Abry-Veitch weights ($0.5848$), but still larger than the Shannon approximation ($0.4949$). Since this approximation 
gets sharper as the number of vanishing moments of the Daubechies wavelet increases, it indicates that one could get a better
variance by increasing $M$, here $M=2$. Notice, however, that the length of the wavelet filters would also increase and thus the
number of available wavelet coefficients decrease for a finite $n$, an effect which is not considered in the asymptotic
variance, see Section~\ref{sec:discussion}.


\section{Conclusion}
\label{sec:conclusion}
We have compared four semi-parametric methods for the estimation of the long-memory parameter $d$ in times series, two
Fourier-based and two wavelet-based. These are the Geweke-Porter Hudak (\GPH) [Regression/Fourier], Local Whittle Fourier
(\LRW)  [Whittle/Fourier],  Local Regression Wavelet (\LRW)  [Regression/Wavelets] and Local Whittle Wavelet (\LWW)
[Whittle/Wavelets].
We have discussed issues related to differencing, tapering and pooling in the case of Fourier-based estimators and choices of
wavelets in the case of wavelet-based estimators. Conditions for the asymptotic normality of the estimators are specified in
Theorems~\ref{thm:LRW:CLT},~\ref{thm:LWW:CLT},~\ref{theo:gphgaus} and~\ref{theo:LWEestimator}.

We have undertaken a Monte Carlo comparison. In the Monte Carlo study, we have focused on ARFIMA(0,$d$,0) and
ARFIMA(1,$d$,0) models with an AR(1) parameter equal to 0.8, a relatively high value, as well as on DARFIMA and subordinated models
defined in Section~\ref{sec:monte-carlo}. All four methods appear to work well with similar performances at the optimal
bandwidth lower scale index.\note{f}{i did not include comments about the table because the results are so model dependent !
see also figure 6 where wavelet methods outperform Fourier's.. i think perhaps we should wait for the last experiments to do
that (i am working on the optimal weights..)}
We have also developed a software package for the benefit of the practitioner which computes the corresponding estimates of
the long-memory parameter $d$ and provides confidence intervals, based on the asymptotically normal distribution of the
estimators. 

We noted that the \LRW\ estimator with Abry-Veitch weights~(\ref{eq:ABweigths}) has the same asymptotic variance as the \LWW\
estimator. This means that the \LRW\ estimator, when used with the optimal weights~(\ref{eq:wopt}), has smaller asymptotic
variance than the \LWW\ estimator.


\bibliographystyle{ims}
{\small
\bibliography{\BIBDIR/lrd}

\begin{thebibliography}{59}
\expandafter\ifx\csname natexlab\endcsname\relax\def\natexlab#1{#1}\fi
\expandafter\ifx\csname url\endcsname\relax
  \def\url#1{\texttt{#1}}\fi
\expandafter\ifx\csname urlprefix\endcsname\relax\def\urlprefix{URL }\fi

\bibitem[{Abadir et~al.(2007)Abadir, Distaso and
  Giraitis}]{abadir:distaso:giraitis:2007}
\textsc{Abadir, K.}, \textsc{Distaso, W.} and \textsc{Giraitis, L.} (2007).
\newblock Nonstationarity-extended local whittle estimation.
\newblock \textit{J. of Econometrics} \textbf{141} 1353--1384.

\bibitem[{Abry et~al.(2000)Abry, Flandrin, Taqqu and
  Veitch}]{abry:flandrin:taqqu:veitch:2000}
\textsc{Abry, P.}, \textsc{Flandrin, P.}, \textsc{Taqqu, M.~S.} and
  \textsc{Veitch, D.} (2000).
\newblock Wavelets for the analysis, estimation and synthesis of scaling data.
\newblock In \textit{Self-Similar Network Traffic and Performance Evaluation}
  (K.~Park and W.~Willinger, eds.). Wiley (Interscience Division), New York.

\bibitem[{Abry et~al.(2003)Abry, Flandrin, Taqqu and
  Veitch}]{abry:flandrin:taqqu:veitch:2003-livre}
\textsc{Abry, P.}, \textsc{Flandrin, P.}, \textsc{Taqqu, M.~S.} and
  \textsc{Veitch, D.} (2003).
\newblock Self-similarity and long-range dependence through the wavelet lens.
\newblock In \textit{{\it Theory and Applications of Long-range Dependence}}
  (P.~Doukhan, G.~Oppenheim and M.~S. Taqqu, eds.). Birkh{\"a}user, 527--556.

\bibitem[{Abry and Veitch(1998)}]{abry:veitch:1998}
\textsc{Abry, P.} and \textsc{Veitch, D.} (1998).
\newblock Wavelet analysis of long-range-dependent traffic.
\newblock \textit{IEEE Trans. Inform. Theory} \textbf{44} 2--15.

\bibitem[{Andrews and Guggenberger(2003)}]{andrews:guggenberger:2003}
\textsc{Andrews, D. W.~K.} and \textsc{Guggenberger, P.} (2003).
\newblock A bias-reduced log-periodogram regression estimator for the
  long-memory parameter.
\newblock \textit{Econometrica} \textbf{71} 675--712.

\bibitem[{Andrews and Sun(2004)}]{andrews:sun:2004}
\textsc{Andrews, D. W.~K.} and \textsc{Sun, Y.} (2004).
\newblock Adaptive local polynomial {W}hittle estimation of long-range
  dependence.
\newblock \textit{Econometrica} \textbf{72} 569--614.

\bibitem[{Bardet(2000)}]{bardet:2000T}
\textsc{Bardet, J.-M.} (2000).
\newblock Testing for the presence of self-similarity of {G}aussian time series
  having stationary increments.
\newblock \textit{Journal of Time Series Analysis} \textbf{21} 497--515.

\bibitem[{Bardet(2002)}]{bardet:2002}
\textsc{Bardet, J.-M.} (2002).
\newblock Statistical study of the wavelet analysis of fractional {B}rownian
  motion.
\newblock \textit{IEEE Trans. Inform. Theory} \textbf{48} 991--999.

\bibitem[{Bardet et~al.(2000)Bardet, Lang, Moulines and
  Soulier}]{bardet:lang:moulines:soulier:2000}
\textsc{Bardet, J.-M.}, \textsc{Lang, G.}, \textsc{Moulines, E.} and
  \textsc{Soulier, P.} (2000).
\newblock Wavelet estimator of long-range dependent processes.
\newblock \textit{Stat. Inference Stoch. Process.} \textbf{3} 85--99.
\newblock 19th ``Rencontres Franco-Belges de Statisticiens'' (Marseille, 1998).

\bibitem[{Cohen(2003)}]{cohen:2003}
\textsc{Cohen, A.} (2003).
\newblock \textit{Numerical analysis of wavelet methods}, vol.~32 of
  \textit{Studies in Mathematics and its Applications}.
\newblock North-Holland Publishing Co., Amsterdam.

\bibitem[{Dalla et~al.(2006)Dalla, Giraitis and
  Hidalgo}]{dalla:giraitis:hidalgo:2006}
\textsc{Dalla, V.}, \textsc{Giraitis, L.} and \textsc{Hidalgo, J.} (2006).
\newblock Consistent estimation of the memory parameter for nonlinear time
  series.
\newblock \textit{J. Time Ser. Anal.} \textbf{27} 211--251.

\bibitem[{Daubechies(1992)}]{daubechies:1992}
\textsc{Daubechies, I.} (1992).
\newblock \textit{Ten lectures on wavelets}, vol.~61 of \textit{CBMS-NSF
  Regional Conference Series in Applied Mathematics}.
\newblock Society for Industrial and Applied Mathematics (SIAM), Philadelphia,
  PA.

\bibitem[{Deo et~al.(2006{\natexlab{a}})Deo, Hsieh, Hurvich and
  Soulier}]{deo:hsieh:hurvich:soulier:2006}
\textsc{Deo, R.}, \textsc{Hsieh, M.}, \textsc{Hurvich, C.~M.} and
  \textsc{Soulier, P.} (2006{\natexlab{a}}).
\newblock Long memory in nonlinear processes.
\newblock In \textit{Dependence in probability and statistics}, vol. 187 of
  \textit{Lecture Notes in Statist.} Springer, New York, 221--244.

\bibitem[{Deo et~al.(2006{\natexlab{b}})Deo, Hurvich and
  Lu}]{deo:hurvich:lu:2006}
\textsc{Deo, R.}, \textsc{Hurvich, C.~M.} and \textsc{Lu, Y.}
  (2006{\natexlab{b}}).
\newblock Forecasting realized volatility using a long-memory stochastic
  volatility model: estimation, prediction and seasonal adjustment.
\newblock \textit{J. Econometrics} \textbf{131} 29--58.

\bibitem[{Fa{\"y} et~al.(2004)Fa{\"y}, Moulines and
  Soulier}]{fay:moulines:soulier:2004}
\textsc{Fa{\"y}, G.}, \textsc{Moulines, E.} and \textsc{Soulier, P.} (2004).
\newblock Edgeworth expansions for linear statistics of possibly
  long-range-dependent linear processes.
\newblock \textit{Statist. Probab. Lett.} \textbf{66} 275--288.

\bibitem[{Fa{\"y} et~al.(2007)Fa{\"y}, Roueff and
  Soulier}]{fay:roueff:soulier:2007}
\textsc{Fa{\"y}, G.}, \textsc{Roueff, F.} and \textsc{Soulier, P.} (2007).
\newblock Estimation of the memory parameter of the infinite-source {P}oisson
  process.
\newblock \textit{Bernoulli} \textbf{13} 473--491.

\bibitem[{Fox and Taqqu(1986)}]{fox:taqqu:1986}
\textsc{Fox, R.} and \textsc{Taqqu, M.~S.} (1986).
\newblock Large-sample properties of parameter estimates for strongly dependent
  stationary {G}aussian time series.
\newblock \textit{Ann. Statist.} \textbf{14} 517--532.

\bibitem[{Geweke and Porter-Hudak(1983)}]{geweke:porter-hudak:1983}
\textsc{Geweke, J.} and \textsc{Porter-Hudak, S.} (1983).
\newblock The estimation and application of long memory time series models.
\newblock \textit{J. Time Ser. Anal.} \textbf{4} 221--238.

\bibitem[{Giraitis et~al.(1997)Giraitis, Robinson and
  Samarov}]{giraitis:robinson:samarov:1997}
\textsc{Giraitis, L.}, \textsc{Robinson, P.~M.} and \textsc{Samarov, A.}
  (1997).
\newblock Rate optimal semiparametric estimation of the memory parameter of the
  {G}aussian time series with long range dependence.
\newblock \textit{J. Time Ser. Anal.} \textbf{18} 49--61.

\bibitem[{Granger and Joyeux(1980)}]{granger:joyeux:1980}
\textsc{Granger, C.~W.~J.} and \textsc{Joyeux, R.} (1980).
\newblock An introduction to long-memory time series models and fractional
  differencing.
\newblock \textit{J. Time Ser. Anal.} \textbf{1} 15--29.

\bibitem[{Hannan and Nicholls(1977)}]{hannan:nicholls:1977}
\textsc{Hannan, E.~J.} and \textsc{Nicholls, D.~F.} (1977).
\newblock The estimation of the prediction error variance.
\newblock \textit{J. Amer. Statist. Assoc.} \textbf{72} 834--840.

\bibitem[{Hurvich and Chen(2000)}]{hurvich:chen:2000}
\textsc{Hurvich, C.~M.} and \textsc{Chen, W.~W.} (2000).
\newblock An efficient taper for potentially overdifferenced long-memory time
  series.
\newblock \textit{J. Time Ser. Anal.} \textbf{21} 155--180.

\bibitem[{Hurvich et~al.(2005{\natexlab{a}})Hurvich, Lang and
  Soulier}]{hurvich:lang:soulier:2005}
\textsc{Hurvich, C.~M.}, \textsc{Lang, G.} and \textsc{Soulier, P.}
  (2005{\natexlab{a}}).
\newblock Estimation of long memory in the presence of a smooth nonparametric
  trend.
\newblock \textit{J. Amer. Statist. Assoc.} \textbf{100} 853--871.

\bibitem[{Hurvich et~al.(2002)Hurvich, Moulines and
  Soulier}]{hurvich:moulines:soulier:2002}
\textsc{Hurvich, C.~M.}, \textsc{Moulines, E.} and \textsc{Soulier, P.} (2002).
\newblock The {FEXP} estimator for potentially non-stationary linear time
  series.
\newblock \textit{Stoch. Proc. App.} \textbf{97} 307--340.

\bibitem[{Hurvich et~al.(2005{\natexlab{b}})Hurvich, Moulines and
  Soulier}]{hurvich:moulines:soulier:2005}
\textsc{Hurvich, C.~M.}, \textsc{Moulines, E.} and \textsc{Soulier, P.}
  (2005{\natexlab{b}}).
\newblock Estimating long memory in volatility.
\newblock \textit{Econometrica} \textbf{73} 1283--1328.

\bibitem[{Hurvich and Ray(1995)}]{hurvich:ray:1995}
\textsc{Hurvich, C.~M.} and \textsc{Ray, B.~K.} (1995).
\newblock Estimation of the memory parameter for nonstationary or noninvertible
  fractionally integrated processes.
\newblock \textit{J. Time Ser. Anal.} \textbf{16} 17--41.

\bibitem[{Johnson and Kotz(1970)}]{johnson:kotz:1970}
\textsc{Johnson, N.~L.} and \textsc{Kotz, S.} (1970).
\newblock \textit{Distributions in statistics. {C}ontinuous univariate
  distributions. 2.}
\newblock Houghton Mifflin Co., Boston, Mass.

\bibitem[{Kaplan and Kuo(1993)}]{kaplan:kuo:1993}
\textsc{Kaplan, L.~M.} and \textsc{Kuo, C.-C.~J.} (1993).
\newblock Fractal estimation from noisy data via discrete fractional {G}aussian
  noise ({DFGN}) and the {H}aar basis.
\newblock \textit{IEEE Trans. Signal Process.} \textbf{41} 3554--3562.

\bibitem[{K{\"u}nsch(1987)}]{kunsch:1987}
\textsc{K{\"u}nsch, H.} (1987).
\newblock Statistical aspects of self-similar processes.
\newblock In \textit{Proceedings of the 1st World Congress of the Bernoulli
  Society, Vol. 1 (Tashkent, 1986)}. VNU Sci. Press, Utrecht.

\bibitem[{Lahiri(2003)}]{lahiri:2003}
\textsc{Lahiri, S.~N.} (2003).
\newblock A necessary and sufficient condition for asymptotic independence of
  discrete {F}ourier transforms under short- and long-range dependence.
\newblock \textit{Ann. Statist.} \textbf{31} 613--641.

\bibitem[{Mallat(1998)}]{mallat:1998}
\textsc{Mallat, S.} (1998).
\newblock \textit{A wavelet tour of signal processing}.
\newblock Academic Press Inc., San Diego, CA.

\bibitem[{McCoy and Walden(1996)}]{mccoy:walden:1996}
\textsc{McCoy, E.~J.} and \textsc{Walden, A.~T.} (1996).
\newblock Wavelet analysis and synthesis of stationary long-memory processes.
\newblock \textit{J. Comput. Graph. Statist.} \textbf{5} 26--56.

\bibitem[{Moulines et~al.(2007{\natexlab{a}})Moulines, Roueff and
  Taqqu}]{moulines:roueff:taqqu:2007:fractals}
\textsc{Moulines, E.}, \textsc{Roueff, F.} and \textsc{Taqqu, M.~S.}
  (2007{\natexlab{a}}).
\newblock Central {L}imit {T}heorem for the log-regression wavelet estimation
  of the memory parameter in the {G}aussian semi-parametric context.
\newblock To appear.

\bibitem[{Moulines et~al.(2007{\natexlab{b}})Moulines, Roueff and
  Taqqu}]{moulines:roueff:taqqu:2007:jtsa}
\textsc{Moulines, E.}, \textsc{Roueff, F.} and \textsc{Taqqu, M.~S.}
  (2007{\natexlab{b}}).
\newblock On the spectral density of the wavelet coefficients of long memory
  time series with application to the log-regression estimation of the memory
  parameter.
\newblock \textit{J. Time Ser. Anal.} \textbf{28}.
\newline\urlprefix\url{http://arxiv.org/abs/math.ST/0512635}

\bibitem[{Moulines et~al.(2007{\natexlab{c}})Moulines, Roueff and
  Taqqu}]{moulines:roueff:taqqu:2007:aos}
\textsc{Moulines, E.}, \textsc{Roueff, F.} and \textsc{Taqqu, M.~S.}
  (2007{\natexlab{c}}).
\newblock A wavelet {W}hittle estimator of the memory parameter of a
  non-stationary {G}aussian time series.
\newblock Tech. rep., Ecole Nationale Supérieure des Télécommunications et
  Boston University.
\newblock To appear in the Annals of Statistics.
\newline\urlprefix\url{http://arxiv.org/abs/math.ST/0601070}

\bibitem[{Moulines and Soulier(2003)}]{moulines:soulier:2003}
\textsc{Moulines, E.} and \textsc{Soulier, P.} (2003).
\newblock Semiparametric spectral estimation for fractional processes.
\newblock In \textit{Theory and applications of long-range dependence}.
  Birkh\"auser Boston, Boston, MA, 251--301.

\bibitem[{Percival et~al.(2000)Percival, Sardy and
  Davison}]{percival:sardy:davison:2000}
\textsc{Percival, D.~B.}, \textsc{Sardy, S.} and \textsc{Davison, A.~C.}
  (2000).
\newblock Wavestrapping time series: adaptive wavelet-based bootstrapping.
\newblock In \textit{Nonlinear and nonstationary signal processing (Cambridge,
  1998)}. Cambridge Univ. Press, Cambridge, 442--471.

\bibitem[{Percival and Walden(2006)}]{percival:walden:2006}
\textsc{Percival, D.~B.} and \textsc{Walden, A.~T.} (2006).
\newblock \textit{Wavelet methods for time series analysis}, vol.~4 of
  \textit{Cambridge Series in Statistical and Probabilistic Mathematics}.
\newblock Cambridge University Press, Cambridge.
\newblock Reprint of the 2000 original [MR1770693].

\bibitem[{Robinson(1994)}]{robinson:1994:NS}
\textsc{Robinson, P.~M.} (1994).
\newblock Efficient tests of nonstationary hypotheses.
\newblock \textit{J. Amer. Statist. Assoc.} \textbf{89} 1420--1437.

\bibitem[{Robinson(1995{\natexlab{a}})}]{robinson:1995:GSE}
\textsc{Robinson, P.~M.} (1995{\natexlab{a}}).
\newblock {G}aussian semiparametric estimation of long range dependence.
\newblock \textit{Ann. Statist.} \textbf{23} 1630--1661.

\bibitem[{Robinson(1995{\natexlab{b}})}]{robinson:1995:GPH}
\textsc{Robinson, P.~M.} (1995{\natexlab{b}}).
\newblock Log-periodogram regression of time series with long range dependence.
\newblock \textit{The Annals of Statistics} \textbf{23} 1048--1072.

\bibitem[{Robinson and Henry(2003)}]{robinson:henry:2003}
\textsc{Robinson, P.~M.} and \textsc{Henry, M.} (2003).
\newblock Higher-order kernel semiparametric {$M$}-estimation of long memory.
\newblock \textit{J. Econometrics} \textbf{114} 1--27.

\bibitem[{Roueff and Taqqu(2007)}]{roueff:taqqu:2007}
\textsc{Roueff, F.} and \textsc{Taqqu, M.~S.} (2007).
\newblock Asymptotic normality of wavelet estimators of the memory parameter:
  the linear case.
\newblock Tech. rep.

\bibitem[{Samorodnitsky and Taqqu(1994)}]{samorodnitsky:taqqu:1994}
\textsc{Samorodnitsky, G.} and \textsc{Taqqu, M.~S.} (1994).
\newblock \textit{Stable non-{G}aussian processes: stochastic models with
  infinite variance}.
\newblock Chapman and Hall.

\bibitem[{Serroukh et~al.(2000)Serroukh, Walden and
  Percival}]{serroukh:walden:percival:2000}
\textsc{Serroukh, A.}, \textsc{Walden, A.~T.} and \textsc{Percival, D.~B.}
  (2000).
\newblock Statistical properties and uses of the wavelet variance estimator for
  the scale analysis of time series.
\newblock \textit{J. Amer. Statist. Assoc.} \textbf{95} 184--196.

\bibitem[{Shimotsu and Phillips(2005)}]{shimotsu:phillips:2005}
\textsc{Shimotsu, K.} and \textsc{Phillips, P.~C.~B.} (2005).
\newblock Exact local {W}hittle estimation of fractional integration.
\newblock \textit{Ann. Statist.} \textbf{33} 1890--1933.

\bibitem[{Shimotsu and Phillips(2006)}]{shimotsu:phillips:2006}
\textsc{Shimotsu, K.} and \textsc{Phillips, P.~C.~B.} (2006).
\newblock Local {W}hittle estimation of fractional integration and some of its
  variants.
\newblock \textit{J. Econometrics} \textbf{130} 209--233.

\bibitem[{Stoev et~al.(2006)Stoev, Taqqu, Park, Michailidis and
  Marron}]{stoev:taqqu:park:2006}
\textsc{Stoev, S.}, \textsc{Taqqu, M.~S.}, \textsc{Park, C.},
  \textsc{Michailidis, G.} and \textsc{Marron, J.~S.} (2006).
\newblock L{ASS}: a tool for the local analysis of self-similarity.
\newblock \textit{Comput. Statist. Data Anal.} \textbf{50} 2447--2471.

\bibitem[{Tanaka(1999)}]{tanaka:1999}
\textsc{Tanaka, K.} (1999).
\newblock The nonstationary fractional unit root.
\newblock \textit{Econometric Theory} \textbf{15} 549--582.

\bibitem[{Teyssi{\`e}re and Abry(2007)}]{teyssiere:abry:2007}
\textsc{Teyssi{\`e}re, G.} and \textsc{Abry, P.} (2007).
\newblock Wavelet analysis of nonlinear long-range dependent processes.
  {A}pplications to financial time series.
\newblock In \textit{Long memory in economics}. Springer, Berlin, 173--238.

\bibitem[{Veitch and Abry(1999)}]{veitch:abry:1999}
\textsc{Veitch, D.} and \textsc{Abry, P.} (1999).
\newblock A wavelet-based joint estimator of the parameters of long-range
  dependence.
\newblock \textit{IEEE Trans. Inform. Theory} \textbf{45} 878--897.

\bibitem[{Veitch et~al.(2003)Veitch, Abry and Taqqu}]{veitch:abry:taqqu:2003}
\textsc{Veitch, D.}, \textsc{Abry, P.} and \textsc{Taqqu, M.~S.} (2003).
\newblock On the automatic selection of the onset of scaling.
\newblock \textit{Fractals} \textbf{11} 377--390.

\bibitem[{Veitch et~al.(2000)Veitch, Taqqu and Abry}]{veitch:taqqu:abry:2000m}
\textsc{Veitch, D.}, \textsc{Taqqu, M.~S.} and \textsc{Abry, P.} (2000).
\newblock Meaningful {MRA} initialisation for discrete time series.
\newblock \textit{Signal Processing} \textbf{80} 1971--1983.

\bibitem[{Velasco(1999{\natexlab{a}})}]{velasco:1999}
\textsc{Velasco, C.} (1999{\natexlab{a}}).
\newblock Gaussian semiparametric estimation of non-stationary time series.
\newblock \textit{J. Time Ser. Anal.} \textbf{20} 87--127.

\bibitem[{Velasco(1999{\natexlab{b}})}]{velasco:1999:JE}
\textsc{Velasco, C.} (1999{\natexlab{b}}).
\newblock Non-stationary log-periodogram regression.
\newblock \textit{J. Econometrics} \textbf{91} 325--371.

\bibitem[{Velasco(2000)}]{velasco:2000:NG}
\textsc{Velasco, C.} (2000).
\newblock Non-{G}aussian log-periodogram regression.
\newblock \textit{Econometric Theory} \textbf{16} 44--79.

\bibitem[{Velasco and Robinson(2000)}]{velasco:robinson:2000}
\textsc{Velasco, C.} and \textsc{Robinson, P.~M.} (2000).
\newblock Whittle pseudo-maximum likelihood estimation for nonstationary time
  series.
\newblock \textit{J. Am. Statist. Assoc.} \textbf{95} 1229--1243.

\bibitem[{Wornell and Oppenheim(1992)}]{wornell:oppenheim:1992}
\textsc{Wornell, G.~W.} and \textsc{Oppenheim, A.~V.} (1992).
\newblock Estimation of fractal signals from noisy measurements using wavelets.
\newblock \textit{IEEE Trans. Signal Process.} \textbf{40} 611 -- 623.

\bibitem[{{\v{Z}}urbenko(1979)}]{zurbenko:1979}
\textsc{{\v{Z}}urbenko, I.} (1979).
\newblock On the efficiency of estimates of a spectral density.
\newblock \textit{Scand. J. Statist.} \textbf{6} 49--56.

\end{thebibliography}
}


\begin{sidewaystable}[p]
  \centering
\begin{scriptsize}
 
 \begin{tabular}{c||cccc|cccc|cccc|cccc|} 
      &\multicolumn{4}{|c}{\textbf{GPH}}& \multicolumn{4}{|c}{\textbf{LWF}}& \multicolumn{4}{|c|}{\textbf{LRW}} & \multicolumn{4}{|c|}{\textbf{LWW}}\\ 
\textbf{Model} & bias & std & RMSE & $m_n^{\text{opt}}$ & bias &       std & RMSE & $m_n^{\text{opt}}$ & bias & std & RMSE &       $L_{n}^{\text{opt}}$ &bias & std & RMSE & $L_{n}^{\text{opt}}$       \\ 
\hline ARFIMA(0,-1.2,0) & 0.007 & 0.105 & \textbf{0.105} & 26& -0.108 & 0.071 & 0.129 & 234& 0.047 & 0.106 & 0.116 & 2& 0.105 & 0.083 & 0.134 & 2\\ 
\hline ARFIMA(1,-1.2,0) & 0.000 & 0.161 & 0.161 & 12& -0.138 & 0.128 & 0.188 & 72& -0.093 & 0.108 & 0.142 & 2& -0.048 & 0.083 & \textbf{0.096} & 2\\ 
\hline ARFIMA(0,0.0,0) & -0.022 & 0.103 & 0.105 & 26& -0.099 & 0.073 & 0.123 & 234& -0.026 & 0.058 & 0.064 & 1& -0.002 & 0.046 & \textbf{0.046} & 1\\ 
\hline ARFIMA(1,0.0,0) & -0.073 & 0.154 & 0.170 & 12& -0.175 & 0.134 & 0.220 & 72& -0.169 & 0.104 & 0.198 & 2& -0.058 & 0.143 & \textbf{0.154} & 3\\ 
\hline ARFIMA(0,0.3,0) & -0.029 & 0.104 & 0.108 & 26& -0.094 & 0.071 & 0.118 & 234& -0.065 & 0.060 & 0.088 & 1& -0.045 & 0.046 & \textbf{0.065} & 1\\ 
\hline ARFIMA(1,0.3,0) & -0.076 & 0.151 & 0.169 & 12& -0.194 & 0.105 & 0.221 & 108& -0.172 & 0.100 & 0.199 & 2& -0.060 & 0.143 & \textbf{0.154} & 3\\ 
\hline ARFIMA(0,1.5,0) & -0.049 & 0.097 & 0.109 & 26& -0.077 & 0.072 & 0.105 & 234& -0.085 & 0.110 & 0.139 & 2& -0.045 & 0.093 & \textbf{0.103} & 2\\ 
\hline ARFIMA(1,1.5,0) & -0.121 & 0.148 & 0.190 & 12& -0.182 & 0.106 & 0.210 & 108& -0.167 & 0.115 & 0.203 & 2& -0.135 & 0.091 & \textbf{0.163} & 2\\ 
\hline ARFIMA(0,2.5,0) & -0.039 & 0.094 & 0.102 & 26& -0.050 & 0.072 & \textbf{0.087} & 234& -0.093 & 0.120 & 0.152 & 2& -0.047 & 0.097 & 0.108 & 2\\ 
\hline ARFIMA(1,2.5,0) & -0.132 & 0.141 & 0.194 & 12& -0.157 & 0.108 & 0.190 & 108& -0.136 & 0.115 & 0.178 & 2& -0.101 & 0.098 & \textbf{0.141} & 2\\ 
\hline ARFIMA(0,3.5,0) & -0.023 & 0.092 & 0.095 & 26& -0.023 & 0.070 & \textbf{0.074} & 234& -0.077 & 0.110 & 0.134 & 2& -0.037 & 0.089 & 0.097 & 2\\ 
\hline ARFIMA(1,3.5,0) & -0.097 & 0.136 & 0.167 & 12& -0.111 & 0.109 & 0.155 & 108& -0.102 & 0.113 & 0.152 & 2& -0.063 & 0.097 & \textbf{0.116} & 2\\ 
\hline DARFIMA(0,0.0,0) & 0.064 & 0.275 & 0.282 & 5& 0.002 & 0.162 & 0.162 & 72& -0.013 & 0.188 & 0.188 & 3& 0.072 & 0.139 & \textbf{0.157} & 3\\ 
\hline DARFIMA(1,0.0,0) & 0.031 & 0.287 & 0.288 & 5& -0.027 & 0.155 & 0.157 & 72& -0.038 & 0.187 & 0.191 & 3& 0.041 & 0.144 & \textbf{0.149} & 3\\ 
\hline DARFIMA(0,0.3,0) & 0.036 & 0.282 & 0.284 & 5& 0.005 & 0.161 & 0.161 & 72& -0.038 & 0.190 & 0.193 & 3& 0.046 & 0.145 & \textbf{0.152} & 3\\ 
\hline DARFIMA(1,0.3,0) & 0.015 & 0.266 & 0.266 & 5& -0.030 & 0.154 & 0.157 & 72& -0.064 & 0.182 & 0.193 & 3& 0.016 & 0.146 & \textbf{0.147} & 3\\ 
\hline SUBORD1(0,0.0,0) & -0.018 & 0.099 & 0.101 & 26& -0.095 & 0.066 & 0.116 & 234& -0.032 & 0.069 & 0.076 & 1& -0.003 & 0.057 & \textbf{0.057} & 1\\ 
\hline SUBORD1(1,0.0,0) & 0.121 & 0.114 & 0.166 & 26& 0.006 & 0.075 & \textbf{0.076} & 234& -0.043 & 0.097 & 0.106 & 1& -0.015 & 0.086 & 0.087 & 1\\ 
\hline SUBORD1(0,0.3,0) & -0.111 & 0.113 & 0.159 & 26& -0.179 & 0.087 & 0.199 & 234& -0.155 & 0.090 & 0.179 & 1& -0.127 & 0.083 & \textbf{0.152} & 1\\ 
\hline SUBORD1(1,0.3,0) & -0.063 & 0.161 & \textbf{0.173} & 17& -0.133 & 0.157 & 0.206 & 153& -0.122 & 0.174 & 0.212 & 2& -0.032 & 0.185 & 0.188 & 2\\ 
\hline SUBORD2(0,0.0,0) & -0.014 & 0.103 & 0.104 & 26& -0.095 & 0.067 & 0.117 & 234& -0.026 & 0.061 & 0.066 & 1& -0.001 & 0.048 & \textbf{0.048} & 1\\ 
\hline SUBORD2(1,0.0,0) & 0.207 & 0.260 & 0.332 & 5& 0.022 & 0.183 & \textbf{0.184} & 45& 0.160 & 0.212 & 0.266 & 3& 0.293 & 0.197 & 0.353 & 3\\ 
\hline SUBORD2(0,0.3,0) & 0.014 & 0.115 & 0.116 & 26& -0.058 & 0.093 & 0.109 & 234& -0.005 & 0.076 & 0.077 & 1& 0.025 & 0.067 & \textbf{0.072} & 1\\ 
\hline SUBORD2(1,0.3,0) & 0.157 & 0.209 & 0.262 & 8& 0.036 & 0.168 & \textbf{0.172} & 72& 0.174 & 0.088 & 0.195 & 1& 0.214 & 0.085 & 0.230 & 1\\ 
   \end{tabular} 

\end{scriptsize}
\caption{Length of the time series $=512$. Bias, standard deviation, root
  mean-square error and optimal bandwidth/minimal scale for the four estimators
  applied to the three models of Section~\ref{sec:monte-carlo}. The lowest RMSE
  among the four methods appears in boldface.}
  \label{tab:results512}
\end{sidewaystable}

\thispagestyle{empty}
\begin{sidewaystable}[p]
  \centering
\begin{scriptsize}
 
 \begin{tabular}{c||cccc|cccc|cccc|cccc|} 
      &\multicolumn{4}{|c}{\textbf{GPH}}& \multicolumn{4}{|c}{\textbf{LWF}}& \multicolumn{4}{|c|}{\textbf{LRW}} & \multicolumn{4}{|c|}{\textbf{LWW}}\\ 
\textbf{Model} & bias & std & RMSE & $m_n^{\text{opt}}$ & bias &       std & RMSE & $m_n^{\text{opt}}$ & bias & std & RMSE &       $L_{n}^{\text{opt}}$ &bias & std & RMSE & $L_{n}^{\text{opt}}$       \\ 
\hline ARFIMA(0,-1.2,0) & -0.012 & 0.032 & \textbf{0.034} & 224& -0.031 & 0.022 & 0.038 & 2016& 0.020 & 0.037 & 0.043 & 3& 0.038 & 0.032 & 0.050 & 3\\ 
\hline ARFIMA(1,-1.2,0) & -0.024 & 0.061 & 0.065 & 54& -0.081 & 0.041 & 0.091 & 486& -0.046 & 0.023 & 0.051 & 2& -0.037 & 0.021 & \textbf{0.043} & 2\\ 
\hline ARFIMA(0,0.0,0) & -0.016 & 0.031 & 0.035 & 224& -0.027 & 0.022 & 0.035 & 2016& -0.006 & 0.014 & 0.015 & 1& 0.000 & 0.012 & \textbf{0.012} & 1\\ 
\hline ARFIMA(1,0.0,0) & -0.037 & 0.060 & 0.070 & 54& -0.072 & 0.043 & 0.084 & 486& -0.042 & 0.034 & 0.055 & 3& -0.031 & 0.029 & \textbf{0.043} & 3\\ 
\hline ARFIMA(0,0.3,0) & -0.017 & 0.031 & 0.036 & 224& -0.026 & 0.023 & 0.034 & 2016& -0.019 & 0.022 & 0.029 & 2& -0.010 & 0.019 & \textbf{0.021} & 2\\ 
\hline ARFIMA(1,0.3,0) & -0.044 & 0.062 & 0.076 & 54& -0.071 & 0.043 & 0.083 & 486& -0.042 & 0.036 & 0.056 & 3& -0.029 & 0.030 & \textbf{0.041} & 3\\ 
\hline ARFIMA(0,1.5,0) & -0.017 & 0.031 & 0.035 & 224& -0.021 & 0.021 & \textbf{0.029} & 2016& -0.038 & 0.026 & 0.046 & 2& -0.028 & 0.024 & 0.037 & 2\\ 
\hline ARFIMA(1,1.5,0) & -0.051 & 0.057 & 0.077 & 54& -0.062 & 0.041 & 0.074 & 486& -0.038 & 0.042 & 0.057 & 3& -0.022 & 0.037 & \textbf{0.043} & 3\\ 
\hline ARFIMA(0,2.5,0) & -0.013 & 0.030 & 0.033 & 224& -0.014 & 0.022 & \textbf{0.026} & 2016& -0.040 & 0.030 & 0.050 & 2& -0.029 & 0.027 & 0.040 & 2\\ 
\hline ARFIMA(1,2.5,0) & -0.043 & 0.058 & 0.072 & 54& -0.047 & 0.042 & 0.063 & 486& -0.035 & 0.044 & 0.056 & 3& -0.018 & 0.039 & \textbf{0.043} & 3\\ 
\hline ARFIMA(0,3.5,0) & -0.005 & 0.030 & 0.031 & 224& -0.006 & 0.022 & \textbf{0.023} & 2016& -0.033 & 0.028 & 0.043 & 2& -0.019 & 0.027 & 0.034 & 2\\ 
\hline ARFIMA(1,3.5,0) & -0.033 & 0.058 & 0.066 & 54& -0.035 & 0.044 & 0.056 & 486& -0.029 & 0.044 & 0.053 & 3& -0.014 & 0.041 & \textbf{0.043} & 3\\ 
\hline DARFIMA(0,0.0,0) & -0.023 & 0.060 & 0.064 & 54& -0.055 & 0.042 & 0.069 & 486& 0.025 & 0.035 & \textbf{0.043} & 3& -0.002 & 0.044 & 0.044 & 4\\ 
\hline DARFIMA(1,0.0,0) & -0.037 & 0.060 & 0.070 & 54& -0.072 & 0.043 & 0.084 & 486& 0.010 & 0.034 & \textbf{0.035} & 3& 0.029 & 0.029 & 0.041 & 3\\ 
\hline DARFIMA(0,0.3,0) & -0.024 & 0.063 & 0.068 & 54& -0.053 & 0.044 & 0.069 & 486& 0.012 & 0.036 & \textbf{0.038} & 3& 0.033 & 0.030 & 0.044 & 3\\ 
\hline DARFIMA(1,0.3,0) & -0.044 & 0.062 & 0.076 & 54& -0.071 & 0.043 & 0.083 & 486& -0.003 & 0.037 & 0.037 & 3& 0.016 & 0.030 & \textbf{0.034} & 3\\ 
\hline SUBORD1(0,0.0,0) & -0.015 & 0.030 & 0.033 & 224& -0.027 & 0.021 & 0.034 & 2016& -0.006 & 0.014 & 0.016 & 1& 0.000 & 0.013 & \textbf{0.013} & 1\\ 
\hline SUBORD1(1,0.0,0) & 0.028 & 0.089 & 0.093 & 26& 0.032 & 0.028 & 0.043 & 2016& 0.034 & 0.026 & \textbf{0.042} & 1& 0.044 & 0.026 & 0.051 & 1\\ 
\hline SUBORD1(0,0.3,0) & -0.103 & 0.043 & 0.112 & 224& -0.115 & 0.038 & 0.121 & 2016& -0.103 & 0.038 & 0.109 & 2& -0.082 & 0.051 & \textbf{0.096} & 3\\ 
\hline SUBORD1(1,0.3,0) & -0.091 & 0.073 & 0.117 & 110& -0.086 & 0.076 & 0.115 & 990& -0.093 & 0.057 & 0.109 & 2& -0.060 & 0.066 & \textbf{0.089} & 2\\ 
\hline SUBORD2(0,0.0,0) & -0.018 & 0.031 & 0.036 & 224& -0.028 & 0.022 & 0.036 & 2016& -0.006 & 0.014 & 0.016 & 1& -0.001 & 0.013 & \textbf{0.013} & 1\\ 
\hline SUBORD2(1,0.0,0) & 0.043 & 0.096 & 0.106 & 26& -0.008 & 0.071 & \textbf{0.072} & 234& 0.054 & 0.063 & 0.083 & 4& 0.031 & 0.090 & 0.095 & 5\\ 
\hline SUBORD2(0,0.3,0) & 0.032 & 0.045 & 0.055 & 224& 0.022 & 0.040 & 0.046 & 2016& 0.028 & 0.025 & \textbf{0.037} & 1& 0.035 & 0.024 & 0.042 & 1\\ 
\hline SUBORD2(1,0.3,0) & -0.013 & 0.079 & 0.080 & 37& 0.043 & 0.061 & 0.074 & 486& -0.045 & 0.063 & 0.078 & 4& -0.006 & 0.054 & \textbf{0.055} & 4\\ 
   \end{tabular} 

\end{scriptsize}
  \caption{Same as Table~\ref{tab:results512} with length of the time series $=4096$.}
  \label{tab:results4096}
\end{sidewaystable}

\end{document}